\documentclass[11pt]{article}
\usepackage{amsfonts}
\usepackage{bbm}
\usepackage{mathrsfs}
\usepackage{amscd}
\usepackage{color}
\usepackage{amsmath,amsfonts,amssymb,amscd}
\usepackage{indentfirst,graphics,epsfig,psfrag}
\input{epsf}
\usepackage{ifpdf}
\usepackage{enumerate}
\usepackage{appendix}
\usepackage{enumerate}

\usepackage{lineno}

\makeatletter \@addtoreset{figure}{section} \makeatother
\makeatletter
\long\def\@makecaption#1#2{%
   \vskip 10\p@
   \setbox\@tempboxa\hbox{{#1}\ \ #2}%
   \ifdim \wd\@tempboxa >\hsize

       {#1}\ \ #2\par
   \else
       \hbox to\hsize{\hfil\box\@tempboxa\hfil}%
   \fi}
\makeatother

\newtheorem{thm}{Theorem}[section]
\newtheorem{cor}{Corollary}[section]
\newtheorem{lem}{Lemma}[section]
\newtheorem{rem}{Remark}[section]
\newtheorem{con}{Conjecture}[section]

\newtheorem{obs}{Observation}[section]
\newtheorem{pro}{Proposition}[section]
\newtheorem{prop}{Problem}[section]

\begin{document}
\title{\textbf{Steiner Distance in Graphs---A Survey}\footnote{Supported by the National Science Foundation of China
(Nos. 11601254, and 11551001) and the Science Found of Qinghai
Province (Nos. 2016-ZJ-948Q, and 2014-ZJ-907).}}
\author{
\small Yaping Mao\footnote{E-mail: maoyaping@ymail.com}\\[0.2cm]
\small School of Mathematics and Statistics, Qinghai Normal\\
\small University, Xining, Qinghai 810008, China\\[0.2cm]
\small Center for Mathematics and Interdisciplinary Sciences\\
\small  of Qinghai Province, Xining, Qinghai 810008, China}
\date{}
\maketitle
\begin{abstract}
For a connected graph $G$ of order at least $2$ and $S\subseteq
V(G)$, the \emph{Steiner distance} $d_G(S)$ among the vertices of
$S$ is the minimum size among all connected subgraphs whose vertex
sets contain $S$. In this paper, we summarize the known results on
the Steiner distance parameters, including Steiner distance, Steiner diameter, Steiner center, Steiner median, Steiner interval,
Steiner distance hereditary graph, Steiner distance stable graph, average Steiner distance, and Steiner Wiener index. It also contains
some conjectures and open problems
for further studies. \\[2mm]
\end{abstract}

\section{Introduction}

All graphs in this paper are undirected, finite and simple. We refer
to \cite{Bondy} for graph theoretical notation and terminology not
described here.  For a graph $G$, let $V(G)$, $E(G)$, $e(G)$, $L(G)$
and $\overline{G}$ denote the set of vertices, the set
of edges, the size of $G$, the line graph, and the complement, respectively.
The {\it degree\/}, $deg_G(v)$, of a vertex $v$ of $G$ is
the number of edges incident with it. The {\it minimum degree\/} of $G$, $\delta(G)$,
is the smallest of the degrees of vertices in $G$ and the
{\it maximum degree\/}, $\Delta(G)$, of $G$ is the largest of the degrees of the vertices in $G$.
For any subset $X$ of $V(G)$, let $G[X]$ denote the subgraph induced
by $X$; similarly, for any subset
$F$ of $E(G)$, let $G[F]$ denote the subgraph induced by $F$. We use
$G\setminus X$ to denote the subgraph of $G$ obtained by removing all the
vertices of $X$ together with the edges incident with them from $G$;
similarly, we use $G\setminus F$ to denote the subgraph of $G$
obtained by removing all the edges of $F$ from $G$. If $X=\{v\}$ and
$F=\{e\}$, we simply write $G-v$ and $G\setminus e$ for $G-\{v\}$
and $G\setminus \{e\}$, respectively. For two subsets $X$ and $Y$ of
$V(G)$ we denote by $E_G[X,Y]$ the set of edges of $G$ with one end
in $X$ and the other end in $Y$. If $X=\{x\}$, we simply write
$E_G[x,Y]$ for $E_G[\{x\},Y]$. In this paper, we let $K_{n}$, $P_n$ and $C_n$ be the complete graph of order $n$, the path of order $n$, and the cycle of order $n$, respectively.
The \emph{closed neighborhood} of a vertex $v$ in a graph $G$ is the set
$N_G[v]=\{u\in V (G)\,|\,d(u, v)\leq 1\}$, and the \emph{open neighborhood} is $N_G(v)=\{u\in V(G)\,|\,d(u,v)=1\}$.

A graph is {\it planar\/} if it can be drawn in the plane with no crossing edges.
A graph is {\it maximal planar\/} if it is planar, but after the addition of any edge the resulting graph is not planar.
A \emph{subdivision} of $G$ is a graph obtained from $G$
by replacing edges with pairwise internally disjoint paths.
The {\it connectivity\/}, $\kappa(G)$, of $G$ is defined as the
minimum number of vertices whose deletion renders $G$ disconnected or a trivial graph.
A graph $G$ is a {\it threshold graph}, if there exists a weight
function $w: V(G)\rightarrow R$ and a real constant $t$ such that
two vertices $g,g'\in V(G)$ are adjacent if and only if
$w(g)+w(g')\geq t$. A graph is said to be \emph{minimally $k$-connected} if it is
$k$-connected but omitting any of the edges the resulting graph is
no longer $k$-connected.

Let $f(G)$ be a graph invariant and $n$ a positive integer, $n \geq
n$. The {\it Nordhaus--Gaddum Problem\/} is to determine sharp
bounds for $f(G)+f(\overline{G})$ and $f(G)\cdot f(\overline{G})$,
as $G$ ranges over the class of all graphs of order $n$, and to
characterize the extremal graphs, i.e., graphs that achieve the
bounds. Nordhaus--Gaddum type relations have received wide
attention; see the recent survey \cite{AouchicheHansen}.

The join, Cartesian product, lexicographic product, corona, and cluster are defined as
follows.
\begin{itemize}
\item[] The {\it join} or {\it complete product} of two disjoint graphs $G$ and $H$,
denoted by $G\vee H$, is the graph with vertex set $V(G)\cup V(H)$
and edge set $E(G)\cup E(H)\cup \{uv\,|\, u\in V(G), v\in V(H)\}$.

\item[] The \emph{Cartesian product} of two graphs $G$ and $H$, written as
$G\Box H$, is the graph with vertex set $V(G)\times V(H)$, in which
two vertices $(u,v)$ and $(u',v')$ are adjacent if and only if
$u=u'$ and $(v,v')\in E(H)$, or $v=v'$ and $(u,u')\in E(G)$.

\item[] The \emph{lexicographic product} of two graphs $G$ and $H$, written as $G\circ
H$, is defined as follows: $V(G\circ H)=V(G)\times V(H)$, and two
distinct vertices $(u,v)$ and $(u',v')$ of $G\circ H$ are adjacent
if and only if either $(u,u')\in E(G)$ or $u=u'$ and $(v,v')\in
E(H)$.

\item[] The {\it corona} $G*H$ is obtained by taking one copy of $G$
and $|V(G)|$ copies of $H$, and by joining each vertex of the $i$-th
copy of $H$ with the $i$-th vertex of $G$, where
$i=1,2,\ldots,|V(G)|$.

\item[] The {\it cluster} $G\odot H$ is obtained by taking one copy of $G$
and $|V(G)|$ copies of a rooted graph $H$, and by identifying the
root of the $i$-th copy of $H$ with the $i$-th vertex of $G$, where
$i=1,2,\ldots,|V(G)|$.
\end{itemize}

We divide our introduction into the following subsections to state
the motivations and our results of this paper.

\subsection{Distance parameters and their generalizations}

Distance is one of the most basic concepts of
graph-theoretic subjects. If $G$ is a connected graph and $u,v\in
V(G)$, then the \emph{distance} $d_G(u,v)$ between $u$ and $v$ is the
length of a shortest path connecting $u$ and $v$. If $v$ is a vertex
of a connected graph $G$, then the \emph{eccentricity} $e(v)$ of $v$
is defined by $e(v)=\max\{d_G(u,v)\,|\,u\in V(G)\}$. Furthermore, the
\emph{radius} $rad(G)$ and \emph{diameter} $diam(G)$ of $G$ are
defined by $rad(G)=\min\{e(v)\,|\,v\in V(G)\}$ and $diam(G)=\max
\{e(v)\,|\,v\in V(G)\}$. These last two concepts are related by the
inequalities $rad(G)\leq diam(G) \leq 2 rad(G)$. Goddard and
Oellermann gave a survey paper on this subject; see \cite{Goddard}.

\subsubsection{Steiner distance}

The distance between two vertices $u$ and $v$ in a connected graph
$G$ also equals the minimum size of a connected subgraph of $G$
containing both $u$ and $v$. This observation suggests a
generalization of distance. The Steiner distance of a graph,
introduced by Chartrand, Oellermann, Tian and Zou in 1989, is a
natural and nice generalization of the concept of classical graph
distance. For a graph $G(V,E)$ and a set $S\subseteq V(G)$ of at
least two vertices, \emph{an $S$-Steiner tree} or \emph{a Steiner
tree connecting $S$} (or simply, \emph{an $S$-tree}) is a
subgraph $T(V',E')$ of $G$ that is a tree with $S\subseteq V'$. Let
$G$ be a connected graph of order at least $2$ and let $S$ be a
nonempty set of vertices of $G$. Then the \emph{Steiner distance}
$d_G(S)$ among the vertices of $S$ (or simply the distance of $S$)
is the minimum size among all connected subgraphs whose vertex sets
contain $S$. Note that if $H$ is a connected subgraph of $G$ such
that $S\subseteq V(H)$ and $|E(H)|=d_G(S)$, then $H$ is a tree.
Observe that $d_G(S)=\min\{e(T)\,|\,S\subseteq V(T)\}$, where $T$ is
subtree of $G$. Furthermore, if $S=\{u,v\}$, then $d_G(S)=d(u,v)$ is
the classical distance between $u$ and $v$. Set
$d_G(S)=\infty$ when there is no $S$-Steiner tree in $G$.
\begin{obs}\label{obs1-1}
Let $G$ be a graph of order $n$ and $k$ be an integer with $2\leq
k\leq n$. If $S\subseteq V(G)$ and $|S|=k$, then $d_G(S)\geq k-1$.
\end{obs}

The problem of finding the Steiner distance of a set of vertices is
called the {\it Steiner Problem\/} and is $NP$-complete (see
\cite{GareyJohnson}). Steiner trees are well known for their
combinatorial optimization aspects and applications to network
design and transportation.

\subsubsection{Steiner eccentricity, Steiner diameter, and Steiner radius}

Let $n$ and $k$ be two integers with $2\leq k\leq n$. The
\emph{Steiner $k$-eccentricity $e_k(v)$} of a vertex $v$ of $G$ is
defined by $e_k(v)=\max \{d(S)\,|\,S\subseteq V(G), |S|=k,~and~v\in
S \}$. The \emph{Steiner $k$-radius} of $G$ is $srad_k(G)=\min \{
e_k(v)\,|\,v\in V(G)\}$, while the \emph{Steiner $k$-diameter} of
$G$ is $sdiam_k(G)=\max \{e_k(v)\,|\,v\in V(G)\}$. Note for every
connected graph $G$ that $e_2(v)=e(v)$ for all vertices $v$ of $G$
and that $srad_2(G)=rad(G)$ and $sdiam_2(G)=diam(G)$.
{\small
\begin{center}
\begin{tabular}{|c|c|}
\cline{1-2}
\hline  {\bf Distance} & {\bf Steiner distance} \\[0.1cm]
\cline{1-2}
$\left\{
\begin{array}{ll}
d_G(u,v)=\min\{e(P)\,|\,S\subseteq V(P)\},\\[0.1cm]
{\rm where}~P~{\rm is~
subpath~of}~G.\\[0.1cm]
\end{array}
\right.$ & $\left\{
\begin{array}{ll}
d_G(S)=\min\{e(T)\,|\,S\subseteq V(T)\},\\[0.1cm]
{\rm where}~T~{\rm is~
subtree~of}~G.\\[0.1cm]
\end{array}
\right.$  \\[0.1cm]
\cline{1-2}
 {\bf Eccentricity} & {\bf Steiner $k$-eccentricity} \\[0.1cm]
\cline{1-2}
$e(v)=\max\{d_G(u,v)\,|\,u\in V(G)\}$ & $e_k(v)=\max \{d(S)\,|\,S\subseteq V(G),$ \\[0.1cm]
 &    $|S|=k,~and~v\in
S \}$  \\[0.1cm]
\cline{1-2}
{\bf Diameter} & {\bf Steiner $k$-diameter} \\[0.1cm]
\cline{1-2}
$diam(G)=\max
\{e(v)\,|\,v\in V(G)\}$ & $sdiam_k(G)=\max \{e_k(v)\,|\,v\in V(G)\}$  \\[0.1cm]
\cline{1-2}
{\bf Radius} & {\bf Steiner $k$-radius} \\[0.1cm]
\cline{1-2}
 $rad(G)=\min\{e(v)\,|\,v\in V(G)\}$ & $srad_k(G)=\min \{
e_k(v)\,|\,v\in V(G)\}$ \\[0.1cm]
\cline{1-2}
\end{tabular}
\end{center}
\begin{center}
{Table 1.1. Classical distance parameters and Steiner distance parameters}
\end{center}
}

The following observation is immediate.
\begin{obs}\label{obs1-2}
Let $k,n$ be two integers with $2\leq k\leq n$.

$(1)$ If $H$ is a spanning subgraph of $G$, then $sdiam_k(G)\leq
sdiam_k(H)$.

$(2)$ For a connected graph $G$, $sdiam_k(G)\leq sdiam_{k+1}(G)$.
\end{obs}

\subsubsection{$k$-diameter}

Let $G$ be a $k$-connected graph, and let $u$, $v$ be any pair of
vertices of $G$. Let $P_k(u,v)$ be a family of $k$ inner vertex-disjoint
paths between $u$ and $v$, i.e., $P_k(u,v)=\{P_1,P_2,\cdots,P_k\}$,
where $p_1\leq p_2\leq \cdots \leq p_k$ and $p_i$ denotes the number
of edges of path $P_i$. The \emph{$k$-distance} $d_k(u,v)$ between
vertices $u$ and $v$ is the minimum $p_k$ among all $P_k(u,v)$ and
the \emph{$k$-diameter} $d_k(G)$ of $G$ is defined as the maximum
$k$-distance $d_k(u,v)$ over all pairs $u,v$ of vertices of $G$.
The concept of $k$-diameter emerges rather naturally when one looks at
the performance of routing algorithms. Its applications to network
routing in distributed and parallel processing are studied and
discussed by various authors including Chung \cite{Chung}, Du, Lyuu
and Hsu \cite{DuLyuuHsu}, Hsu \cite{Hsu, Hsu2}, Meyer and Pradhan
\cite{MeyerPradhan}.

\subsubsection{Steiner center and Steiner median}

The \emph{center} $C(G)$ of a connected graph $G$ is the subgraph induced by the
vertices $v$ of $G$ with $e(v)=rad(G)$.
As a generalization of the center of a graph, the \emph{Steiner
$k$-center} $C_k(G)\ (k\geq 2)$ of a connected graph $G$ is the
subgraph induced by the vertices $v$ of $G$ with $e_k(v)=srad_k(G)$.
Hence the Steiner $2$-center of a graph is simply its
center. The \emph{Steiner $k$-median} of $G$ is
the subgraph of $G$ induced by the vertices of minimum
Steiner $k$-distance in $G$. Similarly, Steiner $2$-median of a graph is simply its
median. For Steiner centers and Steiner medians, we
refer to \cite{Oellermann, Oellermann2, OellermannTian}.
{\small
\begin{center}
\begin{tabular}{|c|c|}
\cline{1-2}
\hline  {\bf Center} & {\bf Steiner $k$-center} \\[0.1cm]
\cline{1-2}
\hline The~subgraph~induced~by~the & The~subgraph~induced~by~the \\[0.1cm]
vertices~in~$\{v\in V(G)\,|\,e(v)=rad(G)\}$ & vertices~in~$\{v\in V(G)\,|\,e_k(v)=srad_k(G)\}$\\[0.1cm]
\cline{1-2}
{\bf Median} & {\bf Steiner $k$-median} \\[0.1cm]
\cline{1-2}
The~subgraph~induced~by~the~vertices & The~subgraph~induced~by~the~vertices  \\[0.1cm]
of~minimum~distance~in~$G$. & of~minimum
Steiner~$k$-distance~in~$G$.\\[0.1cm]
\cline{1-2}
\end{tabular}
\end{center}
\begin{center}
{Table 1.2. Center, median, and their generalizations.}
\end{center}
}

\subsubsection{Steiner intervals in graphs}

Let $G$ be a graph and $u,v$ two vertices of $G$. Then the {\it interval\/} from $u$ to $v$, $I_G(u,v)$ or
$I(u,v)$, is defined by
$$
I_G(u,v)=\{w\in V(G)\,|\,w~{\rm lies~on~a~shortest}~u-v~{\rm path~in}~G\}.
$$
Thus if $x\in I_G(u,v)$, then $d_G(u,u)=d_G(u,x)+d_G(x,v)$. In the case of a tree $T$ the interval
between two vertices $u$ and $v$ in $T$ consists of the vertices on the unique $u$-$v$ path in
$T$. However, in general the interval between two vertices $u$ and $v$ in a (connected) graph may contain vertices from more than just one shortest $u$-$v$ path. Mulder
\cite{Mulder} devoted an entire monograph to several topics related to intervals in graphs.

Motivated by the notion of the interval between two vertices and the
Steiner distance of a set of vertices, Kubicka, Kubicki, and Oellermann \cite{KubickaKubickiOellermann} defined the \emph{Steiner
interval}, $I_G(S)$ or $I(S)$, of a set $S$ by
$$
I_G(S)=\{w\in V(G)\,|\,w~{\rm lies~on~a~Steiner~tree~for}~S~{\rm in}~G\}.
$$
Thus if $|S|=2$, then the Steiner interval of $S$ is the interval between the two vertices
of $S$.

Let $S$ be an $k$-set and $r\leq k$. Then the {\it $r$-intersection interval\/} of $S$, denoted by $I_r(S)$ is
the intersection of all Steiner intervals of $r$-subsets of $S$. Graphs for which the $2$-
intersection intervals of every $3$-set consist of a unique vertex are
called {\it median graphs\/}. Median graphs were introduced independently by Avann \cite{Avann} and Neb\'{e}sky \cite{Nebesky}.

\subsubsection{Steiner distance
hereditary graphs}

Howorka \cite{Howorka} in 1977 defined a graph $G$ to be {\it
distance hereditary\/} if each connected induced subgraph $F$ of $G$
has the property that $d_F(u,v)=d_G(u,v)$ for each pair $u,v\in
V(F)$. As a generalization of distance hereditary graphs,
Day, Oellermann, and Swart \cite{DayOellermannSwart} introduced the concept of $k$-Steiner distance hereditary. A connected graph $G$ is
defined to be \emph{$k$-Steiner distance hereditary}, $k\geq 2$, if for
every connected induced subgraph $H$ of $G$ of order at least $k$
and a set $S$ of $k$ vertices of $H$, $d_H(S)=d_G(S)$. Thus,
2-Steiner distance hereditary graphs are distance hereditary.

\subsubsection{Steiner distance stable graphs}

In \cite{AliBoalsSherwani} a connected graph is defined to be {\it
vertex (edge) distance stable\/} if the distance between nonadjacent
vertices is unchanged after the deletion of a vertex (edge) of $G$. A
more general definition of an equivalent concept was introduced and
studied in \cite{EntringerJacksonSlater}. It was shown in
\cite{EntringerJacksonSlater} that a graph is vertex distance stable
if and only if it is edge distance stable. We will thus refer to
vertex or edge distance stable graphs as \emph{distance stable graphs}.
From a more general result established in
\cite{EntringerJacksonSlater}, it can be deduced that a graph is
distance stable if and only if for every pair $u,v$ of nonadjacent
vertices, $|N_G(u)\cap N_G(v)|=0$ or $|N_G(u)\cap N_G(v)|\geq 2$.
Thus a graph $G$ is \emph{distance stable} if and only if distances between
pairs of vertices in $G$ at distance $2$ apart remain unchanged
after the deletion of a vertex or an edge.

Further generalizations of distance stable graphs are studied in
\cite{GoddardOellermannSwart2}. In particular, for nonnegative
integers $k$ and $\ell$ not both zero and $D\subseteq N-\{\ell\}$ a
connected graph $G$ is defined to be {\it $(k,\ell,D)$-stable\/} if
for every pair $u,v$ of vertices of $G$ that are at distance
$d_G(u,v)\in D$ apart and every set $A$ consisting of at most $k$
vertices of $G-\{u,v\}$ and at most $\ell$ edges of $G$, the
distance between $u$ and $v$ in $G-A$ equals $d_G(u,v)$. For a
positive integer $m$, let $N_{\geq m}=\{x\in N\,|\,x\geq m\}$. In
\cite{GoddardOellermannSwart2} it is established that a graph is
$(k,\ell,\{m\})$-stable if and only if it is $(k,\ell,N_{\geq
m})$-stable. It is further shown that for a positive integer $x$ a
graph is $(k+x,\ell,\{2\})$-stable if and only if it is
$(k,\ell+x,\{2\})$-stable, but that $(k,\ell+x, \{m\})$-stable
graphs need not be $(k+x,\ell,\{m\})$-stable for $m\geq 4$.

The concepts of Steiner distance in graphs and distance stable
graphs suggest another generalization of distance stable graphs.
Goddard, Oellermann, and Swart \cite{GoddardOellermannSwart} introduced
the concept of Steiner distance stable graphs. We assume that $k,\ell,s$ and $m$ are
nonnegative integers with $m\geq s\geq 2$ and $k$ and $\ell$ not
both zero. If $S$ is a set of $s$ vertices in a connected graph $G$
such that $d_G(S)=m$, then $S$ is called an {\it $(s,m)$-set\/}. A
connected graph $G$ is said to be \emph{$k$-vertex $\ell$-edge
$(s,m)$-Steiner distance stable} if, for every $(s,m)$-set $S$ of $G$
and every set $A$ consisting of at most $k$ vertices of $G-S$ and at
most $\ell$ edges of $G$, $d_{G-A}(S)=d_G(S)$. Thus $k$-vertex
$\ell$-edge $(2,m)$-Steiner distance stable graphs are the $(k,\ell,
\{m\})$-stable graphs. Note that if $S$ is a set of $s$ vertices
such that $d_G(S)=s-1$ then $d_{G-A}(S)=d_G(S)$ for any set $A$ of
at most $k$ vertices of $G-S$ and at most $\ell$ edges of $G$. For
this reason we require that $m\geq s$.

For any integers $k,\ell,m$ and $s$ with $m\geq s\geq 2$ and $k$ and
$\ell$ not both $0$ there exists a $k$-vertex $\ell$-edge
$(s,m)$-Steiner distance stable graph. To see this let $G$ be
obtained from $m-1$ disjoint copies of $K_{k+\ell+1}$, say
$H_1,\ldots, H_{m-1}$, by joining every vertex of $H_i$ to every
vertex of $H_{i+1}$ for $1\leq i<m-1$ and then adding a vertex
$v_0$, and joining it to every vertex of $H_1$ and a vertex $v_{m}$,
and joining it to every vertex of $H_{m-1}$. It is not difficult to
see that $G$ is $k$-vertex $\ell$-edge $(s,m)$-Steiner distance
stable.

\subsubsection{Average Steiner distance}

Let $G=(V,E)$ be a connected graph of order $n$. The {\it average
distance\/} of $G$, $\mu(G)$, is defined to be the average of all
distances between pairs of vertices in $G$, i.e.
$$
\mu(G)={n\choose 2}^{-1}\sum_{\{u,v\}\subseteq V(G)}d_G(u,v).
$$
where $d_G(u,v)$ denotes the length of a shortest $u$-$v$ path in
$G$. This parameter, introduced in architecture as well as by the
chemist Wiener, turned out to be a good measure for analysing
transportation networks. It indicates the average time required to
transport a commodity between two destinations rather than the
maximum time required, as is indicated by the diameter; see \cite{DankelmannSwartOellermann}.

The {\it average Steiner distance\/} $\mu_k(G)$ of a graph $G$,
introduced by Dankelmann, Oellermann, and Swart in
\cite{DankelmannOellermannSwart}, is defined as the average of the Steiner
distances of all $k$-subsets of $V(G)$, i.e.
$$
\mu_k(G)={n\choose k}^{-1}\sum_{S\subseteq V(G), \ |S|=k}d_G(S).
$$
If $G$ represents a network, then the average Steiner $k$-distance indicates
the expected number of communication links needed to connect $k$
processors.

In \cite{BeinekeOellermannPippert}, the \emph{Steiner $k$-distance $(k\geq 2)$ of a vertex $v\in V(G)$} in a connected graph $G$ on $n\geq k$ vertices, denoted by $d_k(v,G)$, is defined by
$$
d_k(v)=\sum_{S\subseteq V(G), \ |S|=k, \ v\in S}d_G(S).
$$

\subsubsection{Steiner distance parameters in chemical graph theory}

In \cite{LiMaoGutman}, Li, Mao, and Gutman proposed a generalization of the Wiener index
concept, using Steiner distance. Thus, the \emph{$k$-th Steiner
Wiener index} $SW_k(G)$ of a connected graph $G$ is defined by
$$
SW_k(G)=\sum_{\overset{S\subseteq V(G)}{|S|=k}} d_G(S)\,.
$$
For $k=2$, the Steiner Wiener index coincides with the ordinary
Wiener index. It is usual to consider $SW_k$ for $2 \leq k \leq
n-1$, but the above definition implies $SW_1(G)=0$ and $SW_n(G)=n-1$
for a connected graph $G$ of order $n$. It should be noted that the average Steiner distance is related to the Steiner
Wiener index via $SW_k(G)/\binom{n}{k}$, that is,
$$
SW_k(G)=\sum_{S\subseteq V(G), \ |S|=k} d_G(S)=k^{-1}\sum_{v\in V(G)}d_k(v,G)=\binom{n}{k}\mu_k(G).
$$
For more details on Steiner Wiener index,
we refer to \cite{GutmanFurtulaLi, Kovse, LiMaoGutman, LiMaoGutman2, MaoWangGutman, MaoWangGutmanLi, MaoWangXiaoYe}.

Gutman \cite{GutmanSDD} offered an analogous
generalization of the concept of degree distance, Eq. (\ref{eq1-1}). Thus,
the {\it $k$-center Steiner degree distance\/} $SDD_k(G)$ of $G$ is
defined as
\begin{equation}                    \label{eq1-1}
SDD_k(G)=\sum_{\overset{S\subseteq V(G)}{|S|=k}}
\left[\sum_{v\in S}deg_G(v)\right] d_G(S),
\end{equation}
see \cite{GutmanSDD, MaoWangGutmanKlobucar, WangMaoDas} for more details.

Furtula, Gutman, and Katani\'{c} \cite{FurtulaGutmanKatanic}
introduced the concept of Steiner Harary index. The \emph{Steiner Harary $k$-index} or \emph{$k$-center Steiner Harary index} $SH_k(G)$ of $G$ is
defined as
$$
SH_k(G)=\sum_{\overset{S\subseteq V(G)}{|S|=k}}\frac{1}{d_G(S)},
$$
see \cite{FurtulaGutmanKatanic, Mao3} for more detials.

Mao and Das \cite{MaoDas} generalized the concept of Gutman index by Steiner distance. The \emph{Steiner
Gutman $k$-index} ${\rm SGut}_k(G)$ of $G$ is defined by
$$
{\rm SGut}_k(G)=\sum_{\overset{S\subseteq V(G)}{|S|=k}}\left(\prod_{v\in
S}deg_G(v)\right) d_G(S),
$$
see \cite{MaoDas, WangMaoDas2} for more detials.

The Steiner Wiener index, Steiner Harary index, the Steiner degree distance, and Steiner Gutman index are shown in the following Table 1.3.
{\small
\begin{center}
\begin{tabular}{|c|c|}
\cline{1-2}
\hline  {\bf Wiener index} & {\bf Steiner Wiener index} \\[0.1cm]
\cline{1-2}
$W(G)=\sum_{\{u,v\}\subseteq V(G)}d_G(u,v)$ & $SW_k(G)=\sum_{\overset{S\subseteq V(G)}{|S|=k}} d_G(S)$\\[0.1cm]
\cline{1-2}
\hline  {\bf Harary index} & {\bf Steiner Harary index} \\[0.1cm]
\cline{1-2}
$H(G)=\sum_{\{u,v\}\subseteq V(G)}\frac{1}{d_G(u,v)}$ & $SW_k(G)=\sum_{\overset{S\subseteq V(G)}{|S|=k}}\frac{1}{d_G(S)}$\\[0.1cm]
\cline{1-2}
{\bf Degree distance} & {\bf $k$-center Steiner degree distance} \\[0.1cm]
\cline{1-2}
$DD(G)=\sum_{\overset{}{\{u,v\}\subseteq V(G)}}[ deg_G(u)+deg_G(v)]d_G(u,v)$ & $SDD_k(G)=\sum_{\overset{S\subseteq V(G)}{|S|=k}}\left[\sum_{v\in
S}deg_G(v)\right] d_G(S)$ \\[0.1cm]
\cline{1-2}
{\bf Gutman index} & {\bf $k$-center Steiner Gutman index} \\[0.1cm]
\cline{1-2}
${\rm SGut}(G)=\sum_{\{u,v\}\subseteq V(G)}[deg_G(u)deg_G(v)]d_G(u,v)$ & ${\rm SGut}_k(G)=\sum_{\overset{S\subseteq V(G)}{|S|=k}}\left[\prod_{v\in
S}deg_G(v)\right] d_G(S)$ \\[0.1cm]
\cline{1-2}
\end{tabular}
\end{center}
\begin{center}
{Table 1.3. Four distance parameters and their generalizations.}
\end{center}
}

In this survey, we only summarize the known results on Steiner Wiener index.

\subsection{Some other Steiner structural parameters}

In this subsection, we introduce some other structural parameters related to Steiner distance.

\subsubsection{Steiner geodetic numbers}

Let $G$ be a connected graph and $u,v$ two vertices of $G$. Then the {\it interval\/} between $u$ and $v$, denoted by $I[u,v]$ is the union of all vertices that belong to some shortest $u$-$v$ path. The {\it closure\/} of a set $S$ of vertices in a connected graph $G$ is $\bigcup_{u,v\in S}I[u,v]$. If the closure of a set $S$ of vertices in $G$ is $V(G)$, then $S$ is called a {\it geodetic set\/}. The {\it geodetic number\/} of $G$, denoted by $g(G)$, is the smallest cardinality of a geodetic set in $G$.
It was shown in \cite{Atici} that the problem of finding the geodetic number of a graph is $NP$-hard.

The \emph{Steiner interval} of a set $S$ of vertices in a connected graph $G$, denoted by $I(S)$, is the union of all vertices of $G$
that lie on some Steiner tree for $S$.
A set $S$ of vertices in a connected graph $G$ is a {\it Steiner geodetic set\/} if $I(S)=V(G)$. The {\it Steiner geodetic number\/} of $G$, denoted by $sg(G)$, is the smallest cardinality of a Steiner geodetic set. As noted in \cite{Pelayo}, the result of \cite{DayOellermannSwart} that stated that $g(G)\leq sg(G)$ for all connected graphs $G$ it is not true.

\subsubsection{Steiner distance and convexity}

Abstract convexity started to develop in the early sixties, with the searching of an axiom system to define a set to be convex, in order to generalize, in some way, the classical concept of Euclidean convex set. These concepts can be found in \cite{Van de Vel}. Among the wide variety of structures that has been studied under this point of view, such as metric spaces, ordered sets or lattices, we are particularly interested in graphs, where several convexities associated to the vertex set are well-known.

Distance optimization properties of Steiner trees have given a way to define the Steiner
distance as a generalization of the usual distance in graphs; see \cite{ChartrandOellermannTianZou}.
Following that, C\'{a}ceres, M\'{a}rquez, and Puertas \cite{CaceresMarquezPuertas} define an abstract convexity in the context of graphs by means of the Steiner distance.

There are several well-known definitions of convex vertex sets in graphs, and these convexities
are usually defined by means of certain paths. In this fashion, a subset $S$ of vertices of a graph
$G$ is {\it monophonically\/} ({\it geodesically\/}) {\it convex\/} (see \cite{FarberJamison}) if $S$ contains every vertex of any chordless (shortest) path between vertices in $S$. These sets are called {\it $m$-convex sets\/} ({\it $g$-convex sets\/}).

The definitions above follow the general scheme of abstract convexities. A family $\mathscr{C}$ of subsets of a set $X$ is called a \emph{convexity} (see [15]) on $X$ if contains the empty set and universal set $X$, is
closed under intersections, and is closed under nested unions; that is, if $\mathscr{D}\subseteq \mathscr{C}$ is non-empty
and totally ordered by inclusion, then $\bigcup \mathscr{D}$ is in $\mathscr{C}$. Note that last property is trivial if $X$ is a
finite set. The elements of $\mathscr{C}$ are called \emph{convex sets}. It is clear that any subset $A$ of a convex
structure is included in a smallest convex set, $CH_{\mathscr{C}}(A)=\bigcap \{C\in \mathscr{C}: A\subseteq C\}$, called the \emph{convex
hull} of $A$. A point $p$ in a convex set $S$ is said to be an \emph{extreme point} if $S\setminus \{p\}$ is convex. The
preceding definitions correspond to $m$-convexity and $g$-convexity in graphs and their convex
hulls are denoted by $CH_m$ and $CH_g$ respectively.

Let $G$ be a connected graph. A subset $S\subseteq V(G)$ is said to be \emph{St-convex} if, for any $A\subseteq S$, all vertices in every Steiner tree of $A$ belong to $S$. The family of all St-convex sets of $V(G)$ defines a convexity called \emph{St-convexity}.

\subsection{Application backgrounds}

Steiner distance parameters have their application background in
both network science and mathematical chemistry.

\subsubsection{Application background of Steiner distance}

The Steiner tree problem in networks, and particularly in graphs,
was formulated in 1971-by Hakimi (see \cite{Hakimi})
and Levi (see \cite{Levi}). In the case of an unweighted, undirected
graph, this problem consists of finding, for a subset of vertices
$S$, a minimal-size connected subgraph that contains the vertices in
$S$. The computational side of this problem has been widely studied,
and it is known that it is an NP-hard problem for general graphs
(see \cite{HwangRW}). The determination of a Steiner tree in a graph
is a discrete analogue of the well-known geometric Steiner problem:
In a Euclidean space (usually a Euclidean plane) find the shortest
possible network of line segments interconnecting a set of given
points. Steiner trees have application to multiprocessor computer
networks. For example, it may be desired to connect a certain set of
processors with a subnetwork that uses the least number of
communication links. A Steiner tree for the vertices, corresponding
to the processors that need to be connected, corresponds to such a
desired subnetwork.

Steiner distance has application to multiprocessor communication. For example, suppose the
primary requirement when communicating a message from a processor $P$ to a collection
$S$ of other processors is to minimize the number of communication links that are
used. Then a Steiner tree for $S\cup \{P\}$ is an optimal way of connecting these vertices.
There are efficient algorithms that find the distance between two vertices or between
the vertices of the entire vertex set. However, if $S$ is a set of $k$ vertices, where
$2 <n <|V(G)|$, the only known algorithms that compute the Steiner distance of
$S$ have complexity which is a polynomial exponential in $k$. Indeed, it is known that the
general problem of finding the Steiner distance of a set is $NP$-hard (see \cite{GareyJohnson}), and
several heuristics for finding approximations to it have been developed (see Winter \cite{Winter}
for an extensive survey). However, if the graph G is a tree, then the Steiner distance of
any set $S$ and a Steiner tree for $S$ (it is unique) can be found efficiently.

\subsubsection{Application backgrounds of Steiner center, Steiner median, and average Steiner distance}

Graphs lend themselves as natural models of transportation networks as well as
communication and computer networks. Consequently, it is natural to study network
problems such as optimal facility location problems for graphs. In almost all such
problems, an optimal location is a point that is in some sense central to the network.
For example, the center of a connected graph is the subgraph induced by those
vertices for which the distance to the most remote vertex is least, and the median is the
subgraph induced by those vertices for which the sum of the distances to all of the
other vertices is least.

Any vertex in the center of a graph would be a suitable location
for an emergency facility, since the distance from the vertex to the furthest vertex from
it is minimized, whereas a vertex in the median is a good location for a service facility
since the average distance from that vertex to all other vertices is minimized. Slater [7]
has given an overview of a variety of other ways of determining centrality, and he
introduced and studied some new measures of centrality.

If $G$ represents a network, then the average Steiner $k$-distance indicates the expected    number of communication links needed to connect $k$ processors. In contrast the Steiner $k$-diumeter of $G$, $sdiam_k(G)$, defined as the maximum distance of the $k$-subsets of  $V(G)$, indicates the number of communication links needed in the worst case.

\subsubsection{Application background in chemical graph theory}

Recall that Wiener index $W(G)$ of the graph $G$ is
$W(G)=\sum_{\{u,v\} \subseteq V(G)} d_G(u,v)$.
Details on this oldest distance--based topological index can be
found in numerous surveys, e.g., in
\cite{DobryninKochetova, Rouvray1, Rouvray2, Xu}.
Li, Mao, and Gutman \cite{LiMaoGutman} put forward a Steiner--distance--based generalization of the
Wiener index concept. A chemical application of $SW_k$ was recently reported in \cite{GutmanFurtulaLi}, where it is shown that the term $W(G)+\lambda SW_k(G)$ provides a better approximation for the boiling points of alkanes than $W(G)$
itself, and that the best such approximation is obtained for $k=7$.
Furtula, Gutman, and Katani\'{c} \cite{FurtulaGutmanKatanic}
introduced the concept
of Steiner Harary index and gave its chemical applications.

\section{Steiner Distance}

The classical distance defined on a connected graph $G$ is a metric on its vertex set.
As such, certain properties are satisfied. Among these are:
\begin{itemize}
\item[] $(i)$ $d(u,v)=0$ for vertices
$u,v$ of $G$ and  $d(u,v)=0$ if and only if $u=v$;

\item[] $(ii)$ $d(u, w)\leq d(u,v)+d(v,w)$ for vertices $u,v,w$ of $G$.
\end{itemize}

The extensions of these properties to the Steiner distance are given by Chartrand, Oellermann, Tian, and Zou \cite{ChartrandOellermannTianZou}.
\begin{itemize}
\item[] $(1)$ Let $G$ be a connected graph and let $S\subseteq V(G)$, where $S\neq \emptyset$.  Then $d_G(S)\geq 0$.
Furthermore, $d_G(S)=0$ if and only if $|S|=1$. This is an extension of $(i)$.

\item[] $(2)$  To provide an extension of $(ii)$, let $S,S_1$ and $S_2$ be subsets of $V(G)$ such that $\emptyset \neq S  \subseteq S_1\cup S_2$ and $S_1\cap S_2\neq \emptyset$. Then $d_G(S)\leq d_G(S_1)+d_G(S_2)$ To  see this, let $T_i \ (i=1,2)$ be a tree of
size $d_G(S_i)$ such that $S_i\subseteq V(T_i)$. Let $H$ be the graph with vertex set $V(T_1)  \cup V(T_2)$ and edge set $E(T_1)\cup E(T_2)$. Since $T_1$ and $T_2$ are connected and $V(T_1)\cap V(T_2) \neq \emptyset$, the graph $H$ is connected. Since $S\subseteq V(H)$,
$d_G(S)\leq e(H)\leq d_G(S_1)+d_G(S_2)$.
\end{itemize}

It is useful to observe that if $T$ is a nontrivial
tree and $S\subseteq V(T)$, where $|S|\geq 2$, then there is a unique subtree $T_s$ of size  $d_{T}(S)$ containing the vertices of $S$. We refer to such a tree as the \emph{tree generated by $S$}.
\begin{obs}{\upshape \cite{ChartrandOellermannTianZou}}\label{obs2-1}
If $S$ is a set of vertices of a tree $T$ and $v$ is a vertex in $V(T)-S$,
then the tree generated by $S\cup \{v\}$ contains the tree generated by $S$.
\end{obs}

\begin{obs}{\upshape \cite{ChartrandOellermannTianZou}}\label{obs2-2}
Let $w$ be the
(necessarily unique) vertex of $T_S$ whose distance from $v$ is a minimum. Then $T_{S\cup \{v\}}$ contains the unique $v$-$w$ path and
$$
d_T(S\cup \{v\})=d_T(S)+d_T(v,w).
$$
\end{obs}

\begin{obs}{\upshape \cite{ChartrandOellermannTianZou}}\label{obs2-3}
If $H$ is a subgraph of a graph $G$ and $v$ is a vertex of $G$, then
$$
d_G(S\cup \{v\})=d_G(S)+d_G(v,T_S),
$$
$d(v,H)$ denotes the minimum distance from $v$ to a vertex of $H$.
\end{obs}

For a tree $T$, we denote by $V_1(T)$ the set of end-vertices of $T$ and $n_1=|V_1(T)|$.
If $S=V_1(T)$, then $T_S=T$ so that $d_T(S)=e(T)$ and $d_T(S\cup \{v\})=e(T)$ for all
$v\in V(T)$. Hence if $T$ is a tree and $k\geq 2$ and integer with $n_1<k$, then
$e_k(v)=e(T)$ for all $v\in V(T)$.

Chartrand, Oellermann, Tian, and Zou \cite{ChartrandOellermannTianZou} derived the following results for Steiner distance.
\begin{pro}{\upshape \cite{ChartrandOellermannTianZou}}\label{pro2-1}
Let $k\geq 2$ be an integer and suppose that $T$ is a tree of order $n$
with $n_1\geq k$. Let $v\in V(T)$. If $S\subseteq V(T)$ such that $v\notin S$, $|S|=k-1$ and $d_T(S\cup \{v\})=e_k(v)$, $S\subseteq V_1(T)$.
\end{pro}

\begin{cor}{\upshape \cite{ChartrandOellermannTianZou}}\label{cor2-1}
Let $k\geq 2$ be an integer and $T$ a tree with $n_1\geq k$. Then $sdiam_k(T)=d_T(S)$, where $S$ is a set of $k$ end-vertices of $T$.
\end{cor}

\begin{pro}{\upshape \cite{ChartrandOellermannTianZou}}\label{pro2-1}
Let $k\geq 3$ be an integer and suppose that $T$ is a tree with $n_1\geq k$
end-vertices. If $v$ is a vertex of $T$ with $e_k(v)=srad_k(T)$, then there exists a set $S$  of $k-1$ end-vertices of $T$ such that $d_T(S\cup \{v\})=e_k(v)$ and $v\in V(T_S)$.
\end{pro}

\begin{cor}{\upshape \cite{ChartrandOellermannTianZou}}\label{cor2-2}
Let $k\geq 3$ be an integer and suppose that $T$ is a tree with at least $k$
end-vertices. If $v$ is a vertex of $T$ with $e_k(v)=srad_k(T)$, then $v$ is not an  end-vertex of $T$.
\end{cor}

Let $k$ and $n$ be integers with $2\leq k\leq n$. A graph $G$ of order $n$ is called an \emph{$(k; n)$-graph} if it is of minimum size with the property that $d_G(S)=k-1$ for all sets $S$
of vertices of $G$ with $|S|=k$.

Chartrand, Oellermann, Tian, and Zou \cite{ChartrandOellermannTianZou} determined the size of an $(k;n)$ graph for each pair $k,n$ of integers with $2\leq k\leq n$.

\begin{thm}{\upshape \cite{ChartrandOellermannTianZou}}\label{th2-1}
Let $k$ and $n$ be integers with $2\leq k\leq n$. The size of an $(k;n)$-graph
is $k-1$ if $k=n$ and $[(n-k+1)n/2]$ if $n>k$.
\end{thm}

\subsection{Steiner distance in some graph classes}

The following observation is easy to make from the definition of a
threshold graph.
\begin{obs}\label{obs2-1}
Let $G([n],E)$ be a threshold graph with a weight function $w:
V(G)\rightarrow R$. Let the vertices be labelled so that $w(1)\geq
w(2)\geq \cdots \geq w(n)$. Then

$(1)$ $d_1\geq d_2\geq \cdots \geq d_n$, where $d_i$ is the degree
of vertex $i$.

$(2)$ $I=\{i\in V(G): d_i\leq i-1\}$ is a maximum independent set of
$G$ and $G\setminus I$ is a clique in $G$.

$(3)$ $N(i)=\{1,2,\cdots,d_i\}$ for every $i\in I$. Thus, the
neighborhoods of vertices in $I$ form a linear order under set
inclusion. Furthermore, if $G$ is connected, then every vertex in
$G$ is adjacent to $1$.
\end{obs}

Let $C_r$ and $I_{n-r}$ denote the clique and the maximum
independent set of $G$, respectively, with
$V(C_r)=\{u_1,u_2,\ldots,u_r\}$ and
$V(I_{n-r})=\{u_1',u_2',\ldots,u_{n-r}'\}$ such that $deg_G(u_1)\geq
deg_G(u_2)\geq \cdots \geq deg_G(u_r)$ and $deg_G(u_1')\geq deg_G(u_2')\geq
\cdots \geq deg_G(u_{n-r}')$.

Wang, Mao, Cheng, and Melekian \cite{WangMaoChengMelekian} derived the following
results for the Steiner distance of threshold graphs.
\begin{pro}{\upshape \cite{WangMaoChengMelekian}}\label{pro2-3}
Let $k,n$ be two integers with $3\leq k\leq n$, and let $G$ be a
threshold graph of order $n$. Let $S$ be a set of distinct vertices
of $G$ such that $|S|=k$. Let $u_i$ be the vertex in $S\cap V(C_r)$
with the minimum subscript, and $u_j'$ be the vertex in $S\cap
V(I_{n-r})$ with the maximum subscript.

$(1)$ If $S\subseteq V(C_r)$, then $d_{G}(S)=k-1$.

$(2)$ If $S\subseteq V(I_{n-r})$, then $d_{G}(S)=k$.

$(3)$ If $S\cap V(C_r)\neq \emptyset$, $S\cap V(I_{n-r})\neq
\emptyset$, and $u_iu_j'\in E(G)$, then $d_{G}(S)=k-1$.

$(4)$ If $S\cap V(C_r)\neq \emptyset$, $S\cap V(I_{n-r})\neq
\emptyset$, and $u_iu_j'\notin E(G)$, then $d_{G}(S)=k$.
\end{pro}

For two graphs $G$ and $H$ with
$V(G)=\{u_1,u_2,\ldots,u_{n}\}$ and $V(H)=\{v_1,v_2,\ldots,v_{m}\}$,
from the definition of corona graphs, $V(G*H)=V(G)\cup
\{(u_i,v_j)\,|\,1\leq i\leq n, \ 1\leq j\leq m\}$, where $*$ denotes
the corona product operation. For $u\in V(G)$, we use $H(u)$ to
denote the subgraph of $G*H$ induced by the vertex set
$\{(u,v_j)\,|\,1\leq j\leq m\}$. For fixed $i \ (1\leq i\leq n)$, we
have $u_i(u_i,v_j)\in E(G*H)$ for each $j \ (1\leq j\leq m)$. Then
$V(G*H)=V(G)\cup V(H(u_1))\cup V(H(u_2))\cup \ldots \cup V(H(u_n))$.

Wang, Mao, Cheng, and Melekian \cite{WangMaoChengMelekian} also derived the following
results on the Steiner distance of corona graphs.
\begin{thm}{\upshape \cite{WangMaoChengMelekian}}\label{th2-2}
Let $k,m,n$ be three integers with $3\leq k\leq n(m+1)$, and let
$G,H$ be two connected graphs with $V(G)=\{u_1,u_2,\ldots,u_{n}\}$
and $V(H)=\{v_1,v_2,\ldots,v_{m}\}$. Let $S$ be a set of distinct
vertices of $G*H$ such that $|S|=k$.
$$
d_{G*H}(S)=d_{G}(S_G')+k-t,
$$
where $|S\cap V(G)|=t$, and $S_G'$ is the maximum subset of $V(G)$
such that $S\cap (V(H(u))\cup \{u\})\neq \emptyset$ for each $u\in
S_G'$.
\end{thm}

From the definition of cluster, $V(G\odot
H)=\{(u_i,v_j)\,|\,1\leq i\leq n, \ 1\leq j\leq m\}$, where $\odot$
denotes the cluster product operation. For $u\in V(G)$, we use
$H(u)$ to denote the subgraph of $G\odot H$ induced by the vertex
set $\{(u,v_j)\,|\,1\leq j\leq m\}$. Without loss of generality, we
assume $(u_i,v_1)$ is the root of $H(u_i)$ for each $u_i\in V(G)$.
Let $G(h_1)$ be the graph induced by the vertices in
$\{(u_i,v_1)\,|\,1\leq i\leq n\}$. Clearly, $G(u_1)\simeq G$, and
$V(G\odot H)=V(H(u_1))\cup V(H(u_2))\cup \ldots \cup V(H(u_n))$.

Wang, Mao, Cheng, and Melekian \cite{WangMaoChengMelekian} obtained the following
results on the Steiner distance of cluster graphs.
\begin{thm}{\upshape \cite{WangMaoChengMelekian}}\label{th2-3}
Let $k,m,n$ be three integers with $3\leq k\leq n(m+1)$, and let
$G,H$ be two connected graphs with $V(G)=\{u_1,u_2,\ldots,u_{n}\}$
and $V(H)=\{v_1,v_2,\ldots,v_{m}\}$. Let
$S=\{(u_{i_1},v_{j_1}),(u_{i_2},v_{j_2}),\ldots,(u_{i_{k}},v_{j_{k}})\}$
be a set of distinct vertices of $G\odot H$. Let
$S_G=\{u_{i_1},u_{i_2},\ldots,u_{i_{k}}\}$ and
$S_H=\{v_{j_1},v_{j_2},\ldots,v_{j_{k}}\}$.

$(1)$ If $S\subseteq V(G(v_1))$, then $d_{G\odot H}(S)=d_G(S_G)$.

$(2)$ If there exists some $H(u_i)\ (1\leq i\leq n)$ such that
$S\subseteq V(H(u_i))$, then $d_{G\odot H}(S)=d_H(S_H)$.

$(3)$ If there is no $H(u_i)\ (1\leq i\leq n)$ such that
$S\subseteq V(H(u_i))$, then
$$
d_G(S_G')+k-t\leq d_{G\odot H}(S)\leq \left\{
\begin{array}{ll}
rd_H(S_H)+d_G(S_G') &\mbox v_1\in S_H,\\[0.2cm]
rd_H(S_H\cup \{h_1\})+d_G(S_G') &\mbox v_1\notin S_H,
\end{array}
\right.
$$
where $|S\cap V(G(v_1))|=t$, $|S_G'|=r$, and $S_G'$ is the maximum
subset of $V(G)$ such that $S\cap V(H(u))\neq \emptyset$ for each
$u\in S_G'$.

Moreover, the upper and lower bounds are sharp.
\end{thm}

To show the sharpness of the above lower and upper bounds, Wang, Mao, Cheng, and Melekian \cite{WangMaoChengMelekian}
considered the following examples.

\noindent{\bf Example 2.1.} \cite{WangMaoChengMelekian} Let $G=P_n=u_1u_2\ldots u_n$ and
$H=P_m=v_1v_2\ldots v_m$ with $3\leq k\leq mn$. Note that
$H(u_i)\cong P_m$ for each $u_i \ (1\leq i\leq n)$, and $G(v_1)\cong
P_n$. For $v_1\notin S_H$, if $k\leq n$, then we choose
$S=\{(u_1,v_m),(u_2,v_m),\ldots,(u_{k-1},v_m)\}\cup
\{(u_{n},v_m)\}$. Then $r=k$, $d_H(S_H\cup \{v_1\})=m-1$,
$d_G(S_G')=n-1$. Since the tree induced by the edges in
$E(G(v_1))\cup E(H(u_1))\cup E(H(u_2))\cup \ldots \cup
E(H(u_{k-1}))\cup E(H(u_{n}))$ is the unique $S$-Steiner tree, it
follows that $d_{G\odot H}(S)\geq k(m-1)+(n-1)$. From Theorem
\ref{th2-3}, $d_{G\odot H}(S)\leq rd_H(S_H\cup
\{v_1\})+d_G(S_G')=k(m-1)+(n-1)$. So, the upper bound for $v_1\notin
S_H$ is sharp. For $v_1\in S_H$, if $k\leq n$, then we choose
$S=\{(u_1,v_1),(u_2,v_1),\ldots,(u_{k},v_1)\}$. Then $r=k$,
$d_H(S_H)=0$, $d_G(S_G')=k-1$. Then $d_{G\odot H}(S)\geq k-1$. From
Theorem \ref{th2-3}, $d_{G\odot H}(S)\leq rd_H(S_H)+d_G(S_G')=k-1$.
So, the upper bound for $v_1\in S_H$ is sharp.

\vskip 0.3cm

\noindent{\bf Example 2.2.} \cite{WangMaoChengMelekian} Let $G=P_n=u_1u_2\ldots u_n$ and $H=K_m$
with $3\leq k\leq mn$, where $V(H)=\{v_1,v_2,\ldots,v_n\}$. Note
that $H(u_i)\cong K_m$ for each $u_i \ (1\leq i\leq n)$, and
$G(v_1)\cong P_n$. Choose
$S=\{(u_1,v_m),(u_2,v_m),\ldots,(u_{k-1},v_m)\}\cup \{(u_{n},v_m)\}
\ (m\geq 2)$. Then $d_G(S_G')=n-1$ and $t=0$, and hence $d_{G\odot
H}(S)\geq n-1+k$. Clearly, the tree induced by the edges in
$E(G(v_1))\cup E_{G\odot H}[V(G(v_1)),S]$ is an $S$-Steiner tree in
$G\odot H$, and hence $d_{G\odot H}(S)\leq n-1+k$. So, we have
$d_{G\odot H}(S)=n-1+k$, which implies that the lower bound is
sharp.

\begin{cor}{\upshape \cite{WangMaoChengMelekian}}\label{cor2-3}
Let $k,m,n$ be three integers with $3\leq k\leq n(m+1)$, and let
$G,H$ be two connected graphs with $V(G)=\{u_1,u_2,\ldots,u_{n}\}$
and $V(H)=\{v_1,v_2,\ldots,v_{m}\}$. Let $(u_i,v_1)$ be the root of
$H(u_i)$ for each $u_i \ (1\leq i\leq n)$. Let $S$ be a set of
distinct vertices of $G\odot H$ such that $|S|=k$. Then
$$
d_{G\odot H}(S)\leq rd_H(S_H\cup \{v_1\})+d_G(S_G'),
$$
where $|S_G'|=r$, and $S_G'$ is the maximum subset of $V(G)$ such
that $S\cap V(H(u))\neq \emptyset$ for each $u\in S_G'$.
\end{cor}

\subsection{Steiner distance of graph products}

Product networks were proposed based upon the idea of using the
product as a tool for ``combining'' two known graphs with
established properties to obtain a new one that inherits properties
from both \cite{DayA}.

Wang, Mao, Cheng, and Melekian \cite{WangMaoChengMelekian} gave the exact value for Steiner distance of joined graphs.
\begin{pro}{\upshape \cite{WangMaoChengMelekian}}\label{pro2-4}
Let $k,m,n$ be three integers with $3\leq k\leq m+n$, and let $G,H$
be two connected graphs with $n,m$ vertices, respectively. Let $S$
be a set of distinct vertices of $G\vee H$ such that $|S|=k$.

$(1)$ If $S\cap V(G)\neq \emptyset$ and $S\cap V(H)\neq \emptyset$,
then $d_{G\vee H}(S)=k-1$.

$(2)$ If $S\cap V(H)=\emptyset$ and $G[S]$ is connected, then
$d_{G\vee H}(S)=k-1$; if $S\cap V(H)=\emptyset$ and $G[S]$ is not
connected, then $d_{G\vee H}(S)=k$.

$(3)$ If $S\cap V(G)=\emptyset$ and $H[S]$ is connected, then
$d_{G\vee H}(S)=k-1$; if $S\cap VG)=\emptyset$ and $H[S]$ is not
connected, then $d_{G\vee H}(S)=k$.
\end{pro}

In \cite{Gologranc}, Gologranc obtained the following lower bound for Steiner distance.
\begin{thm}{\upshape \cite{Gologranc}}\label{th2-4}
Let $k\geq 2$ be an integer, and let $G,H$ be two connected graphs.
Let
$S=\{(u_{i_1},v_{j_1}),(u_{i_2},v_{j_2}),\ldots,(u_{i_k},v_{j_k})\}$
be a set of distinct vertices of $G\Box H$. Let
$S_G=\{u_{i_1},u_{i_2},\ldots,u_{i_k}\}$ and
$S_H=\{v_{j_1},v_{j_2},\ldots,v_{j_k}\}$. Then
$$
d_{G\Box H}(S)\geq d_G(S_G)+d_H(S_H).
$$
\end{thm}

\begin{figure}[!hbpt]
\begin{center}
\includegraphics[scale=0.7]{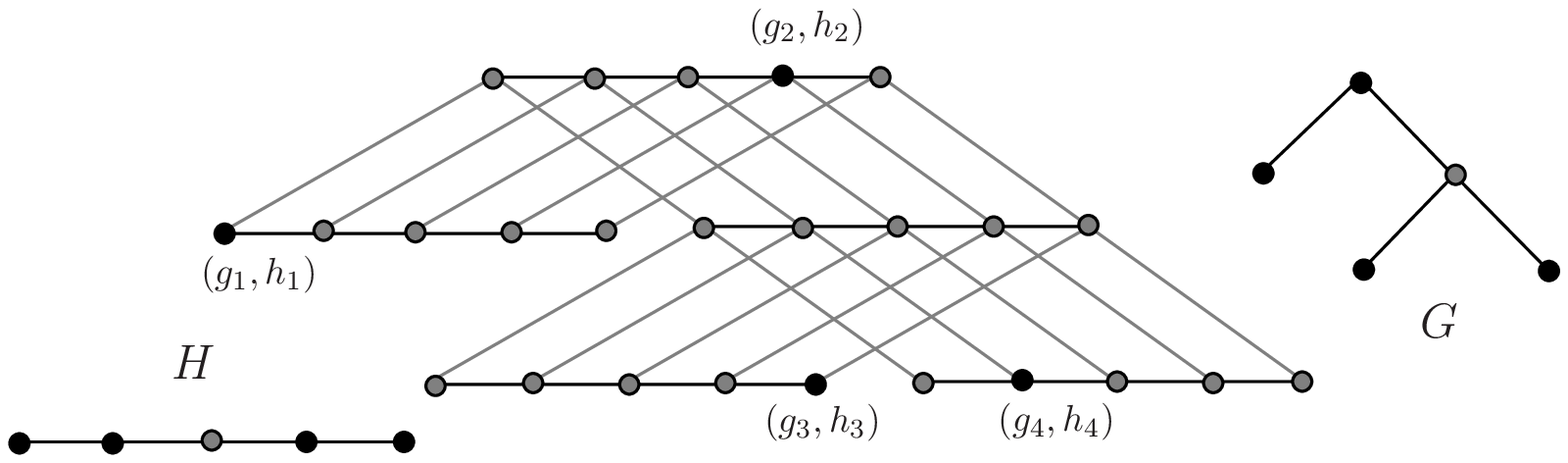}
\end{center}
\begin{center}
Figure 2.1: Graphs for Remark 2.1.
\end{center}\label{fig1}
\end{figure}

Mao, Cheng, and Wang \cite{MaoChengWang} showed that the inequality in Theorem \ref{th2-4} can be
equality if $k=3$; shown in following Corollary \ref{cor2-5}. But,
for general $k \ (k\geq 4)$, from Theorem \ref{th2-4} and Corollary
\ref{cor2-5}, one may conjecture that for two connected graphs
$G,H$, $d_{G\Box H}(S)=d_G(S_G)+d_H(S_H)$, where
$S=\{(u_{i_1},v_{j_1}),(u_{i_2},v_{j_2}),\ldots,(u_{i_k},v_{j_k})\}\subseteq
V(G\Box H)$, $S_G=\{u_{i_1},u_{i_2},\ldots,u_{i_k}\}\subseteq V(G)$
and $S_H=\{v_{j_1},v_{j_2},\ldots,v_{j_k}\}\subseteq V(H)$.

\begin{rem}{\upshape \cite{WangMaoChengMelekian}}\label{rem2-1}
Actually, the
equality $d_{G\Box H}(S)=d_G(S_G)+d_H(S_H)$ is not true for $|S|\geq 4$. For example, let $G$ be a tree
with degree sequence $(3,2,1,1,1)$ and $H$ be a path of order $5$.
Let $S=\{(u_1,v_1),(u_2,v_2),(u_3,v_3),(u_4,v_4)\}$ be a vertex set
of $G\Box H$ shown in Figure 2.1. Then $d_G(S_G)=4$ for
$S_G=\{u_1,u_2,u_3,u_4\}$, and $d_H(S_H)=4$ for
$S_H=\{v_1,v_2,v_3,v_4\}$. One can check that there is no
$S$-Steiner tree of size $8$ in $G\Box H$, which implies $d_{G\Box
H}(S)\geq 9$.
\end{rem}

Although the conjecture of such an ideal formula is not correct, it
is possible to give a strong upper bound for general $k \ (k\geq 3)$. Remark 2.1 also indicates
that obtaining a nice formula for the general case may be difficult.

Mao, Cheng, and Wang \cite{MaoChengWang} got such an upper bound of $d_{G\Box H}(S)$ for $S\subseteq
V(G\Box H)$ and $|S|=k$.
\begin{thm}{\upshape \cite{MaoChengWang}}\label{th2-5}
Let $k,m,n$ be three integers with $3\leq k\leq mn$, and let $G,H$
be two connected graphs with $V(G)=\{u_1,u_2,\ldots,u_{n}\}$ and
$V(H)=\{v_1,v_2,\ldots,v_{m}\}$. Let
$S=\{(u_{i_1},v_{j_1}),(u_{i_2},v_{j_2}),\ldots,(u_{i_k},v_{j_k})\}$
be a set of distinct vertices of $G\Box H$, $S_G=\{u_{i_1},u_{i_2},$
$\ldots,u_{i_k}\}$, and $S_H=\{v_{j_1},v_{j_2},\ldots,v_{j_k}\}$,
where $S_G\subseteq V(G)$, $S_H\subseteq V(H)$ ($S_G,S_H$ are both
multi-sets). Then
\begin{eqnarray*}
d_G(S_G)+d_H(S_H)&\leq &d_{G\Box H}(S)\\[0.1cm]
&\leq&\min\{d_G(S_G)+(r+1)d_H(S_H),d_H(S_H)+(t+1)d_G(S_G)\},
\end{eqnarray*}
where $r,t \ (0\leq r,t\leq k-3)$ are defined as follows.

$\bullet$ Let $X_G^{i} \ (1\leq i\leq {k\choose 3})$ be all the
$(k-3)$-multi-subsets of $\{u_{i_1},u_{i_2},,\ldots,u_{i_k}\}$ in
$G$, and let $r_i$ be the numbers of distinct vertices in $X_G^{i} \
(1\leq i\leq {k\choose 3})$, and let $r=\max\{r_i\,|\,1\leq i\leq
{k\choose 3}\}$.

$\bullet$ Let $Y_H^{j} \ (1\leq j\leq {k\choose 3})$ be all the
$(k-3)$-multi-subsets of $\{v_{j_1},v_{j_2},\ldots,v_{j_k}\}$ in
$H$, and let $t_j$ be the numbers of distinct vertices in $Y_H^{j} \
(1\leq j\leq {k\choose 3})$, and let $t=\max\{t_j\,|\,1\leq j\leq
{k\choose 3}\}$.
\end{thm}

The following corollaries are immediate from Theorem \ref{th2-5}.
\begin{cor}{\upshape \cite{MaoChengWang}}\label{cor2-4}
Let $G,H$ be two connected graphs of order $n,m$, respectively. Let
$k$ be an integer with $3\leq k\leq mn$. Let
$S=\{(u_{i_1},v_{j_1}),(u_{i_2},v_{j_2}),\ldots,(u_{i_k},v_{j_k})\}$
be a set of distinct vertices of $G\Box H$. Let
$S_G=\{u_{i_1},u_{i_2},\ldots,u_{i_k}\}$ and
$S_H=\{v_{j_1},v_{j_2},\ldots,v_{j_k}\}$. Then
\begin{eqnarray*}
d_G(S_G)+d_H(S_H)&\leq &d_{G\Box H}(S)\\[0.1cm]
&\leq&\min\{d_G(S_G)+(k-2)d_H(S_H),d_H(S_H)+(k-2)d_G(S_G)\}\\[0.1cm]
&=&d_G(S_G)+d_H(S_H)+(k-3)\min\{d_H(S_H),d_G(S_G)\}.
\end{eqnarray*}
\end{cor}

\begin{cor}{\upshape \cite{MaoChengWang}}\label{cor2-5}
Let $G,H$ be two connected graphs, and let $(u,v)$, $(u',v')$ and
$(u'',v'')$ be three vertices of $G\Box H$. Let $S_G=\{u,u',u''\}$,
$S_H=\{v,v',v''\}$, and $S=\{(u,v),(u',v'),$ $(u'',v'')\}$. Then
$$
d_{G\Box H}(S)=
d_{G}(S_G)+
d_{H}(S_H)
$$
\end{cor}

To show the sharpness of the above upper and lower bound, Mao, Cheng, and Wang \cite{MaoChengWang}
considered the following example.

\noindent{\bf Example 2.3.} \cite{MaoChengWang} $(1)$ For $k=3$, from Corollary
\ref{cor2-5}, we have $d_{G\Box H}(S)=d_G(S_G)+d_H(S_H)$, which
implies that the upper and lower bounds in Corollary \ref{cor2-4}
and Theorem \ref{th2-5} are sharp.

$(2)$ Let $G,H$ be two connected graphs. Suppose $4\leq k\leq
|V(H)|$. Choose $k$ vertices in $H(u)$ for fixed $u\in V(G)$, say
$(u,v_1),(u,v_2),\cdots,(u,v_k)$. Clearly, $d_G(S_G)=0$. From
Corollary \ref{cor2-4}, $d_{G\Box H}(S)\geq
d_G(S_G)+d_H(S_H)=d_H(S_H)$ and $d_{G\Box H}(S)\leq
d_G(S_G)+d_H(S_H)+(k-3)\min\{d_H(S_H),d_G(S_G)\}=d_H(S_H)$, and
hence $d_{G\Box H}(S)=d_H(S_H)$. This implies that the upper and
lower bounds in Corollary \ref{cor2-4} are sharp. From the
definition of $r$, we have $r=1$ and $d_{G\Box H}(S)\leq
d_G(S_G)+d_H(S_H)+r\min\{d_H(S_H),d_G(S_G)\}=d_H(S_H)$, and hence
$d_{G\Box H}(S)=d_H(S_H)$. This also implies that the upper and
lower bounds in Theorem \ref{th2-5} are sharp.

$(3)$ Let $G=P_n$ and $H=P_m$ with $n\leq m$, where
$P_n=u_1u_2\cdots u_n$ and $P_m=v_1v_2\cdots v_m$. Choose
$S=\{(u_1,v_1),(u_1,v_m),(u_n,v_1),(u_n,v_m)\}$. Then
$d_G(S_G)=n-1$, $d_H(S_H)=m-1$ and $d_{G\Box
H}(S)=2(n-1)+(m-1)=d_G(S_G)+d_H(S_H)+(4-3)\min\{d_H(S_H),d_G(S_G)\}$,
which implies that the upper bound in Corollary \ref{cor2-4} are
sharp.

From the definition, the lexicographic product graph $G\circ H$ is a
graph obtained by replacing each vertex of $G$ by a copy of $H$ and
replacing each edge of $G$ by a complete bipartite graph $K_{m,m}$,
where $m=|V(H)|$.
\begin{thm}{\upshape \cite{Hammack}}\label{th2-6}
Let $(u,v)$ and $(u',v')$ be two vertices of $G\circ H$. Then
$$
d_{G\circ H}((u,v),(u',v'))=\left\{
\begin{array}{ll}
d_{G}(u,u'), &\mbox {\rm if} \ u\neq u';\\[0.2cm]
d_{H}(v,v'), &\mbox {\rm if} \ u=u'~and~deg_{G}(u)=0;\\[0.2cm]
\min \{d_{H}(v,v'),2\}, &\mbox {\rm if} \ u=u'~and~deg_{G}(u)\neq 0.
\end{array}
\right.
$$
\end{thm}

Anand, Changat, Klav\v{z}ar, and Peterin \cite{AnandChangatKlavzarPeterin} obtained the following formula.
\begin{thm}{\upshape \cite{AnandChangatKlavzarPeterin}}\label{th2-7}
Let $k\geq 2$. Let
$S=\{(u_{i_1},v_{j_1}),(u_{i_2},v_{j_2}),\ldots,(u_{i_k},v_{j_k})\}$
be a set of distinct vertices of $G\circ H$ such that $u_{i_p}\neq
u_{i_q} \ (1\leq p,q\leq k)$. Let
$S_G=\{u_{i_1},u_{i_2},\ldots,u_{i_k}\}$. Then
$$
d_{G\circ H}(S)=d_{G}(S_G).
$$
\end{thm}

For general case, Mao, Cheng, and Wang \cite{MaoChengWang} had the following formula for Steiner distance
of lexicographic product graphs.
\begin{thm}{\upshape \cite{MaoChengWang}}\label{th2-8}
Let $k,n,m$ be three integers with $2\leq k\leq mn$. Let $G$ be a
connected graph of order $n$, and $H$ be a graph of order $m$. Let
$S=\{(u_{i_1},v_{i_1}),(u_{i_2},v_{i_2}),\ldots,(u_{i_k},v_{i_k})\}$
be a set of distinct vertices of $G\circ H$. Let
$S_G=\{u_{i_1},u_{i_2},\ldots,u_{i_k}\}$ and
$S_H=\{v_{j_1},v_{j_2},\ldots,v_{j_k}\}$ (note that $S_G,S_H$ are
both multi-set). Let $r$ be the number of distinct vertices in
$S_G$, where $1\leq r\leq k$.

$(1)$ If $r=1$ and $H[S_H]$ is connected in $H$, then $d_{G\circ
H}(S)=k-1$.

$(2)$ If $r=1$ and $H[S_H]$ is not connected in $H$, then $d_{G\circ
H}(S)=k$.

$(3)$ If $r\geq 2$, then $d_{G\circ H}(S)=d_{G}(S_G)+k-r$.
\end{thm}

In Theorem \ref{th2-8}, they assumed that $G$ is a connected graph. For
$k=3$, they had the following by assuming that $G$ is not connected.
\begin{pro}{\upshape \cite{MaoChengWang}}\label{pro2-5}
Let $G$ and $H$ be two graphs, and let $(u,v)$, $(u',v')$ and
$(u'',v'')$ be three vertices of $G\circ H$. Let
$S=\{(u,v),(u',v'),(u'',v'')\}$, $S_G=\{u,u',u''\}$ and
$S_H=\{v,v',v''\}$. Then
$$
d_{G\circ H}(S)=\left\{
\begin{array}{ll}
d_{H}(S_H), &\mbox {\rm if} \ u=u'=u''~and~deg_{G}(u)=0;\\[0.2cm]
\min \{d_{H}(S_H),3\}, &\mbox {\rm if} \ u=u'=u''~and~deg_{G}(u)\neq 0;\\[0.2cm]
\infty, &\mbox {\rm if} \ u\neq u',~u'=u''~and~d_G(u,u')=\infty;\\[0.2cm]
d_{G}(u,u')+1, &\mbox {\rm if} \ u\neq u',~u'=u''~and~d_G(u,u')\neq \infty;\\[0.2cm]
d_{G}(S_G), &\mbox {\rm if} \ u\neq u',~u\neq u''~and~u'\neq u''.
\end{array}
\right.
$$
\end{pro}

\subsection{Steiner distance and converxity}

In \cite{CaceresMarquezPuertas}, C\'{a}ceres, M\'{a}rquez, and Puertas formulated the following relation between three convex hulls in any graph.
\begin{pro}{\upshape \cite{CaceresMarquezPuertas}}\label{pro2-6}
Let $G$ be a connected graph and let $S\subseteq V(G)$. Then the following chain of
inclusions holds: $CH_g(S)\subseteq CH_{St}(S)\subseteq CH_{m}(S)$.
\end{pro}

A graph is called \emph{house-hole-domino free ($HHD$-free)} (see \cite{Brandstadt}) if it contains no induced house, domino, or induced cycle $C_k$, $k\geq 5$.
\begin{thm}{\upshape \cite{CaceresMarquezPuertas}}\label{th2-9}
For any set of vertices $S$ of a connected $HHD$-free graph $G$, any Steiner tree $T$ of
$S$ is contained in the geodesical convex hull of $S$.
\end{thm}

\begin{thm}{\upshape \cite{CaceresMarquezPuertas}}\label{th2-10}
Let $G$ be a connected $HHD$-free graph and let $S\subseteq V(G)$. Then $S$ is a $g$-convex set if and only if it is a St-convex set.
\end{thm}

A convexity is said to be a \emph{convex geometry} (see \cite{Van de Vel}) if every convex set is the convex hull of its extreme points. This property gives good behavior to a convexity in graphs, because in this case we can keep all information about a convex vertex sets just in its extreme points.  To
find conditions under which St-convexity shares this property, we firstly need to characterize
extreme points of a St-convex set. Recall that a vertex is called \emph{simplicial} if its neighborhood is
a complete subgraph. The next lemma shows that St-extreme points are simplicial vertices, the
same condition as in case of $g$-convex sets.

Note that, in any case, the extreme points of the convex hull of a set of vertices A belong to
$A$, because if $p\in CH(A)\setminus A$ is an extreme point; then $CH(A)\setminus \{p\}$ is a convex set containing $A$, which contradicts the minimality of convex hull.

C\'{a}ceres, M\'{a}rquez, and Puertas \cite{CaceresMarquezPuertas} characterized the class of graphs in which St-convexity becomes a convex geometry.
\begin{thm}{\upshape \cite{CaceresMarquezPuertas}}\label{th2-11}
The St-convexity in a connected graph $G$ is a convex geometry if and only if $G$ is
chordal and contains no induced $3$-fan.
\end{thm}

Let $G$ be a connected graph of order $n\geq 2$ and $k$ be an
integer with $2\leq k\leq n$. Let $S\subseteq V(G)$ and $v\in S$, the $k$-eccentricity $e_{k,S}(v)$ de $v$ en $S$ is defined
by $e_{k,S}(v)=\max\{d_S(K): K\subseteq S, |K|=k, v\in K\}$. In case $S=V(G)$, we denote $e_{k,S}$ simply
by $e_k$.
\begin{pro}{\upshape \cite{CaceresMarquezPuertas}}\label{pro2-10}
Let $G$ be a connected graph and let $uv$ be an edge of $G$, then $e_k(u)-1\leq
e_k(v)\leq e_k(u)+1$.
\end{pro}

Let $G$ be a connected graph and $S\subseteq V(G)$. A vertex $v\in S$ is called \emph{$k$-contour}
of $S$ if it satisfies $e_{k,S}(v)\geq e_{k,S}(u)$, for any $u\in N_S[v]$. The set $Ct_{k,S}(G)$ of $k$-contour vertices of $S$ is
called the \emph{$k$-contour set} of $S$.

The $k$-contour set is an enlargement of St-extreme
point set.
\begin{pro}{\upshape \cite{CaceresMarquezPuertas}}\label{pro2-10}
Let $G$ be a connected graph and $S\subseteq V(G)$. Then $Ct_{k,S}(G)$ contains all St-extreme points of $S$.
\end{pro}

C\'{a}ceres, M\'{a}rquez, and Puertas \cite{CaceresMarquezPuertas} proved the main result for $k$-contour vertices: these vertices can rebuild any St-convex set by means of a Steiner convex hull operation. This result provides, in some sense,
a generalization of Theorem \ref{th2-11}, using a vertex set bigger than St-extreme points, that works in
any connected graph.
\begin{thm}{\upshape \cite{CaceresMarquezPuertas}}\label{th2-12}
Let $G$ be a connected graph and $S\subseteq V(G)$ a St-convex set. Then
$CH_{St}(Ct_{k,S}(S))=S$.
\end{thm}

\section{Steiner Diameter}

For a graph $G$ of order $n\geq 2$, the Steiner diameter sequence of $G$ is defined as the
sequence
$$
{\rm sdiam}_2(G),{\rm sdiam}_3(G),\ldots,{\rm sdiam}_n(G),
$$
while the Steiner radius sequence is the sequence
$$
{\rm srad}_2(G),{\rm srad}_3(G),\ldots,{\rm srad}_n(G).
$$

\subsection{Steiner diameter of some graph classes}

The following results are immediate.
\begin{pro}{\upshape \cite{Mao2, Mao3}}\label{pro3-1}
Let $k,n$ be two integers with $2\leq k\leq n$.

$(1)$ For a complete graph $K_n$, ${\rm sdiam}_k(K_n)=k-1$;

$(2)$ For a path $P_n$, ${\rm sdiam}_k(P_n)=n-1$;

$(3)$ For a cycle $C_n$, ${\rm sdiam}_k(C_n)=\big
\lfloor\frac{n(k-1)}{k}\big \rfloor$.

$(4)$ For complete $r$-partite graph $K_{n_1,n_2,\ldots,n_r} \ (n_1\leq n_2\leq \cdots \leq n_r)$,
$$
{\rm sdiam}_k(K_{n_1,n_2,\ldots,n_r})=\left\{
\begin{array}{ll}
k-1,&\mbox{{\rm if}~$k>n_r$;}\\
k,&\mbox{{\rm if}~$k\leq n_r$.}
\end{array}
\right.
$$
\end{pro}

For Steiner diameter of threshold graphs, Wang, Mao, Cheng, and Melekian
\cite{WangMaoChengMelekian} derived the following results.
\begin{pro}{\upshape \cite{WangMaoChengMelekian}}\label{pro3-2}
Let $k,n$ be two integers with $3\leq k\leq n$, and let $G$ be a
threshold graph of order $n$. Let $i$ be the subscript of vertices
in $V(C_r)$ such that $g_ig_{n-r}'\in E(G)$ but
$g_{i+1}g_{n-r}'\notin E(G)$.

$(i)$ If $3\leq k\leq n-i$, then $siam_k(G)=k$.

$(i)$ If $n-i+1\leq k\leq n$, then $siam_k(G)=k-1$.
\end{pro}

Chartrand, Oellermann, Tian, and Zou \cite{ChartrandOellermannTianZou} established a relation between the Steiner $k$-diameter and Steiner $(k-1)$-diameter of
a tree, where $k\geq 3$ is an integer.
\begin{thm}{\upshape \cite{ChartrandOellermannTianZou}}\label{th3-1}
Let $k\geq 3$ be an integer and $T$ a tree of order $k\leq n$. Then
$$
{\rm sdiam}_{k-1}(T)\leq {\rm sdiam}_k(T)\leq \frac{k}{k-1}{\rm sdiam}_{k-1}(T).
$$
\end{thm}

The following  proposition will aid us in deriving a relation between the Steiner $k$-diameter and Steiner $k$-radius of a tree.
\begin{pro}{\upshape \cite{ChartrandOellermannTianZou}}\label{pro3-3}
Let $k\geq 3$ be an integer and $T$ a tree of order $k\leq n$. Then
$$
{\rm sdiam}_{k-1}(T)={\rm srad}_{k}(T).
$$
\end{pro}

\begin{cor}{\upshape \cite{ChartrandOellermannTianZou}}\label{cor3-1}
Let $k\geq 2$ be an integer and $T$ a tree of order $k\leq n$. Then
$$
{\rm srad}_{k}(T)\leq {\rm sdiam}_k(T)\leq \frac{k}{k-1}{\rm srad}_{k}(T).
$$
\end{cor}

Chartrand, Oellermann, Tian, and Zou \cite{ChartrandOellermannTianZou} conjectured that Corollary \ref{cor3-1} can be extended to any connected graph.
\begin{con}{\upshape \cite{ChartrandOellermannTianZou}}\label{con3-1}
Let $k\geq 2$ be an integer and $G$ is a connected graph of order $k\leq n$. Then
$$
{\rm srad}_{k}(G)\leq {\rm sdiam}_k(G)\leq \frac{k}{k-1}{\rm srad}_{k}(G).
$$
\end{con}

Chartrand, Oellermann, Tian, and Zou \cite{ChartrandOellermannTianZou} presented the desired characterization of diameter sequences
of trees.
\begin{thm}{\upshape \cite{ChartrandOellermannTianZou}}\label{th3-2}
A sequence $D_2,D_3,\ldots,D_n$ of positive integers is the diameter sequence
of a tree of order $n$ having $r$ end-vertices if and only if

$(1)$ $D_{k-1}<D_k\leq \frac{k}{k-1}D_{k-1}$ for $3\leq k\leq r$;

$(2)$ $D_k=k-1$  for $r\leq k\leq n$;

$(3)$ $D_{k+1}-D_k\leq D_{k}-D_{k-1}$ for $3\leq k\leq n-1$.
\end{thm}

\begin{cor}{\upshape \cite{ChartrandOellermannTianZou}}\label{cor3-2}
A sequence $R_2,R_3,\ldots,R_n$ of positive integers is the radius sequence
of a tree of order $n\geq 2$ having $r$ end-vertices if and only if

$(1)$ $R_{k-1}<R_k\leq \frac{k}{k-1}R_{k-1}$ for $3\leq k\leq r+1$;

$(2)$ $R_k=n-1$  for $r+1\leq k\leq n$;

$(3)$ $R_{k+1}-R_k\leq R_{k}-R_{k-1}$ for $4\leq k\leq n$.
\end{cor}

The maximum number of vertices of maximal
planar graphs of given diameter and maximum degree has been determined. Hell and Seyffarth \cite{HellSeyffarth} have shown that the maximum number of vertices in a planar graph with diameter $2$ and maximum degree $\Delta\geq 8$ is $\lfloor\frac{3}{2}\Delta+1\rfloor$. It was shown
in \cite{Seyffarth} that maximal planar graphs of diameter $2$ and maximum degree $\Delta\geq 8$ have no more than $\frac{3}{2}\Delta+1$ vertices. It was
also shown that there exist maximal planar graphs with diameter two and exactly $\lfloor\frac{3}{2}\Delta+1\rfloor$ vertices. Yang, Lin, and Dai \cite{YangLinDai} have computed the exact maximum number of vertices in planar graphs and maximal planar graphs with diameter two and maximum degree $\Delta$, for $\Delta<8$.

Fulek, Mori\'{c} and Pritchard \cite{FulekMoricPritchard} derived the following upper bound
for diameter.
\begin{thm}{\upshape \cite{FulekMoricPritchard}}\label{th3-3}
For every connected planar graph $G$ of order $n$ and size $m$,
$$
diam(G)\leq \frac{4(n-1)-m}{3}\,.
$$
\end{thm}
Since for $3$, $4$ and $5$-connected maximal planar graphs $m=3n-6$, the bound in Theorem \ref{th3-3} becomes $diam(G)\leq \frac{n+2}{3}$. It is well known that for the ordinary diameter, i.e., for the case $k=2$, if $G$ is a $\ell$-connected graph of order $n$, then
$$
diam(G)\leq \frac{n+\ell-2}{\ell},
$$
which yields
\begin{itemize}
\item[] $diam(G)\leq \frac{n+1}{3}$ for $3$-connected graphs $G$;

\item[] $diam(G)\leq \frac{n+2}{4}$ for $4$-connected graphs $G$;

\item[] $diam(G)\leq \frac{n+3}{5}$ for $5$-connected graphs $G$.
\end{itemize}
So the ordinary diameters of $3$, $4$ and $5$-connected maximal planar graphs do not exceed $\frac{n+1}{3}$, $\frac{n+2}{4}$ and $\frac{n+3}{5}$, respectively.

Ali, Mukwembi, and Dankelmann \cite{AliMukwembiDankelmann} derived the following upper bounds for Steiner $k$-diameter of maximal planar graphs.
\begin{thm}{\upshape \cite{AliMukwembiDankelmann}}\label{th3-4}
$(1)$ Let $G$ be a $3$-connected maximal planar graph of order $n$. If $k$ is an integer with $2\leq k\leq n$, then
$$
sdiam_k(G)\leq \frac{n}{3}+\frac{8k}{3}-5.
$$

$(2)$ Let $G$ be a $4$-connected maximal planar graph of order $n$. If $k$ is an integer with $2\leq k\leq n$, then
$$
sdiam_k(G)\leq \frac{n}{4}+\frac{19k}{4}-9.
$$

$(3)$ Let $G$ be a $5$-connected maximal planar graph of order $n$. If $k$ is an integer with $2\leq k\leq n$, then
$$
sdiam_k(G)\leq \frac{n}{5}+\frac{24k}{5}-9.
$$
\end{thm}

The following example (see Figure 3.1 $(a)$ for an illustration) shows that, for constant $k$, the bound in $(1)$ of Theorem \ref{th2-4} is best possible,
apart from the value of the additive constant.

\noindent {\bf Example 3.1.} {\upshape \cite{AliMukwembiDankelmann}} For an integer $\ell\geq \lceil \frac{n}{3}\rceil$
let $G_1,G_2,\ldots,G_{\ell}$ be disjoint copies of the cycle $C_3$, and
let $a_i, b_i, c_i\in V(G_i)$. Let $G'_{\ell}$ be the graph obtained from the union of $G_1,G_2,\ldots,G_{\ell}$ by adding the edges $a_{i+1}a_i, b_{i+1}b_i, c_{i+1}w_i,a_{i+1}b_i, c_{i+1}b_i, a_{i+1}c_i$ for $i=1,2,\ldots, \ell-1$. Clearly, $|V(G'_{\ell})|=3\ell$ so that $\ell =\frac{|V(G'_{\ell})|}{3}$. Clearly, $diam(G')=\ell-1 =
\frac{|V(G'_{\ell})|}{3}-1$ and
$sdiam_3(G')=d(\{a_1,b_1,c_{\ell}\})=\ell=\frac{|V(G'_{\ell})|}{3}$. For $k\geq 4$, let $k = 3q+r$ with $r\in \{1,2,3\}$ and define $S$ to be the set of $k$ vertices
$\{a_1,b_1,c_1,a_2,b_2,c_2,\ldots,a_q,b_q,c_q\}\cup R$, where $R\subseteq \{u_{\ell}, v_{\ell}, w_{\ell}\}$ is a set with $|R| = r$. It is easy to see that $d(S)=\ell-1+2q +
r-1=\frac{|V(G'_{\ell})|}{3}+ 2\lceil \frac{n}{3}\rceil+r-4\geq \frac{|V(G'_{\ell})|}{3}+ 2\lceil \frac{n}{3}\rceil-3$. Hence $sdiam_k(G'_{\ell})\geq \frac{|V(G'_{\ell})|}{3}+ 2\lceil \frac{n}{3}\rceil-3$ for $k\geq 2$.

The following example (see Figure 3.1 $(b)$ for an illustration) shows that, for constant $k$, the bound in $(2)$ of Theorem \ref{th2-4} is best possible,
apart from the value of the additive constant.

\noindent {\bf Example 3.2.} {\upshape \cite{AliMukwembiDankelmann}} For an integer $\ell\geq k$, let $G_1,G_2,\ldots,G_{\ell}$ be disjoint copies of the $4$-cycle $C_4$, and
let $a_i,b_i,c_i,d_i\in V(G_i)$. Let $G''_{\ell}$ be the graph obtained from the union of $G_1,G_2,\ldots,G_{\ell}$ by adding the edges $a_{i+1}a_i, b_{i+1}b_i,
c_{i+1}c_i, d_{i+1}d_i, a_{i+1}d_i,b_{i+1}a_i,c_{i+1}b_i, d_{i+1}c_i$, for $i=1,2,$ $\ldots,k-1$, $a_{\ell}c_{\ell}$ and $a_1c_1$. Clearly, $|V(G''_{\ell})|=4\ell$ so that $\ell= \frac{|V(G''_{\ell})|}{4}$. Clearly,
$diam(G''_{\ell})=\frac{|V(G''_{\ell})|}{4}-1$.
For $k\geq 3$, let $k=2q+r$ with $r\in \{1,2\}$ and define the set $S$ of $k$ vertices as $\{b_1, d_1, b_3, d_3, b_5,
d_5,\ldots, b_{2q-1}, d_{2q-1}\}\cup R$, where $R\subseteq \{b_{\ell}, d_{\ell}\}$ is a set with $|R|=r$. It is easy to verify that $d(S)=\ell-1+2q+2(r-1) =
\ell-3+k+r\geq \ell-2+k$. Hence we have $sdiam_k(G''_{\ell})\geq
\frac{|V(G''_{\ell})|}{4}+k-2$ for $k\geq 3$.

\begin{figure}[!hbpt]
\begin{center}
\includegraphics[scale=0.9]{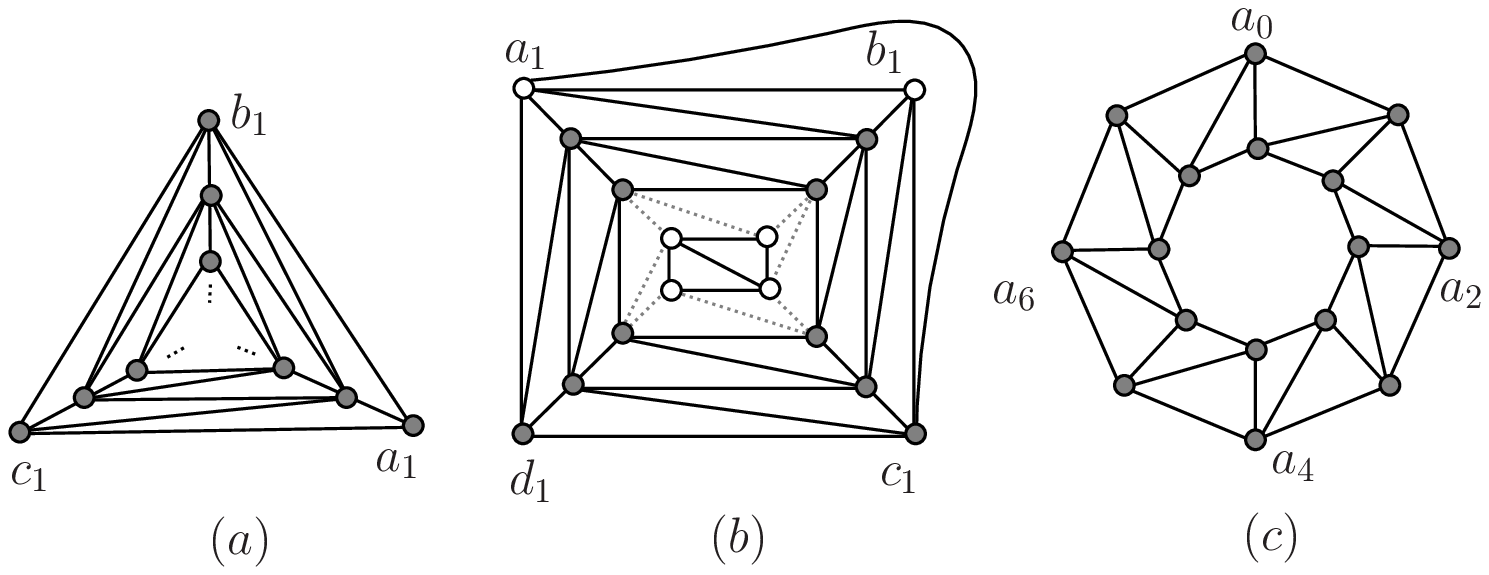}
\end{center}
\begin{center}
Figure 3.1: $(a)$ The graph $G'_{\ell}$; $(b)$ The graph $G''_{\ell}$; $(c)$ A $4$-connected planar graph for which $(1),(2)$ of Theorem \ref{th3-4} do not hold.
\end{center}
\end{figure}\label{fig4-1}

It is essential in $(1),(2)$ of Theorem \ref{th3-4} that $G$ is maximal planar and not just planar, as the following example demonstrates.
(For an illustration see Figure 3.1 $(c)$)

\noindent {\bf Example 3.3.} {\upshape \cite{AliMukwembiDankelmann}} Let $H$ be the Cartesian product of $K_2$ and a cycle $C_{\frac{n}{2}}$, where $n$ is even, i.e. let $V(H)=\{a_0,a_1,\ldots,
a_{\frac{n}{2}-1},b_0,b_1,\ldots,b_{\frac{n}{2}-1}\}$ and $E(H)=\{a_ia_{i+1}\,|\,i=0, 1,\ldots,\frac{n}{2}-1\}\cup \{b_ib_{i+1}\,|\,i=0,1,\ldots,\frac{n}{2}-1\}\cup \{a_ib_i\,|\,i= 0, 1,\ldots,\frac{n}{2}-1\}$
where subscripts are taken modulo $\frac{n}{2}$. Now let $H'$ be the planar graph obtained from $H$ by adding the edges $a_{i+1}b_i$  for $i =
0,1,\ldots,\frac{n}{2}-1$. If $k$ divides $\frac{n}{2}$, then the set $S_k=\{a_i\,|\,i\in \{0,
\frac{n}{2k}, 2\frac{n}{2k}, 3\frac{n}{2k}, (k-1)\frac{n}{2k}\}\}$ has $k$ vertices and $d(S_k)=\frac{k-1}{k}\frac{n}{2}$.
Hence $diam_k(H')\geq \frac{k-1}{2k}n$ which for constant $k\geq 4$ and large $n$ is greater than both $\frac{n}{3}+\frac{8k}{3}-5$ and $\frac{n}{4}+\frac{19k}{4}-9$. Since $H'$
is planar and $4$-connected, this shows that $(1),(2)$ of Theorem \ref{th3-4} do not hold for $4$-connected planar graphs.

The following graphs (see Figure 3.2 $(a)$ for an illustration) show that, for constant $k$, the bound in $(3)$ of Theorem \ref{th3-4} is best possible, apart from the value of the additive constant.

\noindent {\bf Example 3.4.} {\upshape \cite{AliMukwembiDankelmann}} For an integer $\ell \geq k$, let $G_1,G_2,\ldots,G_{\ell}$ be disjoint copies of the $5$-cycle, $C_5$, and
let $a_i,b_i,c_i,d_i,w_i\in V(G_i)$. Let $G'''_{\ell}$ be the graph obtained from the union of $G_1,G_2,\ldots,G_{\ell}$ by adding the edges $a_{i+1}a_i, b_{i+1}b_i,
c_{i+1}c_i, d_{i+1}d_i, w_{i+1}w_i, a_{i+1}w_i, $ $b_{i+1}a_i, c_{i+1}b_i, d_{i+1}c_i, w_{i+1}d_i$ for $i=1,2,\ldots,\ell-1$ and new vertices $v_{\ell}$ adjacent to $a_{\ell}, b_{\ell}, c_{\ell},$ $ d_{\ell}, w_{\ell}$ and $v_1$ adjacent to $a_1,b_1,c_1,d_1,w_1$. Clearly, $|V(G'''_{\ell})|=5\ell+2$ so that $\ell=\frac{|V(G'''_{\ell})|-2}{5}$. Now $diam(G'''_{\ell}) =\ell+1=\frac{|V(G'''_{\ell})|-2}{5}+1$. For $k\geq 3$ let $k=2q+r$ with $r\in \{0,1\}$ and define the set $S$ of $k$ vertices as $\{v_1,v_{\ell}\}\cup \{a_2,c_2,a_4,c_4,a_6,c_6,\ldots,
a_{2(q-1)}, c_{2(q-1)}\}\cup R$, where $R\subseteq \{a_{2q}\}$ is a set with $|R|=r$. It is easy to verify that $d(S)=\ell -1+2(q-1)\geq \ell -1+2q+r-3=k+\ell-4$. Hence we have $sdiam_k(G'''_{\ell})\geq \frac{|V(G'''_{\ell})|-2}{5}+\ell-4$ for $k\geq 3$.

The following example shows that in $(3)$ of Theorem \ref{th3-4} it is essential that $G$ is maximal planar, and not just planar. (For an illustration see Figure 3.2 $(b)$)

\noindent {\bf Example 3.5.} {\upshape \cite{AliMukwembiDankelmann}} For $n$ a multiple of $4k$, let $H''$ be the graph of order $n$ obtained from the disjoint union of two cycles of length $n/4$ with vertices $a_0, a_1,\ldots, a_{n/4-1}$ and $c_0,c_1,\ldots,c_{n/4-1}$, respectively, and a cycle of length $n/2$ with vertices $b_0, b_1,\ldots,b_{n/2-1}$, by adding edges $a_ib_{2i-1}, a_ib_{2i}, a_ib_{2i+1}, c_ib_{2i}, c_{i}b_{2i+1}, c_ib_{2i+2}$ for $i=0,1,\ldots,n/4$, where the subscripts are taken modulo $n/4$ for $a_j$ and $c_j$, and modulo $n/2$ for $b_j$. Then the set $S_k=\{a_i\,|\,i\in \{0,\frac{n}{4k}, 2\frac{n}{4k}, 3\frac{n}{4k}, (k-1)
\frac{n}{4k}\}\}$ has $k$ vertices and $d(S_k)=\frac{k-1}{k}\frac{n}{4}$. Hence $sdiam_k(H'') \geq \frac{k-1}{4k}n$ which for constant $k\geq 5$ and large $n$ is greater than $\frac{n}{5}
+\frac{24k}{5}-9$. Since $H''$ is planar and
$5$-connected, this shows that $(3)$ of Theorem \ref{th3-4} does not hold for $5$-connected planar graphs.
\begin{figure}[!hbpt]
\begin{center}
\includegraphics[scale=0.8]{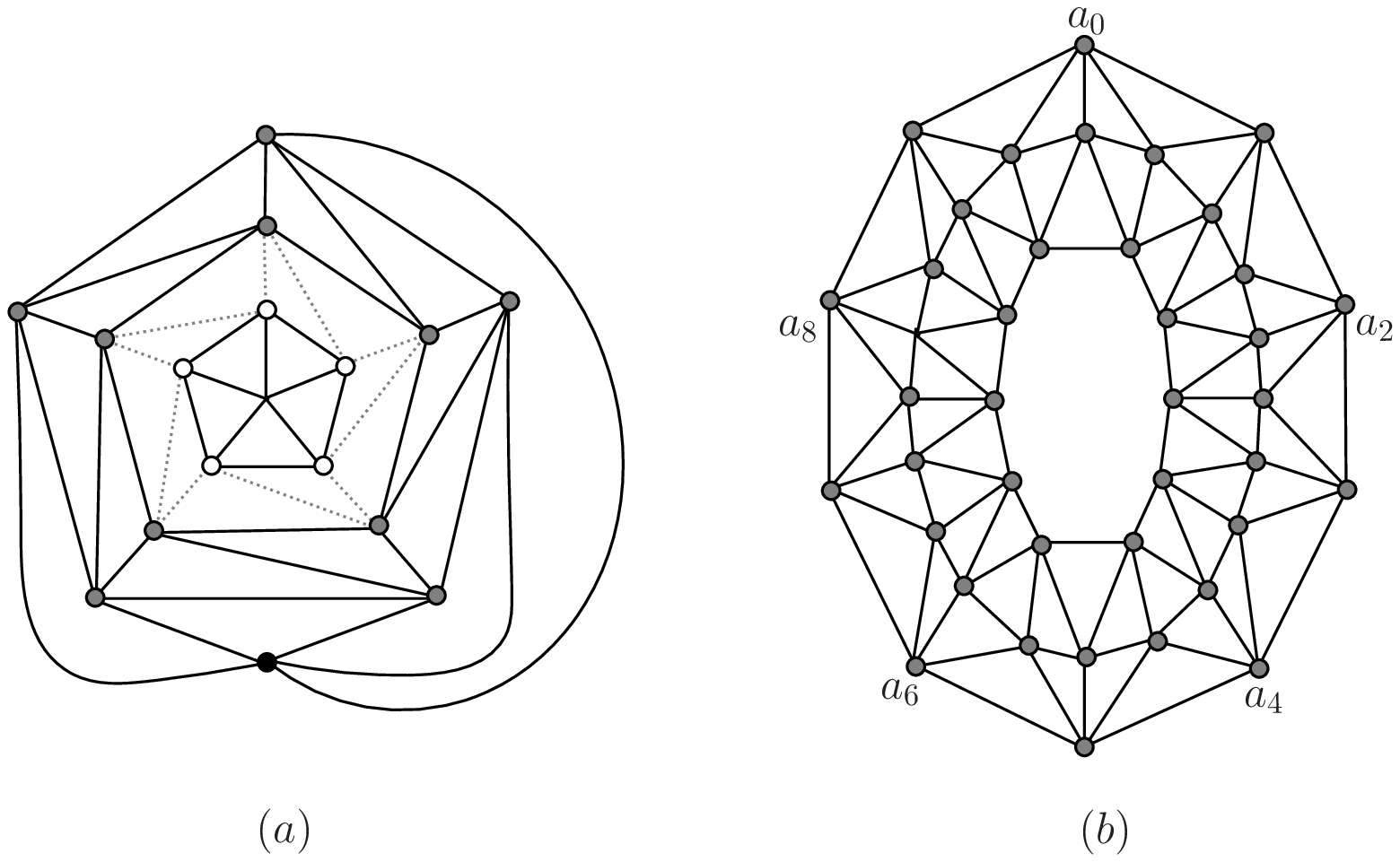}
\end{center}
\begin{center}
Figure 3.2: $(a)$ The graph $G'''_{\ell}$; $(b)$  A $5$-connected planar graph for which $(3)$ of Theorem \ref{th3-4} does not hold.
\end{center}
\end{figure}\label{fig4-1}

\begin{rem}{\upshape \cite{AliMukwembiDankelmann}}\label{rem3-1}
All bounds on the Steiner $k$-diameter given in this paper are sharp except for an additive constant provided that $k$ is
constant. It would be interesting to know if for infinitely many values of $n$ and $k$ there are graphs that come within an additive
constant, neither dependent on $n$ nor on $k$, of our bounds, or if our bounds can be improved by a term that is linear in $k$.
\end{rem}

\subsection{Bounds for Steiner diameter}

In \cite{ChartrandOkamotoZhang}, Chartrand, Okamoto, and Zhang derived the upper
and lower bounds for $sdiam_k(G)$.
\begin{thm}{\upshape\cite{ChartrandOkamotoZhang}}\label{th3-5}
Let $k,n$ be two integers with $2\leq k\leq n$, and let $G$ be a
connected graph of order $n$. Then
$$
k-1\leq sdiam_k(G)\leq n-1.
$$
Moreover, the bounds are sharp.
\end{thm}

\subsubsection{In terms of order and minimum degree}

For the ordinary diameter, Erd\"{o}s, Pach, Pollack, and Tuza \cite{ErdosPachPollackTuza} give the following bounds in terms of order and minimum degree.
\begin{thm}{\upshape\cite{ErdosPachPollackTuza}}\label{th3-6}
$(1)$ For all connected graphs $G$,
$$
diam(G)\leq \frac{3n}{\delta+1}-1.
$$

$(2)$ If $G$ is a triangle-free graph, then
$$
diam(G)\leq 4\left\lceil \frac{n-\delta-1}{2\delta}\right\rceil.
$$

$(3)$ If $G$ is a $C_4$-free graph, then
$$
diam(G)\leq 4\left\lceil \frac{5n}{\delta^2-2[\delta/2]+1}\right\rceil.
$$
\end{thm}

The result for all connected graphs was extended by Dankelmann, Swart, and Oellermann in \cite{DankelmannSwartOellermann2}.
\begin{thm}{\upshape\cite{DankelmannSwartOellermann2}}\label{th3-7}
Let $G$ be a connected graph of order $n$ and minimum degree $\delta\geq 2$.
If $2\leq k\leq n$, then
$$
sdiam_k(G)\leq \frac{3n}{\delta+1}+3k.
$$
\end{thm}

Ali, Mukwembi, and Dankelmann \cite{AliDankelmannMukwembi} improved this upper bound and generalized the corresponding result by Erd\"{o}s, Pach, Pollack, and Tuza \cite{ErdosPachPollackTuza}.
\begin{thm}{\upshape\cite{AliDankelmannMukwembi}}\label{th3-8}
$(1)$ Let $G$ be a connected graph of order $n$ and minimum degree $\delta\geq 2$.
If $k$ is an integer with $2\leq k\leq n$, then
$$
sdiam_k(G)\leq \frac{3n}{\delta+1}+2k-5.
$$

$(2)$ Let G be a connected triangle-free graph of order $n$ and minimum degree $\delta\geq 2$.
If $k$ is an integer with $2\leq k\leq n$, then
$$
sdiam_k(G)\leq \frac{2n}{\delta}+3k-6.
$$

$(3)$ Let $G$ be a connected $C_4$-free graph of order $n$ and minimum degree $\delta\geq 2$. If $k$ is an integer with $2\leq k\leq n$, then
$$
sdiam_k(G)\leq \frac{5n}{\delta^2-2\lfloor \delta/2\rfloor+1}+4k-9.
$$
\end{thm}

The following corollary is derived in the same paper.
\begin{cor}{\upshape\cite{AliDankelmannMukwembi}}\label{cor3-3}
$(1)$ Let $G$ be a connected graph of order $n$ and minimum degree $\delta\geq 2$. If $k$ is an integer with $2\leq k\leq n$, then $G$ has a
spanning tree $T$ with
$$
sdiam_k(T)\leq \frac{3n}{\delta+1}+2k-3.
$$

$(2)$ Let $G$ be a connected triangle-free graph of order $n$ and minimum degree $\delta\geq 2$. If $k$ is an integer with $2\leq k\leq n$, then $G$ has a
spanning tree $T$ with
$$
sdiam_k(T)\leq \frac{2n}{\delta}+3k-3.
$$

$(3)$ Let $G$ be a connected $C_4$-free graph of order $n$ and minimum degree $\delta\geq 2$. If $k$ is an integer with $2\leq k\leq n$, then $G$ has a spanning tree $T$ with
$$
sdiam_k(T)\leq \frac{5n}{\delta^2-2\lfloor \delta/2\rfloor+1}+4k-5.
$$
\end{cor}

\subsubsection{In terms of the girth and minimum degree}

Using methods initiated by Dankelmann
and Entringer \cite{DankelmannEntringer} and methods used in
\cite{AliDankelmannMukwembi}, Ali \cite{Ali} proved in their paper that
\begin{thm}{\upshape\cite{Ali}}\label{th3-9}
Let $G$ be a connected graph of order $n$, minimum degree $\delta\geq 3$, and girth $g$.

$(1)$ If $k$ is an integer with $2\leq k\leq n$ and $g$ is odd, then
$$
sdiam_k(G)\leq g\frac{n}{K}+(g-1)k-2g+1,
$$
where $K=1+\delta\frac{(\delta-1)^{(g-1)/2}-1}{\delta-2}$.

$(2)$ If $k$ is an integer with $2\leq k\leq n$ and $g$ is even, then
$$
sdiam_k(G)\leq g\frac{n}{L}+(g-1)k-2g+2,
$$
where $L=2\delta\frac{(\delta-1)^{g/2}-1}{\delta-2}$.
\end{thm}

Ali \cite{Ali} showed that the bounds in Theorem \ref{th3-9} are asymptotically sharp, apart from an additive constant.

\noindent {\bf Example 3.6.} If $\delta$ and $g$ are such that there exists a Moore graph of minimum degree $\delta$ and girth $g$, i.e., a graph with minimum degree $\delta$, girth $g$, and order $1+\delta+\delta(\delta-1)+
\delta(\delta-1)2+\cdots+\delta(\delta-1)(g-3)/2$ (if $g$ is odd), or $2(1+(\delta-1)+(\delta-1)^2+(\delta-1)^3+\cdots+(\delta-1)(g-2)/2)$ (if $g$ is even).
For a given integer $\ell>0$, let $G_1,G_2,\cdots,G_{\ell}$ be disjoint copies of the $(\delta,g)$-Moore graph, and let $u_iv_i\in E(G_i)$. Let $G_{\ell,\delta,g}$ be the graph obtained from the union of $G_1,G_2,\ldots,G_{\ell}$ by deleting the edges $u_iv_i$ for $i =2,3,\ldots,\ell-1$ and adding the edges
$u_{i+1}v_i$ for $i=1,2,\ldots,\ell-1$.

$(1)$ If $g$ is odd, then $|V(G_{\ell,\delta,g})|=\ell K$, so $\ell=\frac{|V(G_{\ell,\delta,g})|}{K}$. If $2\leq n\leq 2\delta$ then, by a simple calculation, $diam_k(G_{\ell,\delta,g})=g\ell+k-5$, and so $sdiam_k(G_{\ell,\delta,g})=g\frac{|V(G_{\ell,\delta,g})|}{K}+k-5$. Note that the set $S$ of $k$ vertices of $G_{\ell,\delta,g}$ contains $\delta$ vertices of $G_1$ and $\delta$ vertices
of $G_{\ell}$. In this case, the difference between $diam_k(G_{\ell,\delta,g})$ and the bound in Theorem \ref{th3-9} is $(g-2)k-2g+6$.
For $k>2\delta$, we use the estimate $diam_k(G_{\ell,\delta,g})\geq diam(G_{\ell,\delta,g})=g\frac{|V(G_{\ell,\delta,g})|}{K}
-3$. In this case, the difference between the Steiner
$k$-diameter of $G_{\ell,\delta,g}$ and the bound in Theorem \ref{th3-9} is bounded by the additive constant $(g-1)k-2g+4$.

$(2)$ If $g$ is even, then $|V(G_{\ell,\delta,g})|=L\ell$, and thus $\ell=\frac{|V(G_{\ell,\delta,g})|}{L}$. If $k\leq 2\delta-2$ then, by a simple calculation, $sdiam_k(G_{\ell,\delta,g})=g\ell+k-7$, and thus $sdiam_k(G_{\ell,\delta,g})=g\frac{|V(G_{\ell,\delta,g})|}{L}
+k-7$. Note that the set $S$ of $k$ vertices of $G_{\ell,\delta,g}$ contains $\delta-1$
vertices of $G_1$ and $\delta-1$ vertices of $G_{\ell}$. In this case, the difference between $sdiam_k(G_{\ell,\delta,g})$ and the bound in Theorem \ref{th3-9} is $(g-2)k-2g+9$. For $k>2\delta-2$, we use the estimate $sdiam_k(G_{\ell,\delta,g})\geq diam(G_{\ell,\delta,g})=g \frac{|V(G_{\ell,\delta,g})|}{L}-5$. In this case, the
difference between the Steiner $k$-diameter of $G_{\ell,\delta,g}$ and the bound in Theorem \ref{th3-9} is bounded by the additive constant $(g-1)k-2g+7$.

A slight modification of the proof of Theorem \ref{th3-9} yields the following result.
\begin{cor}{\upshape\cite{Ali}}\label{cor3-4}
Let $G$ be a connected graph of order $n$, minimum degree $\delta\geq 3$, and girth $g$.

$(1)$ If $k$ is an integer with $2\leq k\leq n$ and $g$ is odd, then $G$ has a spanning tree $T$ with
$$
sdiam_k(G)\leq g\frac{n}{K}+(g-1)k-g,
$$
where $K=1+\delta\frac{(\delta-1)^{(g-1)/2}-1}{\delta-2}$.

$(2)$ If $k$ is an integer with $2\leq k\leq n$ and $g$ is even, then $G$ has a spanning tree $T$ with
$$
sdiam_k(G)\leq g\frac{n}{L}+(g-1)k-g+1,
$$
where $L=\delta\frac{(\delta-1)^{g/2}-1}{\delta-2}$.
\end{cor}

\subsubsection{In terms of Steiner radius}

It was conjectured in \cite{HenningOellermannSwart} that for all integers
$k\geq 2$ and every connected graph $G$ of order $n\geq k$,
\begin{equation}                    \label{eq3-5}
sdiam_k(G)\leq \frac{2(k+1)}{2k-1}srad_k(G)\,.
\end{equation}
In \cite{HenningOellermannSwart}, it was shown that for each integer $k\geq 2$
a graph $G$ exists for which equality is attached in $1$ and that (\ref{eq3-5})
is valid if $k\in \{2,3,4\}$, but further progress has been reported.

Ore \cite{Ore} showed that, for a connected graph $G$ of order $n$, size $m$ and diameter
$d$,
$$
e(G)\leq d+\frac{1}{2}(n-d-1)(n-d+4).
$$
Extensions of Ore's result were given be Caccetta and Smyth \cite{CaccettaSmyth, CaccettaSmyth2, CaccettaSmyth3}.

Ore's result is generalized by Dankelmann, Swart, and Oellermann \cite{DankelmannSwartOellermann2} in the next theorem, which can be used to obtain an upper bound
on the Steiner $k$-diameter of $G$ in terms of the order and size of $G$.
\begin{thm}{\upshape\cite{DankelmannSwartOellermann2}}\label{th3-10}
If $G$ is a connected graph of order $n$, size $m$ and $sdiam_k(G)=d_k$, where $2\leq k\leq n-1$, then
$$
e(G)\leq d_k+\binom{k-1}{2}+\binom{k-d_k-1}{2}+(k-d_k-1)(k+1).
$$
\end{thm}

An analysis of the above proof shows that, for $k\geq 4$ and $n\geq 2k+1$ the graph $G$
obtained from the union of three disjoint graphs $G_1,G_2,G_3$ with $G_1\cong K_{k+1}$,
$G_2\cong K_{k}$ and $G_3\cong K_{n-2k-1}$ by joining each vertex in $G_2$ to the same
end vertex of $G_3$ and to each vertex of $G_1$ is such that, for $G$, equality is attained
in $(1)$ of Theorem \ref{th3-10}.

Since the Steiner $k$-diameter of a tree is easy to compute, bounds for the Steiner $k$-diameter
of graphs in classes that do not contain trees are of interest. In the following theorem, they obtained
an upper bound on the Steiner $k$-diameter of a $2$-connected graph of order $n$.
\begin{thm}{\upshape\cite{DankelmannSwartOellermann2}}\label{th3-11}
Let $G$ be a $2$-connected graph of order $n$, and let $2\leq k\leq n$. Then
$$
sdiam_k(G)\leq \left\lfloor \frac{n(k-1)}{k}\right\rfloor
$$
and equality is attained for the cycle $C_n$.
\end{thm}

Theorem \ref{th3-11} shows that among all $2$-connected graphs of given order $n$,
the cycle $C_n$ has the largest possible Steiner $k$-diameter. Dankelmann, Swart, and Oellermann \cite{DankelmannSwartOellermann2} conjectured that more
generally
\begin{con}{\upshape \cite{DankelmannSwartOellermann2}}\label{con3-2}
For a $2k$-connected graph $G$ of order $n$,
$$
sdiam_k(G)\leq sdiam_k(C_{n}^{2k}),
$$
where $C_{n}^{2k}$ denotes the $2k$-th power of the cycle $C_n$.
\end{con}

\subsection{Graphs with given Steiner diameter}

The following observation is immediate.
\begin{obs}\label{obs3-1}
Let $G$ be a connected graph of order $n$. Then

$(1)$ $diam(G)=1$ if and only if $G$ is a complete graph;

$(2)$ $diam(G)=n-1$ if and only if $G$ is a path of order $n$.
\end{obs}

Let $uv$ be an edge in $G$. A \emph{double-star} on $uv$ is a
maximal tree in $G$ which is the union of stars centered at $u$ or
$v$ such that each star contains the edge $uv$. Bloom \cite{Bloom}
characterized the graphs with $diam(G)=2$.
\begin{thm}{\upshape\cite{Bloom}}\label{th3-12}
Let $G$ be a connected graph of order $n$. Then $diam(G)=2$ if
and only if $\overline{G}$ is non-empty and $\overline{G}$ does not
contain a double star of order $n$ as its subgraph.
\end{thm}

In \cite{MaoMelekianCheng}, Mao, Melekian, and Cheng derived the following
result.
\begin{thm}\label{th3-13}
Let $\ell,n$ be two integers with $1\leq \ell\leq n-2$, and let $G$
be a graph of order $n$. Then $\kappa(G)\geq \ell$ if and only if
$sdiam_{n-\ell+1}(G)=n-\ell$.
\end{thm}

In \cite{WangMaoLiYe}, Wang, Mao, Li, and Ye obtained the structural properties of
graphs with $sdiam_k(G)=n-1$.
\begin{thm}{\upshape \cite{WangMaoLiYe}}\label{th3-14}
Let $k,n$ be two integers with $3\leq k\leq n-1$. Let $G$ be a
connected graph of order $n$. Then $sdiam_k(G)=n-1$ if and only if
the number of non-cut vertices in $G$ is at most $k$.
\end{thm}

\subsubsection{Results for small $k$}

In \cite{Mao}, Mao characterized the
graphs with $sdiam_3(G)=2$, which can be seen as an
extension of $(1)$ of Observation \ref{obs3-1}.
\begin{thm}{\upshape \cite{Mao}}\label{th3-15}
Let $G$ be a connected graph of order $n$. Then $sdiam_3(G)=2$ if
and only if $0\leq \Delta(\overline{G})\leq 1$ if and only if
$n-2\leq \delta(G)\leq n-1$.
\end{thm}

A \emph{triple-star} $H_1$ is
defined as a connected graph of order $n$ obtained from a triangle
and three stars $K_{1,a}, K_{1,b}, K_{1,c}$ by identifying the
center of a star and one vertex of the triangle, where $0\leq a\leq
b\leq c$, $c\geq 1$ and $a+b+c=n-3$; see Figure 3.3 $(a)$.
Let $H_2$ be a connected graph of order $n$ obtained from a path
$P=uvw$ and $n-3$ vertices such that for each $x\in V(H_2)-
\{u,v,w\}$, $xu,xv,xw\in E(H_2)$, or $xu,xv\in E(H_2)$ but $xw\notin
E(H_2)$, or $xv,xw\in E(H_2)$ but $xu\notin E(H_2)$, or $xu,xw\in
E(H_2)$ but $xv\notin E(H_2)$, or $xv\in E(H_2)$ but $xu,xw\notin
E(H_2)$; see Figure 3.3 $(b)$.
\begin{figure}[!hbpt]
\begin{center}
\includegraphics[scale=0.9]{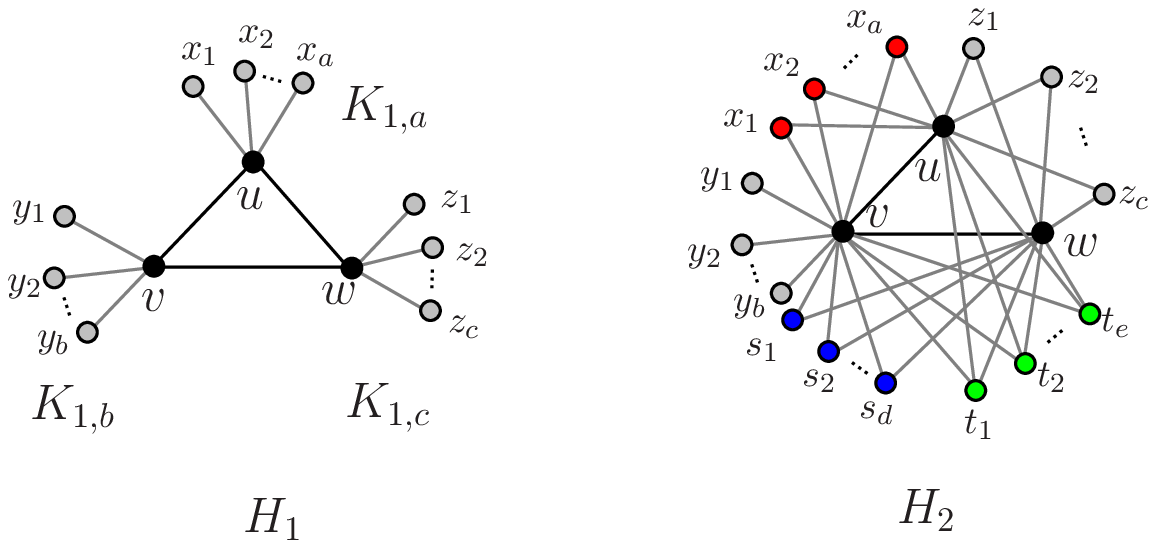}
\end{center}
\begin{center}
Figure 3.3: The graphs $H_1$ and $H_2$.
\end{center}
\end{figure}\label{fig4-1}

\begin{thm}{\upshape \cite{Mao}}\label{th3-16}
Let $G$ be a connected graph of order $n$. Then $sdiam_3(G)=3$ if
and only if $G$ satisfies the following conditions.

$\bullet$ $\Delta(\overline{G})\geq 2$;

$\bullet$ $\overline{G}$ does not contain a triple-star $H_1$ as its
subgraph;

$\bullet$ $\overline{G}$ does not contain $H_2$ as its subgraph.
\end{thm}

Denote by $T_{a,b,c}$ a tree with a vertex $v$ of degree $3$ such
that $T_{a,b,c}-v=P_{a}\cup P_{b}\cup P_{c}$, where $0\leq a\leq
b\leq c$ and $1\leq b\leq c$ and $a+b+c=n-1$; see Figure 3.4
$(a)$. Observe that $T_{0,b,c}$ where $b+c=n-1$ is a path of order
$n$. Denote by $\bigtriangleup_{p,q,r}$ a unicyclic graph containing
a triangle $K_{3}$ and satisfying $\bigtriangleup_{p,q,r}-
V(K_{3})=P_{p}\cup P_{q}\cup P_{r}$, where $0\leq p\leq q\leq r$ and
$p+q+r=n-3$; see Figure 3.4 $(b)$.
\begin{figure}[!hbpt]
\begin{center}
\includegraphics[scale=0.8]{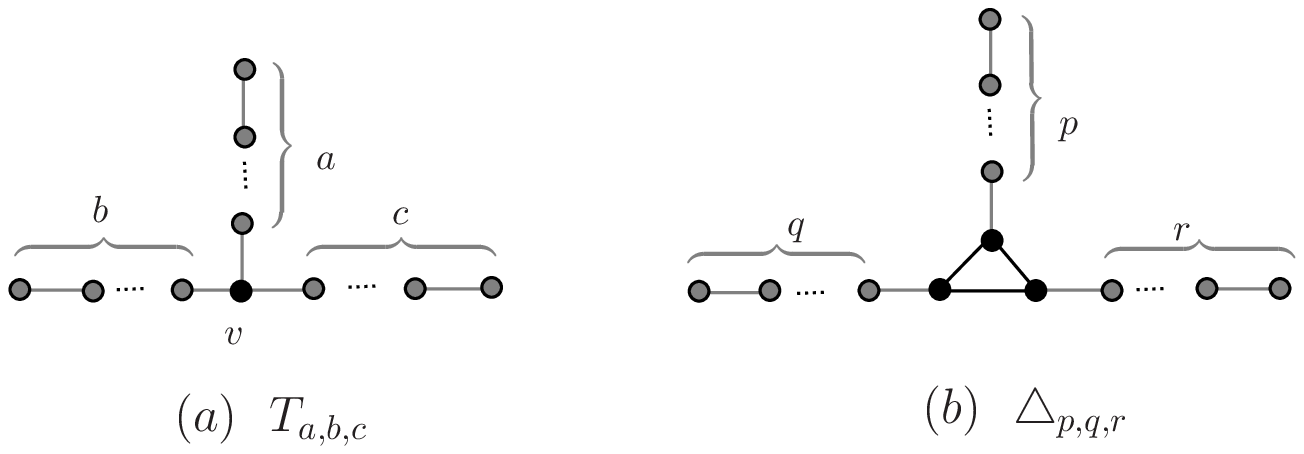}
\end{center}
\begin{center}
Figure 3.4: The graphs $T_{a,b,c}$ and $\bigtriangleup_{p,q,r}$.
\end{center}\label{fig4-2}
\end{figure}

\begin{thm}{\upshape \cite{Mao}}\label{th3-17}
Let $G$ be a connected graph of order $n \ (n\geq 3)$. Then
$sdiam_3(G)=n-1$ if and only if $G=T_{a,b,c}$ where $a\geq 0$ and
$1\leq b\leq c$ and $a+b+c=n-1$, or $G=\bigtriangleup_{p,q,r}$ where
$0\leq p\leq q\leq r$ and $p+q+r=n-3$.
\end{thm}

In \cite{WangMaoLiYe}, Wang, Mao, Li, and Ye characterized the graphs with $sdiam_4(G)=3,4,n-1$, respectively.
\begin{thm}{\upshape \cite{WangMaoLiYe}}\label{th3-18}
Let $G$ be a connected graph of order $n \ (n\geq 4)$.

$(i)$ If $n=4$, then $sdiam_4(G)=3$;

$(ii)$ If $n\geq 5$, then $sdiam_4(G)=3$ if
and only if $n-3\leq\delta(G)\leq n-1$ and $C_4$ is not a subgraph of $\overline{G}$.
\end{thm}

A graph $A_1$ is
defined as a connected graph of order $n \ (n\geq 5)$ obtained from a $K_4$ with vertex set $\{u_1,u_2,u_3,u_4\}$
and four stars $K_{1,a}, K_{1,b}, K_{1,c}, K_{1,d}$ by identifying the
center of one star and one vertex in $\{u_1,u_2,u_3,u_4\}$, where $0\leq a\leq
b\leq c\leq d$, $d\geq 1$, and $a+b+c+d=n-4$; see Figure 3.3.
\begin{figure}[!hbpt]
\begin{center}
\includegraphics[scale=0.8]{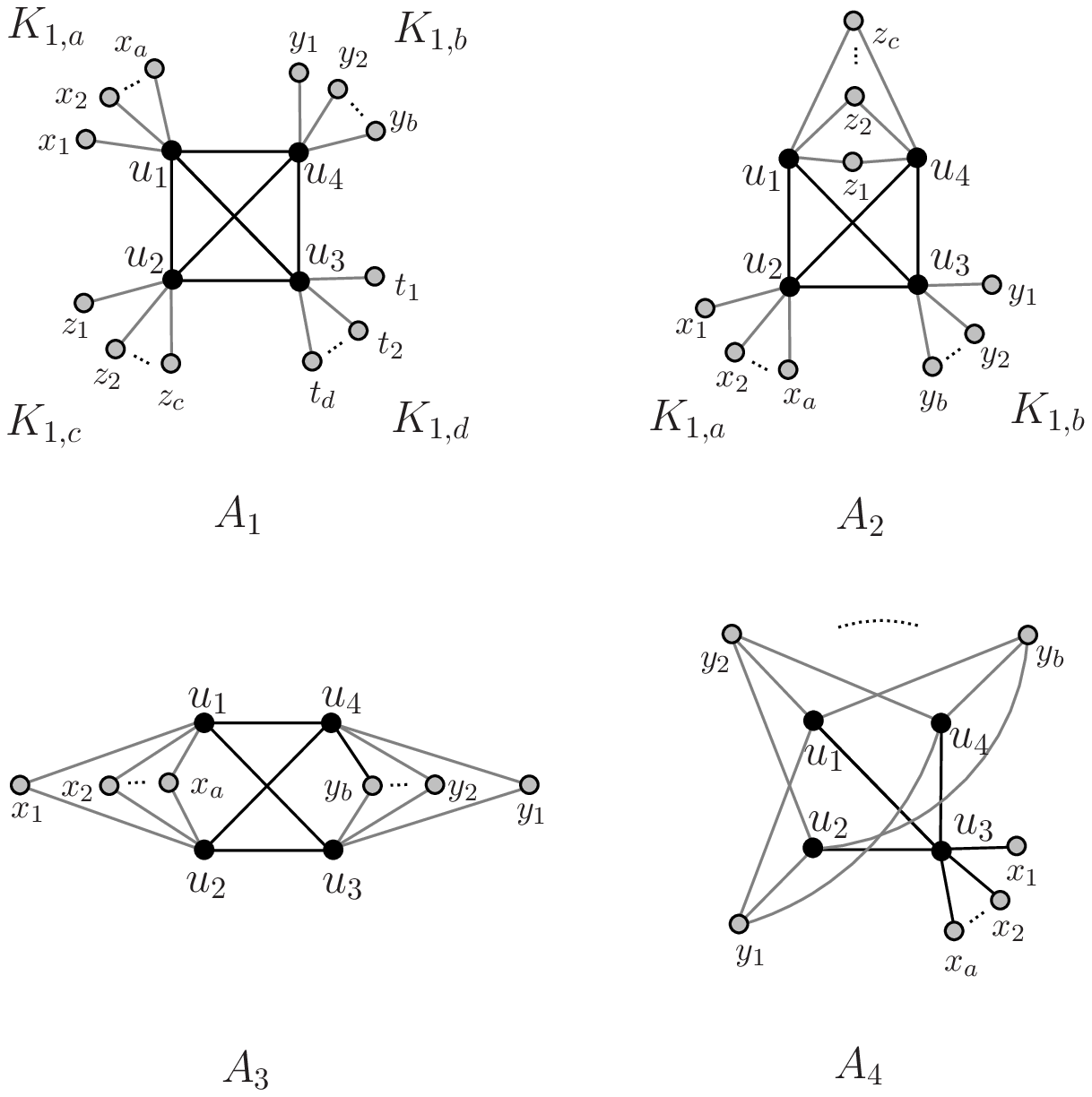}
\end{center}
\begin{center}
Figure 3.5: The graphs $A_i \ (1\leq i\leq 4)$.
\end{center}\label{fig4-3}
\end{figure}

A graph $A_2$ is
defined as a connected graph of order $n \ (n\geq 5)$ obtained from $K_4-e$ with vertex set
$\{u_1,u_2,u_3,u_4\}$, $e=u_1u_4$ and two stars $K_{1,a}$, $K_{1,b}$ by identifying the
center of a star and one vertex in $\{u_2,u_3\}$, and then adding the paths $u_1z_iu_4 \ (1\leq i\leq c)$, where $0\leq a\leq
b$, $b\geq 0$, $c\geq 0$ and $a+b+c=n-4$; see Figure 3.5.

A graph $A_3$ is
defined as a connected graph of order $n \ (n\geq 5)$ obtained from a cycle $C_4=u_1u_2u_3u_4u_1$ by adding the paths $u_1x_iu_2 \ (1\leq i\leq a)$ and the paths $u_3y_ju_4 \ (1\leq j\leq b)$, where $0\leq a\leq
b$, $b\geq 1$ and $a+b=n-4$; see Figure 3.5.

A graph $A_4$ is
defined as a connected graph of order $n \ (n\geq 5)$ obtained from a star $K_{1,3}$ with vertex set $\{u_1,u_2,u_3,u_4\}$
and a star $K_{1,a}$ by identifying $u_3$ and the center of $K_{1,a}$, where $u_3$ is the center of $K_{1,3}$, and then
adding the vertices $y_i$ and the edges $y_iu_j \ (1\leq i\leq b, \ j=1,2,4)$, where $0\leq a\leq
b$, $b\geq 1$ and $a+b=n-4$; see Figure 3.5.
\begin{thm}{\upshape \cite{WangMaoLiYe}}\label{th3-19}
Let $G$ be a connected graph of order $n \ (n\geq 5)$. Then $sdiam_4(G)=4$ if
and only if $G$ satisfies one of the following conditions.

$(i)$ $\delta(G)=n-3$ and $C_4$ is a subgraph of $\overline{G}$;

$(ii)$ $\delta(G)\leq n-4$ and each $A_i\ (1\leq i\leq 4)$ is not a spanning subgraph of $\overline{G}$ (see Figure 3.3).
\end{thm}

In \cite{WangMaoLiYe}, Wang, Mao, Li, and Ye also defined the following graph classes.
\begin{itemize}
\item Let $T_{a,b,c,d} \ (0\leq a,b,c,d\leq n-1, a+b+c+d\leq n-1)$ be a tree of order $n  \ (n\geq 5)$ obtained from three paths
$P_1,P_2,P_3$ of length $n-b-c-1,b,c$ respectively by identifying
the $(a+1)$-th vertex of $P_1$ and one endvertex of $P_2$, and then
identifying the $(n-b-c-d)$-th vertex of $P_1$ and one endvertex of
$P_3$ (Note that $u$ and $v$ can be the same vertex);
\begin{figure}[!hbpt]
\begin{center}
\includegraphics[scale=0.8]{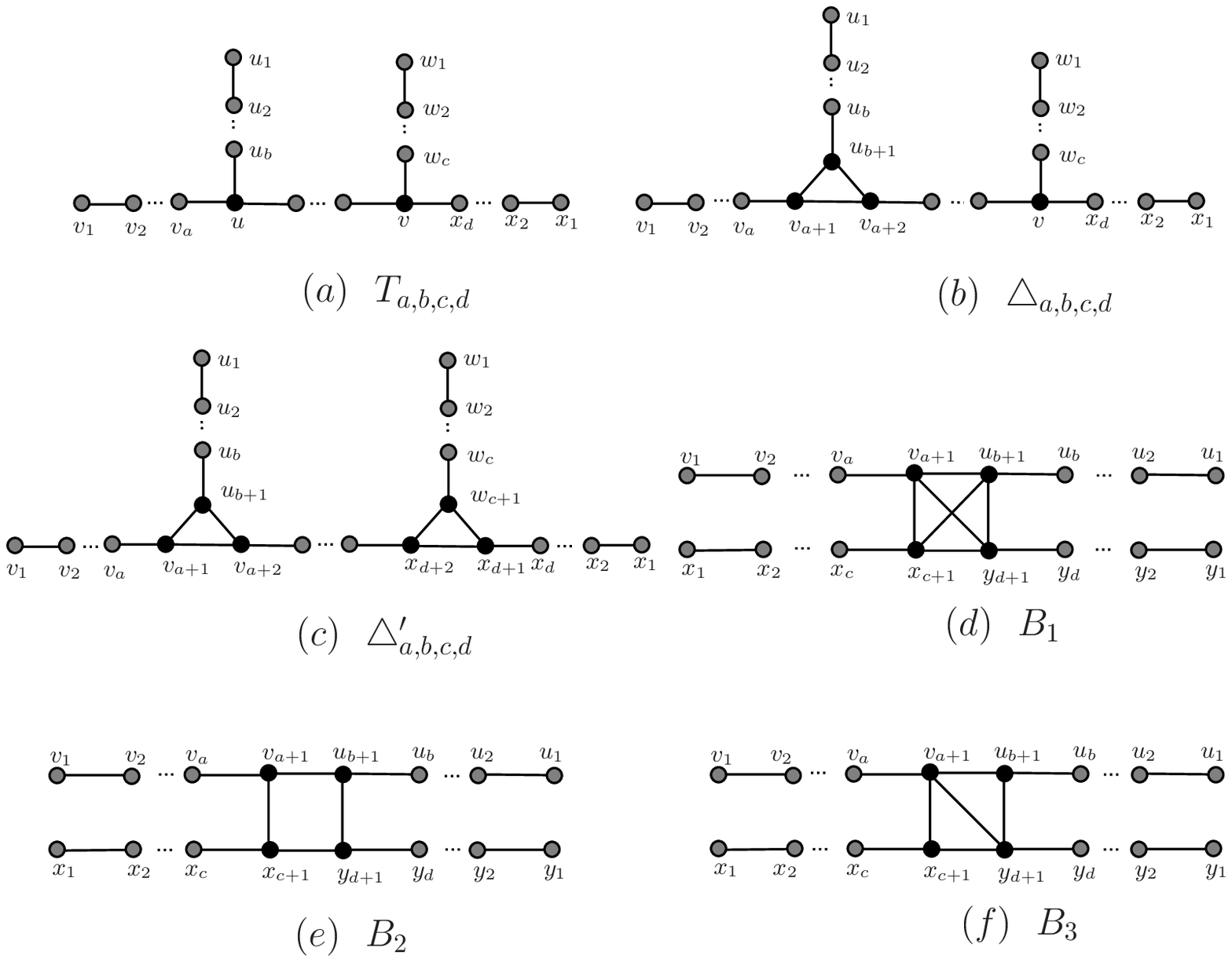}\\[0.5cm]
Figure 3.6: The graphs $T_{a,b,c,d},\triangle_{a,b,c,d},\triangle_{a,b,c,d}'$ and $B_i \ (1\leq i\leq 3)$.
\end{center}\label{fig4-4}
\end{figure}

\item Let $\triangle_{a,b,c,d} \ (0\leq a,b,c,d\leq n-2, a+b+c+d\leq n-2)$ be an unicyclic graph of order $n  \ (n\geq 5)$ obtained from three paths
$P_1,P_2,P_3$ of length $n-b-c-1,b+1,c$ respectively by identifying
the $(a+1)$-th vertex of $P_1$ and one endvertex of $P_2$, and then
identifying the $(n-b-c-d)$-th vertex of $P_1$ and one endvertex of
$P_3$, and then adding an edge $u_{b+1}v_{a+2}$  (Note that $v_{a+2}$ and $v$ can be the same vertex).

\item Let $\triangle_{a,b,c,d}' \ (0\leq a,b,c,d\leq n-3, \ a+b+c+d\leq n-3)$ be an bicyclic graph of order $n  \ (n\geq 5)$ obtained from three paths
$P_1,P_2,P_3$ of length $n-b-c-1,b+1,c+1$ respectively by identifying
the $(a+1)$-th vertex of $P_1$ and one endvertex of $P_2$, and then
identifying the $(n-b-c-d)$-th vertex of $P_1$ and one endvertex of
$P_3$, and then adding two edges $u_{b+1}v_{a+2}$ and $w_{c+1}x_{d+2}$ (Note that $v_{a+2}$ and $v$ can be the same vertex).

\item Let $B_1$ be a graph of order $n  \ (n\geq 5)$ obtained from a clique of order $4$ and four paths
$P_1,P_2,P_3,P_4$ of length $a,b,c,d \ (0\leq a,b,c,d\leq n-4, \ a+b+c+d=n-4)$ respectively by identifying
each vertex of this clique with an endvertex of one of the four
paths.

\item Let $B_2$ be a graph of order $n  \ (n\geq 5)$ obtained from a cycle of order $4$ and four paths
$P_1,P_2,P_3,P_4$ of length $a,b,c,d \ (0\leq a,b,c,d\leq n-4, \ a+b+c+d=n-4)$ respectively by identifying
each vertex of this cycle with an endvertex of one of the four
paths.

\item Let $B_3$ be a graph of order $n  \ (n\geq 5)$ obtained from $K_4^-$ and four paths
$P_1,P_2,P_3,P_4$ of length $a,b,c,d \ (0\leq a,b,c,d\leq n-4, \ a+b+c+d=n-4)$ respectively by identifying
each vertex of $K_4^-$ with an endvertex of one of the four paths,
where $K_4^-$ denotes the graph obtained from a clique of order $4$
by deleting one edge.
\end{itemize}

\begin{thm}{\upshape \cite{WangMaoLiYe}}\label{th3-20}
Let $G$ be a connected graph of order $n \ (n\geq 5)$. Then $sdiam_4(G)=n-1$ if
and only if $G=T_{a,b,c,d}$ or $G=\triangle_{a,b,c,d}$ or
$G=\triangle^\prime_{a,b,c,d}$ or $G=B_1$ or $G=B_2$ or $G=B_3$.
\end{thm}

\subsubsection{Results for large $k$}

The following result is immediate.
\begin{obs}{\upshape \cite{MaoMelekianCheng}}\label{obs3-2}
Let $G$ be a graph of order $n$. Then $sdiam_n(G)=n-1$ if and only
if $G$ is connected.
\end{obs}

Mao, Melekian, and Cheng \cite{MaoMelekianCheng} characterized the graphs with
$sdiam_{n-1}(G)=d \ (n-2\leq d\leq n-1)$ and $sdiam_{n-2}(G)=d \ (n-3\leq d\leq n-1)$.
\begin{pro}{\upshape \cite{MaoMelekianCheng}}\label{pro3-4}
Let $G$ be a connected graph of order $n$. Then

$(1)$ $sdiam_{n-1}(G)=n-2$ if and only if $G$ is $2$-connected;

$(2)$ $sdiam_{n-1}(G)=n-1$ if and only if $G$ contains at least one
cut vertex.
\end{pro}

\begin{thm}{\upshape \cite{MaoMelekianCheng}}\label{th3-21}
Let $G$ be a connected graph of order $n \ (n\geq 5)$. Then

$(1)$ $sdiam_{n-2}(G)=n-3$ if and only if $\kappa(G)\geq 3$.

$(2)$ $sdiam_{n-2}(G)=n-2$ if and only if $\kappa(G)=2$ or $G$
contains only one cut vertex.

$(3)$ $sdiam_{n-2}(G)=n-1$ if and only if there are at least two cut
vertices in $G$.
\end{thm}

\subsection{Nordhaus-Gaddum-type results}

It is well known that if $diam(G)\geq 3$, then $diam(G)\leq 3$,
a result first proved by Harary and Robinson \cite{HararyRobinson}. A similar result,
namely that if $diam(G)\geq 4$, then $diam(G)\leq 2$, is due to Straffin
\cite{Straffin}. A common generalization of both results is given in the next theorem.
\begin{thm}{\upshape \cite{DankelmannSwartOellermann2, Mao}}\label{th3-22}
$(1)$ Let $G$ be a connected graph. If $sdiam_k(G)\geq k+r$, $1\leq r\leq k-1$, then
$sdiam_k(\overline{G})\leq 2k-r$.

$(2)$ Let $G$ be a connected graph. If $sdiam_k(G)\geq 2k$, then
$sdiam_k(\overline{G})\leq k$.
\end{thm}

Xu \cite{Xu} obtained the Nordhaus-Gaddum results for the diameter of graphs. In \cite{Mao}, Mao got the
Nordhaus-Gaddum results for the Steiner $k$-diameter of graphs.
\begin{thm}{\upshape \cite{Mao}}\label{th3-23}
Let $G\in \mathcal {G}(n)$ be a connected graph with a connected complement. Let $k$ be an integer with $3\leq
k\leq n$. Then

$(i)$ $2k-1-x\leq sdiam_k(G)+sdiam_k(\overline{G})\leq
\max\{n+k-1,4k-2\}$;

$(ii)$ $(k-1)(k-x)\leq sdiam_k(G)\cdot sdiam_k(\overline{G})\leq
\max\{k(n-1),(2k-1)^2\}$,

where $x=0$ if $n\geq 2k-2$ and $x=1$ if $n<2k-2$.
\end{thm}

For $k=n,n-1,n-2,3$, Mao \cite{Mao} improved the above Nordhaus-Gaddum results of
Steiner $k$-diameter and derived the following results.
\begin{obs}{\upshape \cite{Mao}}\label{obs3-3}
Let $G\in \mathcal {G}(n) \ (n\geq 3)$ be a connected graph with a connected complement. Then

$(i)$ $sdiam_n(G)+sdiam_n(\overline{G})=2n-2$;

$(ii)$ $sdiam_n(G)\cdot sdiam_n(\overline{G})=(n-1)^2$.
\end{obs}

\begin{pro}{\upshape \cite{Mao}}\label{pro3-5}
Let $G\in \mathcal {G}(n) \ (n\geq 5)$ be a connected graph with a connected complement. Then

$(i)$ $2n-4\leq sdiam_{n-1}(G)+sdiam_{n-1}(\overline{G})\leq 2n-2$;

$(ii)$ $(n-2)^2\leq sdiam_{n-1}(G)\cdot
sdiam_{n-1}(\overline{G})\leq (n-1)^2$.

Moreover,

$(a)$ $sdiam_{n-1}(G)+sdiam_{n-1}(\overline{G})=2n-4$ or
$sdiam_{n-1}(G)\cdot sdiam_{n-1}(\overline{G})=(n-2)^2$ if and only
if both $G$ and $\overline{G}$ are $2$-connected;

$(b)$ $sdiam_{n-1}(G)+sdiam_{n-1}(\overline{G})=2n-3$ or
$sdiam_{n-1}(G)\cdot sdiam_{n-1}(\overline{G})=(n-1)(n-2)$ if and
only if $\lambda(G)=1$ and $\overline{G}$ are $2$-connected, or
$\lambda(\overline{G})=1$ and $G$ are $2$-connected.

$(c)$ $sdiam_{n-1}(G)+sdiam_{n-1}(\overline{G})=2n-2$ or
$sdiam_{n-1}(G)\cdot sdiam_{n-1}(\overline{G})=(n-1)^2$ if and only
if $G$ satisfies one of the following conditions.

$\bullet$ $\kappa(G)=1$, $\Delta(G)=n-2$;

$\bullet$ $\kappa(G)=1$, $\Delta(G)\leq n-3$ and $G$ has a cut
vertex $v$ with pendant edge $e$ and pendant vertex $u$ such that
$G-u$ contains a spanning complete bipartite subgraph.
\end{pro}

\begin{pro}{\upshape \cite{Mao}}\label{pro3-6}
Let $G\in \mathcal {G}(n) \ (n\geq 5)$ be a connected graph with a connected complement. If both $G$ and
$\overline{G}$ contains at least two cut vertices, then

$(i)$ $2n-6\leq sdiam_{n-2}(G)+sdiam_{n-2}(\overline{G})\leq 2n-2$;

$(ii)$ $(n-3)^2\leq sdiam_{n-2}(G)\cdot
sdiam_{n-2}(\overline{G})\leq (n-1)^2$.

Otherwise,

$(iii)$ $2n-6\leq sdiam_{n-2}(G)+sdiam_{n-2}(\overline{G})\leq
2n-3$;

$(iv)$ $(n-3)^2\leq sdiam_{n-2}(G)\cdot
sdiam_{n-2}(\overline{G})\leq (n-1)(n-2)$.

Moreover, the upper and lower bounds are sharp.
\end{pro}

\begin{pro}{\upshape \cite{Mao}}\label{pro3-7}
Let $G\in \mathcal {G}(n) \ (n\geq 10)$ be a connected graph with a connected complement. Then

$(i)$ $6\leq sdiam_3(G)+sdiam_3(\overline{G})\leq n+2$;

$(ii)$ $9\leq sdiam_3(G)\cdot sdiam_3(\overline{G})\leq 3(n-1)$.

Moreover, the bounds are sharp.
\end{pro}

\subsection{Steiner diameter of graph products}

For Steiner diameter of joined, corona, cluster graphs, Wang, Mao, Cheng, and Melekian \cite{WangMaoChengMelekian} had the following.
\begin{pro}{\upshape \cite{WangMaoChengMelekian}}\label{pro3-8}
Let $k,n,m$ be two integers with $3\leq k\leq n(m+1)$, and let $G,H$
be two connected graphs with $n,m$ vertices, respectively.

$(i)$ If $3\leq k\leq n$, then $siam_k(G*H)=siam_k(G)+k$.

$(ii)$ If $n+1\leq k\leq mn$, then $siam_k(G*H)=n-1+k$.

$(iii)$ If $mn+1\leq k\leq (m+1)n$, then $siam_k(G*H)=n-1+mn$.
\end{pro}

\begin{pro}{\upshape \cite{WangMaoChengMelekian}}\label{pro3-9}
Let $G$ be a connected graph with $n$ vertices, and let $H$ be a
connected graph with $m \ (n\leq m)$ vertices. Let $k$ be an integer
with $3\leq k\leq n+m$.

$(i)$ If $k>m$, then $siam_k(G\vee H)=k-1$.

$(ii)$ If $n<k\leq m$ and $siam_k(H)=k-1$, then $siam_k(G\vee
H)=k-1$; if $n<k\leq m$ and $siam_k(H)\geq k$, then $siam_k(G\vee
H)=k$.

$(iii)$ If $3\leq k\leq n$, and $siam_k(G)\geq k$ or $siam_k(H)\geq
k$, then $siam_k(G\vee H)=k$; If $3\leq k\leq n$ and
$siam_k(G)=siam_k(H)=k-1$, then $siam_k(G\vee H)=k-1$.
\end{pro}

\begin{pro}{\upshape \cite{WangMaoChengMelekian}}\label{pro3-10}
Let $k,n,m$ be two integers with $3\leq k\leq nm$, and let $G,H$ be
two connected graphs with $n,m$ vertices, respectively.

$(i)$ If $m>n$ and $3\leq k\leq n$, then
$$
siam_k(G)+k\leq siam_k(G\odot H)\leq k \cdot
siam_{k+1}(H)+siam_k(G).
$$

$(ii)$ If $m>n$ and $n+1\leq k\leq m-1$, then
$$
n-1+k\leq siam_k(G\odot H)\leq n \cdot siam_{k+1}(H)+n-1.
$$

$(iii)$ If $m>n$ and $m\leq k\leq nm-n$, then
$$
n-1+k\leq siam_k(G\odot H)\leq mn-1.
$$

$(iv)$ If $m>n$ and $nm-n\leq k\leq nm$, then
$$
siam_k(G\odot H)=nm-1.
$$

$(v)$ If $m\leq n$ and $3\leq k<m$, then
$$
siam_k(G)+k\leq siam_k(G\odot H)\leq k \cdot
siam_{k+1}(H)+siam_k(G).
$$

$(vi)$ If $m\leq n$ and $m\leq k\leq n$, then
$$
siam_k(G)+k\leq siam_k(G\odot H)\leq k(m-1)+siam_k(G).
$$

$(vii)$ If $m\leq n$ and $n<k\leq mn-n$, then
$$
n-1+k\leq siam_k(G\odot H)\leq mn-1.
$$

$(viii)$ If $m\leq n$ and $mn-n<k\leq mn$, then
$$
siam_k(G\odot H)=mn-1.
$$
\end{pro}

For Steiner $k$-diameter of Cartesian product graphs, Mao, Cheng, and Wang
\cite{MaoChengWang} had the following.
\begin{thm}{\upshape \cite{MaoChengWang}}\label{th3-24}
Let $k,m,n$ be an integer with $3\leq k\leq mn$ and $n\leq m$. Let
$G,H$ be two connected graphs of order $n,m$, respectively.

$(1)$ If $k\leq n$, then
\begin{eqnarray*}
&&sdiam_k(G)+sdiam_k(H)\\
&\leq&sdiam_k(G\Box H)\\[0.1cm]
&\leq&sdiam_k(G)+sdiam_k(H)+(k-3)\min\{sdiam_k(G),sdiam_k(H)\}.
\end{eqnarray*}

$(2)$ If $n<k\leq m$, then
\begin{eqnarray*}
n-1+sdiam_k(H)
&\leq&sdiam_k(G\Box H)\\[0.1cm]
&\leq&n-1+sdiam_k(H)+(k-3)\min\{n-1,sdiam_k(H)\}.
\end{eqnarray*}

$(3)$ If $m<k\leq mn$, then
$$
n+m-2\leq sdiam_k(G\Box H)\leq m-1+(k-2)(n-1).
$$

$(4)$ If $mn-\kappa(G\Box H)+1\leq k\leq mn$, then $sdiam_k(G\Box
H)=k-1$.
\end{thm}

The following corollary is immediate from Theorem \ref{th3-24}.
\begin{cor}{\upshape \cite{MaoChengWang}}\label{cor3-5}
Let $G,H$ be two connected graphs of order at least $3$. Then
$$
sdiam_3(G\Box H)=sdiam_3(G)+sdiam_3(H).
$$
\end{cor}

To show the sharpness of the above upper and lower bound, we consider the following example.

\noindent{\bf Example 3.7.} \cite{MaoChengWang} $(1)$ For $k=3$, from Corollary
\ref{cor3-5}, we have $sdiam_k(G\Box H)=sdiam_k(G)+sdiam_k(H)$,
which implies that the upper and lower bounds in Theorem \ref{th3-24}
are sharp.

$(2)$ Let $G=P_n$ and $H=P_m$ with $5\leq n\leq m$. Then
$sdiam_4(G)=n-1$, $sdiam_4(H)=m-1$ and $sdiam_4(G\Box
H)=2(n-1)+(m-1)$, which implies that all the upper bounds in Theorem
\ref{th3-24} are sharp.

\vskip 0.3cm

By Theorem \ref{th2-8}, Mao, Cheng, and Wang \cite{MaoChengWang} derived the following results for Steiner diameter of lexicographic product graphs.
\begin{thm}{\upshape \cite{MaoChengWang}}\label{th3-25}
Let $k,n,m$ be three integers with $2\leq k\leq mn$. Let $G$ be a
connected graph of order $n$, and $H$ be a graph of order $m$.
Then
$$
sdiam_k(G\circ H)\leq \left\{
\begin{array}{ll}
sdiam_k(G)+k-2, &\mbox {\rm if} \ 2\leq k\leq n;\\[0.2cm]
\max\{n+k-3,k\}, &\mbox {\rm if} \ n<k\leq mn.
\end{array}
\right.
$$
and
$$
sdiam_k(G\circ H)\geq \left\{
\begin{array}{ll}
sdiam_k(G), &\mbox {\rm if} \ m+1\leq k\leq n;\\[0.2cm]
n-1, &\mbox {\rm if} \ \max\{n,m+1\}\leq k\leq mn;\\[0.2cm]
k-1, &\mbox {\rm if} \ 2\leq k\leq m.
\end{array}
\right.
$$
\end{thm}

To show the sharpness of the upper and lower bound in Theorem
\ref{th3-25}, they \cite{MaoChengWang} considered the following example.

\noindent{\bf Example 3.8.} \cite{MaoChengWang} Let $G=P_n$, and $H$ be a graph of order
$m$. If $k\leq \min\{2m,n\}$, then $sdiam_k(G\circ
H)=n+k-3=sdiam_k(G)+k-2$. If $\max\{n,m+1\}\leq k\leq 2m$, then
$sdiam_k(G\circ H)=n+k-3=\max\{n+k-3,k\}$. These implies that the
upper bounds in Theorem \ref{th3-25} are sharp. Let $G=K_n$ and
$H=K_m$. Then $G\circ H$ is a complete graph of order $mn$. If
$2\leq k\leq m$, then $sdiam_k(G\circ H)=k-1$. If $m+1\leq k\leq n$,
then $sdiam_k(G\circ H)=k-1=sdiam_k(G)$. These implies that the
lower bounds in Theorem \ref{th3-25} are sharp.

The following result is immediate from Proposition \ref{pro2-5}.
\begin{pro}{\upshape \cite{MaoChengWang}}\label{pro3-11}
Let $G,H$ be two connected graphs. Then
$$
sdiam_3(G\circ H)=\left\{
\begin{array}{ll}
diam(G)+1 &\mbox {\rm if} \ G=P_n, \ diam(G)\geq 2;\\[0.2cm]
sdiam_3(G) &\mbox {\rm if} \ G\neq P_n, \ diam(G)\geq 2;\\[0.2cm]
\min \{sdiam_3(H),3\} &\mbox {\rm if} \ G=K_n.
\end{array}
\right.
$$
\end{pro}

\section{Average Steiner Distance and Steiner Wiener Index}

Average Steiner distance is related to the Steiner
Wiener index via $SW_k(G)/\binom{n}{k}$.

The following results are due to Dankelmann, Oellermann,
Swart \cite{DankelmannOellermannSwart}, and Li, Mao, Gutman \cite{LiMaoGutman, LiMaoGutman2}.
\begin{pro}{\upshape \cite{DankelmannOellermannSwart, LiMaoGutman, LiMaoGutman2}}\label{pro4-1}
Let $k,n$ be two integers with $2\leq k\leq n$.

$(1)$ For a complete graph $K_n$, $SW_k(K_n)=
\binom{n}{k}(k-1)$.

$(2)$ For a path $P_n$, $\mu_k(P_n)=\frac{(k-1)(n+1)}{k+1}$; $SW_k(P_{n})=(k-1)\binom{n+1}{k+1}$.

$(3)$ For a star $S_n$, $SW_k(S_n)=\binom{n-1}{k-1}(n-1)$.

$(4)$ For a complete bipartite graph $K_{a,b} \ (2\leq k\leq a+b)$,
$$
SW_k(K_{a,b})=\left\{
\begin{array}{ll}
(k-1)\binom{a+b}{k}+\binom{a}{k}+\binom{b}{k}, &{\rm if}~1\leq k\leq
a;\\[4mm]
(k-1)\binom{a+b}{k}+\binom{b}{k}, &{\rm if}~a< k\leq
b;\\[4mm]
(k-1)\binom{a+b}{k}, &{\rm if}~b< k\leq a+b.
\end{array}
\right.
$$

$(5)$ Let $T$ be a graph obtained from a path $P_{t}$ and a star $S_{n-t+1}$ by identifying a pendant vertex of $P_{t}$ and the center $v$ of $S_{n-t+1}$, where $1\leq t\leq n-1$ and $k\leq n$. Then
$$
SW_{k}(T)=t\binom{n-1}{k}-\binom{t}{k+1}-\binom{n}{k+1}+\binom{n-t+1}{k+1}+(k-1)\binom{n}{k}.
$$

$(5)$
Let $G$ be a graph obtained from a clique $K_{n-r}$ and a star $S_{r+1}$ by identifying a vertex of $K_{n-r}$ and the center $v$ of $S_{r+1}$. For $k\leq r\leq n-1-k$,
$$
SW_{k}(G)=(n-1)\binom{n-1}{k-1}-\binom{n-r-1}{k}.
$$
\end{pro}

For $2\leq r<k$, Dankelmann, Oellermann, and Swart \cite{DankelmannOellermannSwart} established a relation between $\mu_r(G)$, $\mu_{k+1-r}(G)$, and $\mu_k(G)$.
\begin{thm}{\upshape \cite{DankelmannOellermannSwart}}\label{th4-1}
Let $G$ be connected weighted graph and $2\leq r\leq k-1$. Then
$$
\mu_k(G)\leq \mu_r(G)+\mu_{k+1-r}(G).
$$
\end{thm}

\begin{cor}{\upshape\cite{DankelmannOellermannSwart}}\label{cor4-1}
For $k\geq 3$, $\mu_k(G)\leq (k-1)\mu(G)$.
\end{cor}

The bounds in Theorem \ref{th4-1} and Corollary \ref{cor4-1} are sharp for the complete graph.
\begin{rem}{\upshape\cite{DankelmannOellermannSwart}}\label{rem4-1}
With essentially the same methods as those used in \cite{HenningOellemannSwart2} it can be shown that for
each connected graph $G$ of order $n$ and $3\leq k\leq n$,
$$
\mu_k(G)\leq \frac{k+1}{k-1}\mu_{k-1}(G).
$$
\end{rem}

It remains an open problem to determine a lower bound for $\mu_k(G)$ in terms of $\mu(G)$, but they conjectured that the smallest ratio $\mu_k(G)/\mu(G)$ taken over all connected graphs $G$ of order $n$ where $n\geq k$, is attained if $G$ is the path. More formally:
\begin{con}{\upshape\cite{DankelmannOellermannSwart}}\label{con4-1}
If $G$ is a connected graph of order $n$ and $3\leq k\leq n$, then
$$
\mu_k(G)\geq 3\frac{k-1}{k+1}\mu(G).
$$
\end{con}

In \cite{DankelmannOellermannSwart}, Dankelmann, Oellermann, and Swart  proved that the conjecture
is true for $k=3$ and $k=n$.

\subsection{Results for trees}

It is usual to consider $SW_k$ for $2 \leq k
\leq n-1$, but the definition implies $SW_1(G)=0$ and
$SW_n(G)=n-1$ for a connected graph $G$ of order $n$.
\begin{thm}{\upshape \cite{LiMaoGutman}}\label{th4-2}
Let $T$ be a tree of order $n$, possessing $p$ pendent vertices.
Then
$$
SW_{n-1}(T) = n(n-1)-p,
$$
irrespective of any other structural detail of $T$.
\end{thm}

\begin{thm}{\upshape \cite{LiMaoGutman2}}\label{th4-3}
Let $T$ be a tree of order $n$, possessing $p$ pendant vertices.
Let $q$ be the number of vertices of degree $2$ in $T$ that are adjacent to a pendant vertex. Then
$$
SW_{n-2}(T)=\frac{1}{2} \big( n^3-2n^2+n-2np+2p-2q \big)\,.
$$
\end{thm}

Let $G$ be any graph (not necessarily connected) with vertex set
$V(G)$. Let $e$ be an edge of $G$, connecting the vertices $x$ and
$y$. Define the sets
\begin{eqnarray*}
\mathcal N_1(e) & = & \{u\,|\, u \in V(G), d(u,x)<d(u,y) \} \\
\mathcal N_2(e) & = & \{u\,|\, u \in V(G), d(u,x)>d(u,y) \}
\end{eqnarray*}
and let their cardinalities be $n_1(e) = |\mathcal N_1(e)|$ and
$n_2(e) = |\mathcal N_2(e)|$, respectively. In other words, $n_1(e)$
counts the vertices of $G$, lying closer to one end of the edge $e$
than to its other end, and the meaning of $n_2(e)$ is analogous.

In his seminal paper \cite{Wiener}, Wiener discovered the following
result:
\begin{pro}{\upshape \cite{Wiener}}\label{pro4-2}
If $T$ is a tree, then for its Wiener index holds:
$$
W(T) = \sum_{e \in E(T)} n_1(e)\,n_2(e) .
$$
\end{pro}

Li, Mao, and Gutman \cite{LiMaoGutman} stated the generalization of Proposition \ref{pro4-2} to Steiner Wiener indices:
\begin{thm}{\upshape \cite{LiMaoGutman}}\label{th4-4}
Let $k$ be an integer such that $2 \leq k \leq n$. If $T$ is a tree,
then for its Steiner $k$-Wiener index holds:
\begin{equation}          \label{eq4-3}
SW_k(T) = \sum_{e \in E(T)} \sum_{i=1}^{k-1}
\binom{n_1(e)}{i}\,\binom{n_2(e)}{k-i} .
\end{equation}
\end{thm}

\begin{cor}{\upshape \cite{LiMaoGutman}}\label{cor4-2}
$(1)$ Proposition \ref{pro4-2} is obtained from Eq. (\ref{eq4-3}) by
setting $k=2$.

$(2)$ If $k=3$, then the Steiner $k$-Wiener index of a tree of order $n$
is directly related to the ordinary Wiener index as
$$
SW_3(T) = \frac{n-2}{2}\,W(T).
$$
\end{cor}

Dankelmann, Oellermann, and Swart \cite{DankelmannOellermannSwart} investigated
the average Steiner $k$-distance of trees by establishing sharp upper and lower bounds for this parameter.
\begin{thm}{\upshape \cite{DankelmannOellermannSwart}}\label{th4-5}
Let $T$ be a tree of order $n\geq k$ and $2\leq r\leq k-1$. Then
$$
\mu_k(T)\leq \frac{n}{r}\mu_r(T).
$$
Furthermore, equality holds if and only if $T$ is a star.
\end{thm}

\begin{cor}{\upshape\cite{DankelmannOellermannSwart}}\label{cor4-3}
$(1)$ If $T$ is a tree of order $n\geq k$, then
$\mu_k(T)\leq \frac{k}{2}\mu(T)$.

$(2)$ If $T$ is a tree of order $n\geq k$ and $2\leq r\leq k-2$, then
$\mu_k(T)\leq \mu_r(T)+\mu_{k-r}(T)$.
\end{cor}

The upper and lower bounds for average Steiner distance is also obtained in the same paper \cite{DankelmannOellermannSwart}.
\begin{pro}{\upshape\cite{DankelmannOellermannSwart}}\label{pro4-3}
If $T$ is a tree of order $n \ (2\leq k\leq n)$, then
$$
k\left(1-\frac{1}{n}\right)\leq \mu_k(T)\leq \frac{k-1}{k+1}(n+1)
$$
equality holds if and only if $T$ is a star or path, respectively, or in either case if $k=n$.
\end{pro}

\begin{rem}{\upshape\cite{DankelmannOellermannSwart}}\label{rem4-2}
For $k=2$, Proposition \ref{pro4-3} was already observed in \cite{EntringerJacksonSnyder} and \cite{Lovasz2}.
\end{rem}

\subsection{Algorithmic aspect for average Steiner distance}

An $O(kn^2)$ procedure is developed in \cite{BeinekeOellermannPippert} for calculating the $k$-distances of all vertices of a
tree $T$ of order $n$. Since the Steiner distance of any set $S$ of $k$ vertices contributes the same
amount to the $k$-distance of each vertex in $S$, it follows that
$$
\mu_k(T)=\sum_{v\in V(T)}\frac{d_k(v)}{n} \binom{n}{k}.
$$
Hence, the procedure developed in \cite{BeinekeOellermannPippert} gives a $O(kn^2$) procedure for finding $\mu_k(T)$ for a tree
$T$ which is considerably more efficient than the brute force method of calculating the Steiner
distance of all $\binom{n}{k}$ sets of vertices if $n>k$.

Dankelmann, Oellermann, and Swart \cite{DankelmannOellermannSwart} outlined an even more efficient algorithm that computes the average $k$-distance of
a tree without first computing the $k$-distance of each vertex. For a graph $G$ let $m(G)$ denote
the maximum order among all components of $G$. The algorithm we now describe is based
on the proof of Theorem \ref{th4-5}, where the average $k$-distance of a tree is expressed implicitly
in terms of the values $m(T-e)$ for $e\in E(T)$. It computes the average $k$-distance of a tree
of order $n$ using $O(n)$ graph operations and $O(kn)$ arithmetic operations. In \cite{DankelmannOellermannSwart}, the
\emph{edge weight} $\omega_k(e)$ for $e\in E(T)$ is defined by
$$
\omega_k(e)=\binom{n}{k}-\binom{m(T-e)}{k}-\binom{n-m(T-e)}{k}.
$$
Then $\omega_k(e)$ counts the number of $k$-sets $S\subseteq V(T)$ that have at least one vertex in each
component of $T-e$. Thus $\omega_k(e)$ equals the number of Steiner trees containing $e$ and we have
\begin{eqnarray*}
\mu_k(T)&=&\sum_{e\in E(T)}\omega_k(e)\binom{n}{k}^{-1}\\[0.1cm]
&=&n-1-\sum_{e\in E(T)}\left[\binom{m(T-e)}{k}+\binom{n-m(T-e)}{k}\right]\binom{n}{k}^{-1}.
\end{eqnarray*}
Therefore is suffices to compute the values $m(T-e)$ for all $e\in E(T)$, which is possible in
$O(n)$ time and to apply the above equlity which requires at most $O(kn)$ multiplications and divisions.

\subsection{Upper and lower bounds}

In \cite{DankelmannOellermannSwart}, the range for the average $k$-distance of a connected graph of given order was
determined, generalizing a result for $k=2$ obtained by Entringer, Jackson, and Snyder \cite{EntringerJacksonSnyder}, Doyle and
Graver \cite{DoyleGraver}, and Lov\'{a}sz \cite{Lovasz2}.
\begin{thm}{\upshape\cite{DankelmannOellermannSwart}}\label{th4-6}
Let $G$ be a connected graph of order $n$ and let $2\leq k\leq n-1$. Then
$$
k-1\leq \mu_k(G)\leq \frac{k-1}{k+1}(n+1).
$$
Equality holds on the left (or right) if and only if $G$ is $(n+1-k)$-connected (or if
$G$ is a path, respectively).
\end{thm}

The upper bound for the average $k$-distance given in Theorem \ref{th4-6} can be improved
for $2$-connected graphs. Plesn\'{\i}k \cite{Plesnik} showed that the cycle of order $n$, $C_n$, is the unique $2$-connected graph with given order $n$ and maximum average distance. This result was
generalized for the average $k$-distance in \cite{DankelmannOellermannSwart}.
It is remarkable that Plesn\'{\i}k's result can easily
be generalized for $k=2$ and $2\ell$-connected graphs (see \cite{FavaronKouiderMaheo}), which seems not to be the
case for $k\geq 3$.
\begin{thm}{\upshape\cite{DankelmannSwartOellermann}}\label{th4-7}
Let $G$ be a $2$-connected graph of order $n$ and let $2\leq k\leq n$. Then
$$
\mu_k(G)\leq \mu_k(C_n).
$$
Equulity holds if and only if $G=C_n$, or $k\geq n-1$.
\end{thm}

In \cite{Plesnik}, Plesn\'{\i}k proved that, apart from the obvious restriction $1\leq \mu(G)\leq diam(G)$,
the average distance is independent of the diameter and the radius.
\begin{thm}{\upshape\cite{Plesnik}}\label{th4-8}
Let $r,d$ be positive integers with $d\leq 2r$ and let $t\in \Re$ with $1\leq t\leq d$.
For every $\varepsilon>0$ there exists a graph $G$ with diameter $d$, radius $r$, and
$$
|\mu(G)-t|<\varepsilon.
$$
\end{thm}

It is natural to ask if there is a similar statement for the $k$-diameter and the average
$k$-distance. An answer in the affirmative is stated by Dankelmann, Swart, and Oellermann in \cite{DankelmannSwartOellermann}.
\begin{thm}{\upshape\cite{Bloom}}\label{th4-9}
Let $k,d$ he positive integers, $k\geq 2$, and let $t\in \Re$ with $k-1\leq t\leq d$. For
every $\varepsilon>0$ there exists a graph $G$ with $k$-diameter $d$ and
$$
|\mu_k(G)-t|<\varepsilon.
$$
\end{thm}

Dankelmann, Swart, and Oellermann \cite{DankelmannSwartOellermann} remarked that Theorem \ref{th4-9} is not a generalization of Plesn\'{\i}k's Theorem \ref{th4-8}, since
it does not allow us to prescribe also the Steiner $k$-radius. The problem of finding such a generalization
requires the determination of the possible values for the Steiner $k$-radius of a graph
of given Steiner $k$-diameter. This problem is still unsolved.

In \cite{ChartrandOellermannTianZou} it was conjectured that the
Steiner $k$-diameter of a graph $G$ never exceeds $\frac{k}{k-1}srad_k(G)$. This conjecture was disproved
by Henning, Oellermann, and Swart \cite{HenningOellermannSwart}, where the bound $sdiam_k(G)\leq [2(k+1)/(2k-1)]srad_k(G)$ was
conjectured.

The problem of determining a sharp lower bound for the average $k$-distance of
a connected graph with $n$ vertices and $m$ edges, where $k\geq 3$, is considerably more
difficult than the corresponding problem for $k=2$. The latter one was solved in \cite{EntringerJacksonSnyder}.
The following bound shows that the complete $r$-partite Tur\'{a}n graphs are optimal in
this regard. It remains an open problem to determine the graphs of given order and
size that minimize the average $k$-distance.
\begin{thm}{\upshape\cite{DankelmannSwartOellermann}}\label{th4-10}
Let $G$ be a graph of order $n$ and size $m$. Then
$$
\mu_k(G)\geq k-1+n\binom{n-\frac{2m}{n}-1}{k-1}\binom{n}{k}^{-1},
$$
where for a real number $a$ and a positive integer $b$ the binomial coefficient $\binom{a}{b}$ is
defined as $a(a-1)\ldots (a-b+1)/b!$.
\end{thm}

The bound given is sharp if $k$ is a multiple of $r$. It is attained by the complete
$r$-partite Tur\'{a}n graph.

Tomescu and Melter \cite{TomescuMelter} determined the range for the average distance of a graph
of given order and chromatic number and also the extremal graphs. Dankelmann, Swart, and Oellermann in \cite{DankelmannSwartOellermann} showed in the
following generalization, that the same graphs are also extremal for $k\geq 3$, though there
are other ones as well.

For $r<n$, let $H_{n,r}$ be the graph obtained from a complete graph $K_r$ and
a path of order $n-r$ with end vertices $v_1$ and $v_1'$ by joining $v_1'$ to one vertex of $K_r$.
For $r=n$, let $H_{n,r}$ be the complete graph $K_n$ and let $v_1$ be a vertex of $K_n$.
\begin{thm}{\upshape\cite{DankelmannSwartOellermann}}\label{th4-11}
Let $G$ be a connected graph of order $n \  (2\leq k\leq n)$ and chromatic number $r$
and let $v$ be a vertex of $G$. Then
\begin{itemize}
\item[] $(1)$ $d_k(v,G)\leq d_k(v_1,H_{n,r})$

\item[] $(2)$ $\mu_k(G)\leq \mu_k(H_{n,r})$,
\end{itemize}
with equality if and only if $v=v_1$ and $G=H_{n,r}$, respectively.
\end{thm}

Dankelmann, Swart, and Oellermann in \cite{DankelmannSwartOellermann} remarked that Theorem \ref{th4-10} yields a sharp lower bound for the $k$-distance of a connected
graph $G$ of given order $n$ and chromatic number $r$. From
$$
e(G)\leq n^2\frac{r-1}{2r}
$$
and Theorem \ref{th4-10} we have immediately
$$
\mu_k(G)\geq k-1+n\binom{n/r-1}{k-1}\binom{n}{k}^{-1}.
$$
This bound is sharp if $n$ is a multiple of $r$. Examples for equality in the above equation
are the $r$-partite Tur\'{a}n graph $T_{n,r}$ and, for $k>n/r$, the graph $T_{n,r}-e$.

\subsection{Inverse problem}

The seemingly elementary question: ``{\it which natural numbers are Wiener indices of graphs ?}" was much investigated in the past; see  \cite{FinkLuzarSkrekovski, GutmanYeh, GutmanYehChen, Wagner, WagnerWangYu}.

Li, Mao, and Gutman \cite{LiMaoGutman2} considered the analogous question for Steiner Wiener indices:
\begin{prop}{\upshape\cite{LiMaoGutman2}}\label{prop4-1}
Fixed a positive integer
$k$, for what kind of positive integer $w$ does there exist a connected graph $G$ (or a tree $T$) of order $n\geq k$ such that $SW_k(G)=w$ (or $SW_k(T)=w$) ?
\end{prop}

For $k=2$, the authors have nice results in \cite{GutmanYeh, Wagner2}, completely solved a conjecture
by Lepovi\'{c} and Gutman \cite{LepovicGutman} for trees, which states that for all but 49 positive integers $w$ one can find a tree with Wiener index $w$. This is different from  Problem \ref{prop4-1} for trees, since Li, Mao, and Gutman \cite{LiMaoGutman2} considered graphs or trees with order $n$.

If $G$ is a connected graph or a tree of order $n$, then for $k=n$, $SW_k(G)=n-1$. Thus the
following result is immediate.
\begin{pro}{\upshape \cite{LiMaoGutman2}}\label{pro4-4}
For a positive integer $w$, there exists a connected graph $G$ or a tree $T$ of order $n$
such that $SW_n(G)=w$ or $SW_n(T)=w$ if and only if $w=n-1$.
\end{pro}

For $k=n-1$, Li, Mao, and Gutman \cite{LiMaoGutman2} had the following results.
\begin{pro}{\upshape \cite{LiMaoGutman2}}\label{pro4-5}
For a positive integer $w$, there exists a connected graph $G$ of order $n$ such that
$SW_{n-1}(G)=w$, if and only if $n^2-2n\leq w\leq n^2-n-2$.
\end{pro}

\begin{pro}{\upshape \cite{LiMaoGutman2}}\label{pro4-6}
For a positive integer $w$, there exists a tree $T$ of order $n$ such that
$SW_{n-1}(T)=w$ if and only if $n^2-2n+1\leq w\leq n^2-n-2$.
\end{pro}

For $k=n-2$, Li, Mao, and Gutman \cite{LiMaoGutman2} derived the following result for trees.
\begin{thm}{\upshape \cite{LiMaoGutman2}}\label{th4-12}
For a positive integer $w$, there exists a tree $T$ of order $n \ (n\geq 5)$, possessing $p$ pendant vertices, such that
$SW_{n-2}(T)=w$ if and only if $w=\frac{1}{2} \big( n^3-2n^2+n-2np+2p-2q \big)$, where $q$ is the number of vertices of degree $2$ in $T$ that are adjacent to a pendant vertex, and one of the following holds:
\begin{itemize}
\item[] $(1)$ $2\leq q\leq \lfloor\frac{n-1}{2}\rfloor$ and $q\leq p\leq n-q-1$;

\item[] $(2)$ $q=1$ and $3\leq p\leq n-2$;

\item[] $(3)$ $q=0$ and $4\leq p\leq n-1$.
\end{itemize}
\end{thm}

\begin{pro}{\upshape \cite{LiMaoGutman2}}\label{pro4-7}
For a positive integer $w$, there exists a tree $T$ of order $n$ such that
$SW_{k}(T)=w$ if
$$
w=t\binom{n-1}{k}-\binom{t}{k+1}-\binom{n}{k+1}+\binom{n-t+1}{k+1}+(k-1)\binom{n}{k},
$$
where $1\leq t\leq n-1$ and $k\leq n$.
\end{pro}

\begin{pro}{\upshape \cite{LiMaoGutman2}}\label{pro4-8}
For a positive integer $w$, there exists a connected graph $G$ of order $n$ such that
$SW_{k}(G)=w$ if $w$ satisfies one of the following conditions:
\begin{itemize}
\item[] $(1)$ $w=t\binom{n-1}{k}-\binom{t}{k+1}-\binom{n}{k+1}+\binom{n-t+1}{k+1}+(k-1)\binom{n}{k}$, where $1\leq t\leq n-1$ and $k\leq n$.

\item[] $(2)$ $w=(n-1)\binom{n-1}{k-1}-\binom{n-r-1}{k}$, where $k\leq r\leq n-1-k$ and $k\leq n$.
\end{itemize}
\end{pro}

\subsection{Graph products}

Yeh and Gutman \cite{YehGutman} investigated the Wiener
index of graph products and obtained the following results.
\begin{thm}{\upshape \cite{YehGutman}}\label{th4-13}
Let $G$ be a connected graph with $n$ vertices, and let $H$ be a
connected graph with $m$ vertices. Then

$(1)$ $W(G\vee H)=e(G)+e(H)+mn+2\left[\binom{n}{2}-e(G) + \binom{m}{2}-e(H)\right]$.

$(2)$ $W(G\circ H)=m^2\,\big[ (W(G)+n)-n(e(H)+m) \big]$.

$(3)$ $W(G\Box H)=m^2\,W(G)+n^2\,W(H)$.

$(4)$ $W(G\odot H) = m^2\,W(G)+n\,W(H)+m (n^2-n)\,d(v|H)$,
where $v$ is the root-vertex of $H$ and
$$
d(v|H)=\sum_{u\in V(H)} d(u,v)\,.
$$

$(5)$ $W(G\ominus H)=(m+1)^2\,W(G)+n \big[m^2-e(H) \big] + mn(m+1)(n-1)$.
\end{thm}

In \cite{MaoWangGutman}, Mao, Wang, Gutman studed the $k$-th Steiner Wiener index of the above specified graph products.
\begin{thm}{\upshape \cite{MaoWangGutman}}\label{th4-14}
Let $G$ be a connected graph with $n$ vertices, and let $H$ be a
connected graph with $m \ (n\geq m)$ vertices. Let $k$ be an integer,
$3\leq k\leq n+m$.

$(1)$ If $k>n$, then
$$
SW_k(G\vee H)=(k-1)\binom{n+m}{k}\,.
$$

$(2)$ If $k\leq m$, then
$$
SW_k(G\vee H)=(k-1)\binom{n+m}{k}
+ \binom{n}{k} + \binom{m}{k}-x-y,
$$
where $x$ is the number of the $k$-subsets of $V(G)$ such that the
subgraph induced by each $k$-subset is connected, and $y$ is the
number of the $k$-subsets of $V(H)$ such that the subgraph induced by
each $k$-subset is connected.

$(3)$ If $m< k\leq n$, then
$$
SW_k(G\vee H)=(k-1)\binom{n+m}{k}
+\binom{n}{k}+(k-1)\binom{m}{k}-x\,.
$$
\end{thm}

\begin{thm}{\upshape \cite{MaoWangGutman}}\label{th4-15}
Let $G$ be a connected graph with $n$ vertices, and let $H$ be a
connected graph with $m$ vertices. Let $k$ be an integer, $2\leq
k\leq nm$. Then
\begin{eqnarray*}
SW_k(G\circ H) & = & nk \binom{m}{k}-nx +
\sum_{\ell=2}^{k} \binom{m}{r_1}\binom{m}{r_2} \cdots
\binom{m}{r_{\ell}}SW_{\ell}(G) \nonumber \\[3mm]
& & + \sum_{\ell=2}^k (k-\ell)\binom{n}{\ell}\binom{m\ell-\ell}{k-\ell} \label{Eq-th2.3}
\end{eqnarray*}
where $\sum_{i=1}^{\ell}r_{i}=k$, $r_i\geq 1$ and $x$ is the number
of the $k$-subsets of $V(H)$ such that the subgraph induced by each
$k$-subset is connected in $H$.
\end{thm}

\begin{thm}{\upshape \cite{MaoWangGutman}}\label{th4-16}
Let $G$ be a connected graph with $n$ vertices, and let $H$ be a
connected graph with $m$ vertices. Let $k$ be an integer with $2\leq
k\leq nm$. Then
\begin{eqnarray*}
&&\sum_{x=2}^{k}\binom{m}{r_1}\binom{m}{r_2}\cdots \binom{m}{r_{x}}SW_{x}(G)
+ \sum_{y=2}^{k}\binom{n}{s_1}\binom{n}{s_2}\cdots
\binom{n}{s_{y}}SW_{y}(G) \\[3mm]
&\leq &SW_k(G\Box H)
\leq \frac{k}{2}\left[\sum_{x=2}^{k}\binom{m}{r_1}\binom{m}{r_2}
\cdots \binom{m}{r_{x}}SW_{x}(G) +
\sum_{x=2}^{k}\binom{n}{s_1}\binom{n}{s_2}\cdots \binom{n}{s_{y}}SW_{y}(G)\right]
\end{eqnarray*}
where $\sum_{i=1}^{x}r_{i}=k$ and $r_i\geq 1$, and
$\sum_{i=1}^{y}s_{i}=k$ and $s_i\geq 1$.
\end{thm}

\begin{rem}{\upshape \cite{MaoWangGutman}}\label{rem4-3}
Suppose that $k=2$. Then $x=y=2$,
$r_1=r_2=\ldots=r_{x}=1$, $\sum_{i=1}^{x}r_{i}=2$,
$s_1=s_2=\ldots=s_{y}=1$, $\sum_{i=1}^{y}s_{i}=2$. Therefore,
$$
SW_2(G\Box H) = m^2\,W(G)+n^2\,SW(H)\,.
$$
Thus, the upper and lower bounds in Theorem \ref{th4-16} are sharp.
\end{rem}

Let $v$ is the root vertex of $H$ and
$$
d(v,k|H) = \sum_{\overset{v\in V(H)\,,\,S\subseteq V(H)}{|S|=k}} d(S)\,.
$$

\begin{thm}{\upshape \cite{MaoWangGutman}}\label{th4-17}
Let $G$ be a connected graph with $n$ vertices, and let $H$ be a
connected graph with $m$ vertices. Let $k$ be an integer, $2\leq
k\leq nm$. Then
\begin{eqnarray*}
SW_k(G\odot H) &=&n\,SW_{k}(H)+\sum_{\ell=2}^{k} \binom{m}{r_1}\binom{m}{r_2}
\cdots \binom{m}{r_{\ell}}SW_{\ell}(G) \\[3mm]
& + & \sum_{\ell=2}^{k}\binom{n}{\ell}\left[\sum_{j=1}^{\ell}
\prod_{\overset{x=1}{x\neq j}}^{\ell} \binom{m}{r_x}d(v,k|H)\right]
\end{eqnarray*}
where $\sum_{x=1}^{\ell}r_{x}=k$, $r_x\geq 1$ and $v$ is the
root-vertex of $H$.
\end{thm}

\begin{thm}{\upshape \cite{MaoWangGutman}}\label{th4-18}
Let $G$ be a connected graph with $n$ vertices, and let $H$ be a
connected graph with $m$ vertices. Let $k$ be an integer, $2\leq k\leq nm$.
Then
\begin{eqnarray*}
SW_k(G\ominus H)
&=&\sum_{\ell=2}^{k}\binom{m\!+\!1}{r_1}\binom{m\!+\!1}{r_2}\cdots
\binom{m\!+\!1}{r_{\ell}} SW_{\ell}(G)
+\binom{m}{k\!-\!1}(k\!-\!1)n+kn\binom{m}{k} \nonumber \\[3mm]
&&-xn+\sum_{\ell=2}^{k} \binom{n}{\ell}\left[\sum_{j=1}^{\ell}
\prod_{\overset{x=1}{x\neq j}}^{\ell}
\binom{m+1}{r_x}\left[\binom{m}{r_j-1}(r_j-1)+r_j\binom{m}{r_j}
- x_j\right]\right] \label{eqth2.12}
\end{eqnarray*}
where $\sum_{x=1}^{\ell}r_{x}=k$, $r_x\geq 1$, $x$ is the number of
the $k$-subsets of $V(H)$ such that the subgraph induced by each
$k$-subset is connected in $H$, and $x_j$ is the number of the
$r_j$-subsets of $V(H)$ such that the subgraph induced by each
$r_j$-subset is connected in $H$.
\end{thm}

\begin{rem}
One can see that Theorems Theorem \ref{th4-14}, \ref{th4-15}, \ref{th4-16},
\ref{th4-17}, and \ref{th4-18} are extensions $(1),(2),(3),(4),(5)$ of Theorem \ref{th4-13}, respectively. In all considered case for $k=2$ the new results can be reduced to already
known ones.
\end{rem}

\subsection{Nordhaus-Gaddum-type results}

In \cite{ZhangWu}, Zhang and Wu studied the Nordhaus--Gaddum problem for
the Wiener index and proved that for $G \in \mathcal G(n)$,
$$
3\binom{n}{2}\leq W(G)+W(\overline{G})\leq
\frac{1}{6}(n^3+3n^2+2n-6)\,.
$$

Mao, Wang, Gutman and Li \cite{MaoWangGutmanLi} investigated the analogous problem for the Steiner
Wiener index.
\begin{thm}{\upshape \cite{MaoWangGutmanLi}}\label{th4-19}
Let $G\in \mathcal {G}(n)$ be a connected graph with a connected complement $\overline{G}$. Let $k$ be an integer such that $3\leq
k\leq n$. Then: \\[2mm]
{\rm (1)} \hspace{10mm} $\binom{n}{k}(2k-1-x)\leq SW_k(G)+SW_k(\overline{G})\leq
\max\{n+k-1,4k-2\}\binom{n}{k}$ \\[2mm]
where if $n\geq 2k-2$, then $x=0$; $x=1$ for positive integer $n$. \\[3mm]
{\rm (2)} $(k-1)^2\,\binom{n}{k}^2 \leq SW_k(G)\cdot SW_k(\overline{G})
\leq \max\{k(n-1),(2k-1)^2\}\binom{n}{k}^2$\,.

Moreover, the lower bounds are sharp.
\end{thm}

For $k=n$, the following result is immediate.
\begin{obs}{\upshape \cite{MaoWangGutmanLi}}\label{obs4-1}
Let $G\in \mathcal {G}(n)$ be a connected graph with a connected complement $\overline{G}$. Then \\
{\rm (1)} \hspace{10mm} $SW_n(G)+ SW_n(\overline{G})=2n-2$\,; \\[2mm]
{\rm (2)} \hspace{10mm} $SW_n(G)\cdot SW_n(\overline{G})=(n-1)^2$\,.
\end{obs}

For $k=n-1$, they derived \cite{MaoWangGutmanLi} the following result.
\begin{pro}{\upshape \cite{MaoWangGutmanLi}}\label{prop4-9}
Let $G\in \mathcal {G}(n) \ (n\geq 5)$ be a connected graph with a connected complement $\overline{G}$.

$(1)$ If $G$ and $\overline{G}$ are both $2$-connected, then
$SW_{n-1}(G)+SW_{n-1}(\overline{G})=2n(n-2)$ and $SW_{n-1}(G)\cdot
SW_{n-1}(\overline{G})=n^2\,(n-2)^2$.

$(2)$ If $\kappa(G)=1$ and $\overline{G}$ is $2$-connected, then
$SW_{n-1}(G)+SW_{n-1}(\overline{G})=2n(n-2)+p$ and $SW_{n-1}(G)\cdot
SW_{n-1}(\overline{G})=n(n-2)(n^2-2n+p)$, where $p$ is the number of
cut vertices in $G$.

$(3)$ If $\kappa(G)=\kappa(\overline{G})=1$, $\Delta(G)\leq n-3$, and
$G$ has a cut vertex $v$ with pendent edge $uv$ such that $G-u$
contains a spanning complete bipartite subgraph, and
$\Delta(\overline{G})\leq n-3$ and $\overline{G}$ has a cut vertex
$q$ with pendent edge $pq$ such that $G-p$ contains a spanning
complete bipartite subgraph, then
$SW_{n-1}(G)+SW_{n-1}(\overline{G})=2(n-1)^2$ and $SW_{n-1}(G)\cdot
SW_{n-1}(\overline{G})=(n-1)^4$.

$(4)$ If $\kappa(G)=\kappa(\overline{G})=1$,
$\Delta(\overline{G})=n-2$, $\Delta(G)\leq n-3$ and $G$ has a cut
vertex $v$ with pendent edge $uv$ such that $G-u$ contains a
spanning complete bipartite subgraph, then
$$
SW_{n-1}(G)+SW_{n-1}(\overline{G})=2(n-1)^2 \hspace{3mm} \mbox{or} \hspace{3mm}
SW_{n-1}(G)+SW_{n-1}(\overline{G})=2(n-1)^2+1
$$
and
$$
SW_{n-1}(G)\cdot SW_{n-1}(\overline{G})=(n-1)^4 \hspace{3mm} \mbox{or} \hspace{3mm}
SW_{n-1}(G)\cdot SW_{n-1}(\overline{G})=(n-1)^2\,(n^2-2n+2)\,.
$$

$(5)$ If $\kappa(G)=\kappa(\overline{G})=1$,
$\Delta(G)=\Delta(\overline{G})=n-2$, then
$$
2(n-1)^2\leq SW_{n-1}(G)+SW_{n-1}(\overline{G})\leq 2(n-1)^2+2
$$
and
$$
(n-1)^4 \leq SW_{n-1}(G)\cdot SW_{n-1}(\overline{G})\leq (n^2-2n+2)^2\,.
$$
\end{pro}

For $k=3$
and $n\geq 10$, from Theorem \ref{th4-19}, one can see that
$$
5 \binom{n}{3} \leq SW_3(G)+SW_3(\overline{G}) \leq (n+2) \binom{n}{3}
$$
and
$$
4 \binom{n}{3}^2 \leq SW_3(G)\cdot SW_3(\overline{G}) \leq 3(n-1) \binom{n}{3}^2\,.
$$

Mao, Wang, Gutman and Li \cite{MaoWangGutmanLi} improved these bounds and proved the following result.
\begin{thm}{\upshape \cite{MaoWangGutmanLi}}\label{th4-20}
Let $G\in \mathcal {G}(n) \ (n\geq 4)$ be a connected graph with a connected complement $\overline{G}$. Then \\
{\rm (1)}
\begin{eqnarray*}
&& 5 \binom{n}{3} \leq SW_3(G)+SW_3(\overline{G}) \\[4mm]
&\leq & \left\{
\begin{array}{ll}
7 \binom{n}{3}-3n+8, & {\rm if} \ n=6,7, \ {\rm and} \ sdiam_3(G)=5, \\
&\mbox{\rm ~~or} \ n=6,7, \ {\rm and} \ sdiam_3(\overline{G})=5;\\[3mm]
2 \binom{n+1}{4}+2 \binom{n-3}{3}+\frac{1}{2}(7n^2-35n+48), & {\rm otherwise}.
\end{array}
\right.
\end{eqnarray*}

\noindent
{\rm (2)}
\begin{eqnarray*}
&&6 \binom{n}{3}^2 + (n-2) \binom{n}{3}-(n-2)^2 \\[5mm]
& \leq & SW_3(G)\cdot SW_3(\overline{G}) \\[5mm]
&\leq&\left\{
\begin{array}{ll}
\frac{1}{4}\left[7 \binom{n}{3}-3n+8\right]^2, &{\rm if} \ n=6,7, \ {\rm and} \ sdiam_3(G)=5; \\
& {\rm ~~or} \ n=6,7, \ {\rm and} \ \ sdiam_3(\overline{G})=5;\\[3mm]
\left[\binom{n+1}{4}+ \binom{n-3}{3}+\frac{1}{4}\,(7n^2-35n+48) \right]^2,
&\mbox {\rm otherwise}.
\end{array}
\right.
\end{eqnarray*}

Moreover, the bounds are sharp.
\end{thm}

\subsection{Steiner Wiener index and Steiner betweenness centrality}

The \emph{betweenness centrality} $B(v)$ of a vertex $v\in V(G)$ is defined as the
sum of the fraction of all pairs of shortest paths that pass through $v$ across
all pairs of vertices in a graph:
$$
B(v)=\sum_{x,y\in V(G)\setminus \{v\}, \ x\neq y}\frac{\sigma_{x,y}(v)}{\sigma_{x,y}},
$$
where $\sigma_{x,y}$ denotes the number of all shortest paths between vertices $x$
and $y$ in a graph $G$ and $\sigma_{x,y}(v)$ denotes the number of all shortest paths
between vertices $x$ and $y$ in graph $G$ passing through the vertex $v$.

In a case when a graph models a social or communication network,  as the
name suggests, it measures the centrality of a vertex in a graph, by the
influence of a vertex in the dissemination of information over a network. It
has been independently introduced by Anthonisse in \cite{Anthonisse} and by Freeman in
\cite{Freeman}, and among other applications has been applied to detect communities
in networks \cite{GirvanNewman, NewmanGirvan}.

For a graph $G$, let $n(G)$ denote the number of its vertices. For a forest
(acyclic  graph) $F$  with $p \ (p>1)$ connected components $T_1,T_2,\ldots,T_p$
denote by $N_k(F)$ the sum over all partitions of $k$ into at least two nonzero
parts of products of combinations distributed among the $p$ components of
$F$:
$$
N_k(F)=\sum_{\ell_1+\ell_2+\ldots+\ell_p=k, \ 0\leq \ell_1,\ell_2,\ldots,\ell_p<k}\binom{n(T_1)}{\ell_1}\binom{n(T_2)}{\ell_2}\cdots \binom{n(T_p)}{\ell_p}.
$$
For a tree $T$, we define $N_k(T)=0$. Note that by the definition $\binom{n}{0}=1$ and
$\binom{n}{k}=0$  whenever $n<k$.

Kov\v{s}e \cite{Kovse} derived the following result for Steiner Wiener index.
\begin{thm}{\upshape \cite{Kovse}}\label{th4-21}
$(1)$ Let $T$  be a tree on $n$ vertices. Then
$$
SW_k(T)=\sum_{e\in E(T)}N_k(T-e).
$$

$(2)$ Let $T$  be a tree on $n$ vertices. Then
$$
SW_k(T)=\sum_{v\in V(T)}N_k(T-v)+(k-1)\binom{n}{k}.
$$
\end{thm}

The \emph{Steiner $k$-betweenness centrality} $B_k(v)$ of a vertex $v\in V(G)$ is defined
as the sum of the fraction of all $k$-Steiner trees that include $v$ as its non-terminal vertex across all combinations of $k$ vertices of $G$:
$$
B_k(v)=\sum_{S\in V(G)\setminus \{v\},\ |S|=k}\frac{\sigma_{S}(v)}{\sigma_{S}},
$$
where $\sigma_S$ denotes the number of all Steiner trees between vertices of $S$
in a graph $G$ and $\sigma_S(v)$ denotes the number of all Steiner trees between
vertices of $S$ in a graph $G$ that include also the vertex $v$ as a non-terminal
vertex.
\begin{thm}{\upshape \cite{Kovse}}\label{th4-22}
$(1)$ Let $G$ be a connected graph on $n$ vertices. Then
$$
SW_k(G)=\sum_{v\in V(G)}B_k(v)+(k-1)\binom{n}{k}.
$$
\end{thm}

For a graph $G$ on $n$ vertices the \emph{average $k$-Steiner betweenness} $\overline{B_k}(G)$ is
defined as
$$
\overline{B_k}(G)=\frac{1}{n}\sum_{v\in V(G)}B_k(v).
$$
\begin{cor}{\upshape \cite{Kovse}}\label{cor4-4}
Let $G$ be a connected graph on $n$ vertices. Then
$$
\overline{B_k}(G)=\frac{1}{n}\binom{n}{k}(\mu_k(G)-k+1).
$$
\end{cor}

\section{Steiner Center and Steiner Median}

In a graph $G$, a vertex $x$ is a cut-vertex if deleting $x$ and all edges incident to it increases the number of connected components. A \emph{block} of a graph is a maximal connected vertex-induced subgraph that has no cut vertices. A \emph{block graph} is a connected graph whose blocks are complete graphs. Note that trees are block graphs.

\subsection{Results for Steiner center}

Hedetniemi \cite{BuckleyMillerSlater} verified that every graph $H$
is the center of some graph $G$. As an extension of this result, Oellermann and Tian
\cite{OellermannTian} derived the
following result for Steiner center.
\begin{thm}{\upshape \cite{OellermannTian}}\label{th5-1}
Let $k\geq 2$ be an integer and $H$ be a graph. Then $H$ is the Steiner $k$-center
of some graph $G$.
\end{thm}

The construction given in \cite{OellermannTian} was described in \cite{Oellermann}. Let $H$ be
a given graph and $k\geq 2$. Let $G$ be obtained from $H$ by first
adding $2k$ new vertices $\{v_1,v_2,\ldots,v_k\}\cup \{u_1,u_2,\ldots,
u_k\}$ and then joining $v_i$ to every vertex of $H$ for $1\leq i\leq k$
and next adding the edges $u_iv_i$ for $1\leq i\leq k$.
Then, $e_k(u_i)=2k$, $e_k(v_i)=2k-1$ for $1\leq i\leq k$, and $e_k(u)=2k-2$ for all $u\in V(H)$. Hence, $H$ is the Steiner $k$-center
of $G$.

Even though every graph is the Steiner $k$-center of some graph,
the problem of finding the Steiner $k$-center of any given graph
appears to be quite difficult.

However, for trees, an efficient solution to this problem
was developed in \cite{OellermannTian} and extends the work done by Jordan
\cite{Jordan} in his 1869 paper on centers and centroids of trees.

The key result that leads to a recursive procedure for
finding the Steiner $k$-center of a tree states: \emph{For any tree of order $n\geq 3$
and integer $k$ with $3\leq k\leq n$, if $T$ has at least $k$
end-vertices and $T'$ is obtained by deleting the end-vertices
from $T$, then $C_{k}(T)\subseteq C_{k}(T')$. Moreover, if $T$ has at most $k-1$ end-vertices, then $C_{k}(T)=T$.}

Thus, Steiner $k$-centers of trees of order $3\leq k\leq n$ was
characterized as follows.
\begin{thm}{\upshape \cite{OellermannTian}}\label{th5-2}
A tree $H$ is the $k$-center of some tree if and
only if

$(1)$  $k\geq 3$ and $H$ has at most $k-1$ end-vertices, or

$(2)$  $k=2$ and $H$ is isomorphic to $K_1$ or $K_2$.
\end{thm}

Based on these results, the following procedure for finding the $k$-center of a tree was developed in \cite{OellermannTian}.
\begin{figure}[!hbpt]
\begin{center}
\includegraphics[scale=0.8]{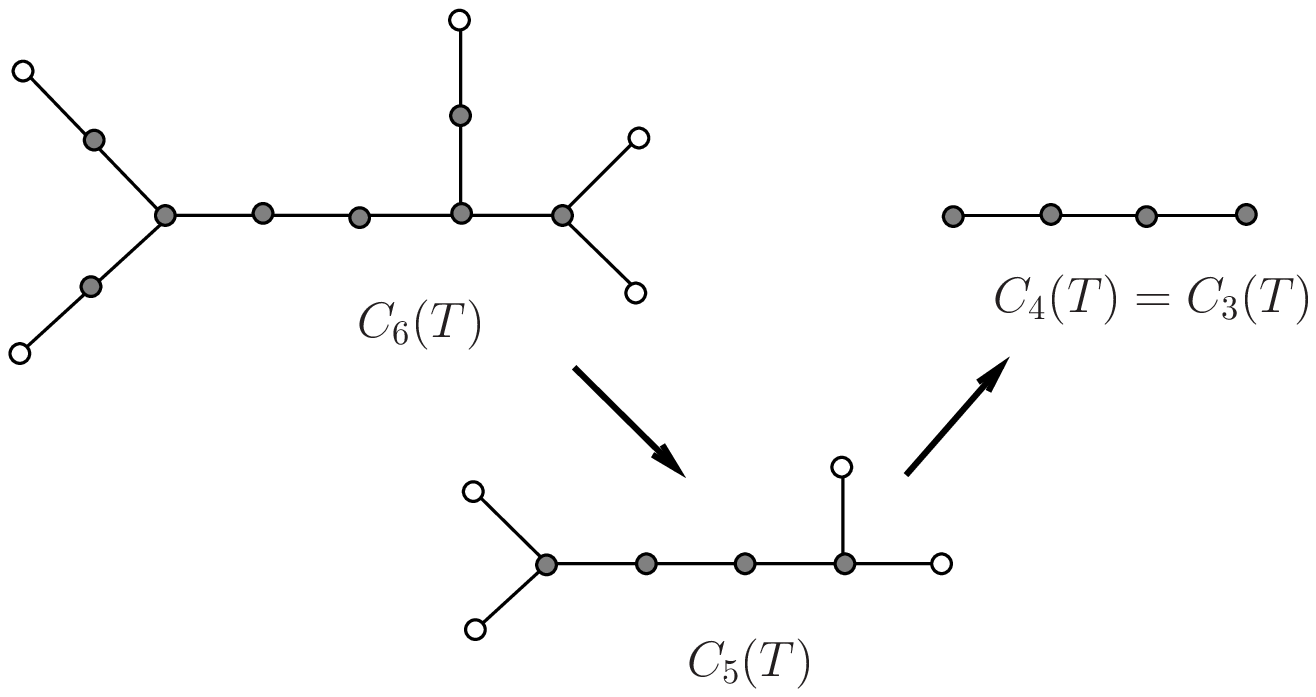}\\[0.5cm]
Figure 5.1: Graphs for Algorithm 5.1.
\end{center}\label{fig7}
\end{figure}

\noindent {\bf Algorithm 5.1.} Finding the Steiner $k$-center of a tree $T$ of
order $n\geq k\geq 2$.

$(i)$  $H\leftarrow T$

$(ii)$  If $H$ has at most $k-1$ end-vertices, or if $H\cong K_1$ or $K_2$
and $k=2$, output $H$ since $H$ is the Steiner $k$-center of $T$ and stop;
otherwise, continue.

$(iii)$  Delete the end-vertices from $H$ and let $H$ be the resulting
tree. Return to $(ii)$.

Figure 5.1 illustrates Algorithm 5.1 with $n=3,4,5$, and $6$.

This procedure gives the following result.
\begin{cor}{\upshape \cite{OellermannTian}}\label{cor5-1}
Let $k \ (k\geq 3)$ be an integer and $T$ a tree of order $n \ (k\leq n)$. Then
$$
C_{k-1}(T)\subseteq C_{k}(T).
$$
\end{cor}

Oellermann \cite{Oellermann} asked whether these containment
relationships hold for general graphs. However, Yeh, Chiang, and Peng \cite{YehChiangPeng} had a tree-like graph $F_1$ which has $C_3(F_1)\nsubseteq C_4(F_1)$, where $F_1$ is a graph obtained from $K_4^-$ by adding a pendant edge. Notice that $H$ is a partial $2$-tree, an
interval graph, and is also a distance-hereditary graph.

The following proposition even shows that graphs $G$ with $C_{n-2}(G)=V(G)$ and $C_{n-1}(G)=\{z\}$ can be constructed
systematically, where $n=|V(G)|$.
\begin{pro}{\upshape \cite{YehChiangPeng}}\label{pro5-1}
Let $G$ be a connected graph having exactly one cut-vertex $z$. Let $n=|V(G)|$. If $G$ has a $2$-vertex cut $\{a,b\}$ such that $z\notin \{a,b\}$, then $C_{n-2}(G)=V(G)$ and $C_{n-1}(G)=\{z\}$.
\end{pro}

\begin{figure}[!hbpt]
\begin{center}
\includegraphics[scale=0.8]{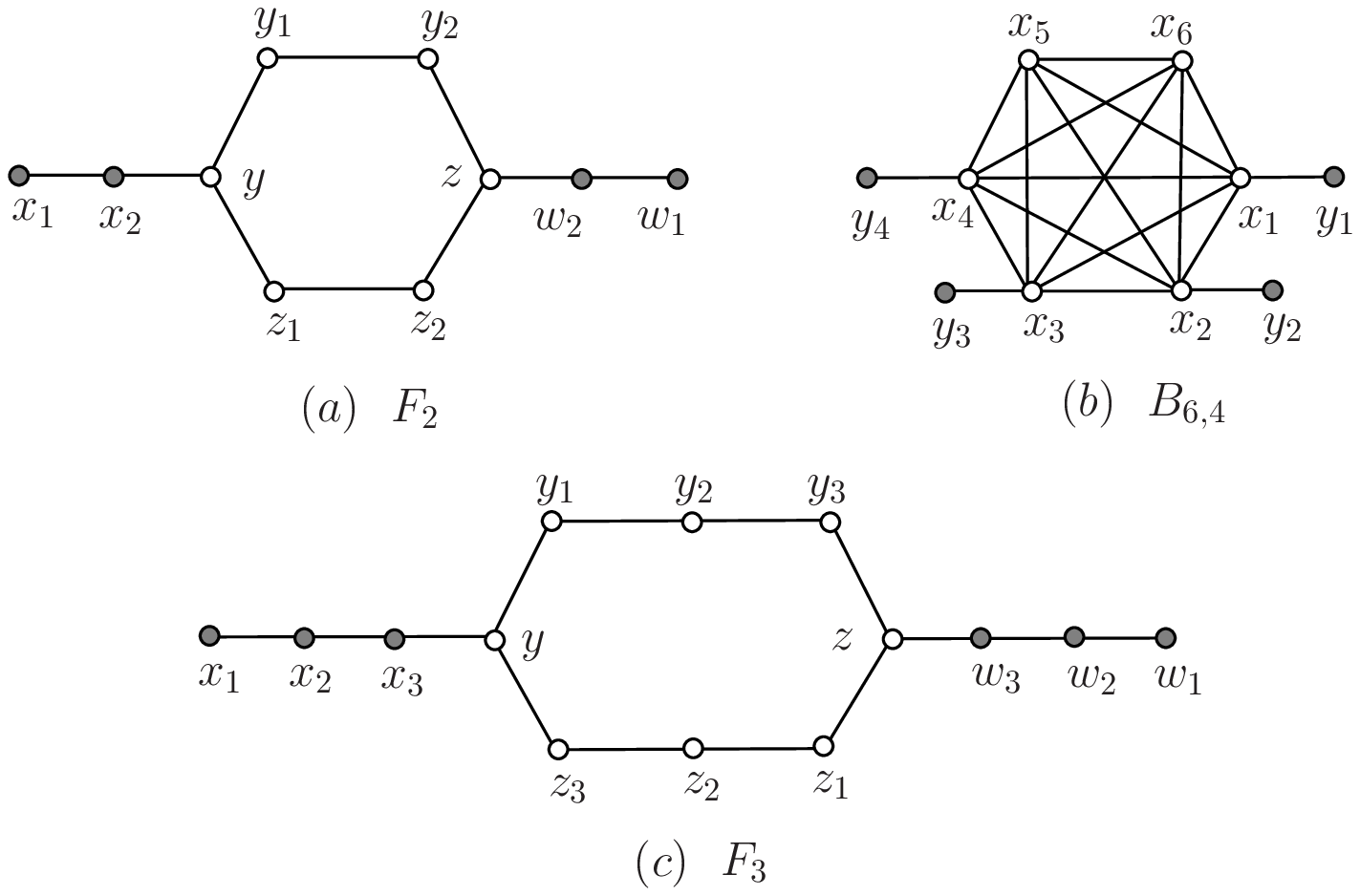}\\
\end{center}
\begin{center}
Figure 5.2: The graphs $F_2,F_3,B_{6,4}$.
\end{center}\label{fig4-3}
\end{figure}

Let $F_{r}$ be a graph such that
$V(F_{r}) =\bigcup_{i=1}^{r}\{y,z,x_i,y_i,z_i,w_i\}$ and
$$
E(F_{r}) =\bigcup_{i=1}^{r}\{w_iw_{i+1}, x_i x_{i+1}, y_iy_{i+1}, z_iz_{i+1}, yx_r, yy_1, yz_r, zy_r, zz_1, zw_r\}.
$$
For positive integers $a>b$, let $B_{a,b}$ be a block graph with $V(B_{a,b})=\{x_1,x_2,\ldots, x_a\}\cup \{y_1, y_2,\ldots, y_b\}$ and $E(B_{a,b})=\{x_ix_j\,|\,1\leq i < j\leq a\}\cup \{x_{s}y_{s}\,|\,1\leq s\leq b\}$.
As examples, graphs $F_2,F_3,B_{6,4}$ are depicted in Figure 5.2.

\begin{pro}{\upshape \cite{YehChiangPeng}}\label{pro5-2}
$(1)$ For each positive integer $r$, $C_2(H_r)$ and $C_4(H_r)$ are disjoint.

$(2)$ For positive integers $a>b\geq 2$, $C_b(B_{a,b})\nsubseteq C_{b+1}(B_{a,b})$.
\end{pro}

\subsection{Results for Steiner median}

It was first shown by Slater \cite{Slater} and later by Miller \cite{Miller},
using a more efficient construction, that every graph is the
median of some graph. However, it is unknown whether
every graph is the Steiner $k$-median of some graph for $k\geq 3$. For
trees, considerably more is known.

With $e=uv$ being a fixed edge of a tree $T$, $T_u$ and $T_v$ will denote the components of
$T-e$ containing $u$ and $v$, respectively. Also, we assume $k\geq 2$. Recall that Steiner $k$-distance of a vertex $v$ in a connected graph $G$ on $n\geq k$ vertices is defined by
$d_k(v)=\{d_G(S)\,|\,S\subseteq V(G),\  v\in S, \ |S|=k\}$.

Beineke, Oellermann, and Pippert \cite{BeinekeOellermannPippert} derived
the following result.
\begin{thm}{\upshape \cite{BeinekeOellermannPippert}}\label{th5-3}
Let $T$ be a tree of order $n\geq k$ and let $e=uv$ he an edge of $T$. If $T_u$ and $T_v$
have orders $r$ and $s$, respectively, then
$$
d_k(v)-d_k(u)=\binom{r-1}{k-1}-\binom{s-1}{k-1}.
$$
\end{thm}

Observe that the difference $d_k(v)-d_k(u)=\binom{r-1}{k-1}-\binom{s-1}{k-1}$ is $0$ if and only if $r=s$
or both $r<n$ and $s<n$. Therefore, they \cite{BeinekeOellermannPippert} had the following corollary.
\begin{cor}{\upshape \cite{BeinekeOellermannPippert}}\label{cor5-2}
$(1)$ Let $T$ be a tree of order $n\geq k$ and let $e=uv$ be an edge of $T$. Let $r$ and $s$ be
the orders of $T_u$ and $T_v$, respectively. Then $d_k(v)<d_k(u)$ if and only if $r>s$ and $r\geq k$.

$(2)$ Let $T$ be a tree of order $n\geq k$. If $v_0v_1\ldots v_r$ is a path in $T$ and
$d_k(v_0)<d_k(v_1)$, then $d_k(v_1)<d_k(v_2)< ... <d_k(v_r)$.

$(3)$ The $k$-median of any tree is connected.

$(4)$ Let $T$ be a tree of order $n>k$. If $T$ has no edge $e=uv$ for which the orders
of $T_u$ and $T_v$, either $(i)$ are equal or $(ii)$ are both less than $k$, then $M_k(T)$ is a single vertex.

$(5)$ If $T$ is a tree of order at least $2k-1$, then its $k$-median is either $K_1$ or $K_2$.
\end{cor}

In \cite{BeinekeOellermannPippert}, the Steiner $k$-medians of
trees are characterized and a linear algorithm for finding the
Steiner $k$-median of a tree is developed.

One of the main results in Corollary \ref{cor5-2} leading to a characterization of
the Steiner $k$-medians of trees states that if $T$ is a tree of order $n\geq k$
and if $v_0,v_1,\ldots,v_{r}$ is a path in $T$ with $d_k(v_0)<d_k(v_1)$, then $d_k(v_1)< d_k(v_2)<\cdots <d_k(v_r)$.

From this result, it follows that the Steiner $k$-median of a tree is
connected and thus a tree. The Steiner $k$-medians of trees are
characterized in \cite{BeinekeOellermannPippert} as follows:
\begin{thm}{\upshape \cite{BeinekeOellermannPippert}}\label{th5-4}
$(1)$ A tree $H$ is the Steiner $k$-median, $k\geq 3$, of some tree of order at
least $2k-1$ if and only if $H\cong K_1$ or $K_2$.

$(2)$ A tree $H$ is a Steiner $k$-median of some tree of order at most $2k-2$ if and only if $H$ is $K_1$ or has exactly $k$ vertices or if
$k_1$ is the number of end-vertices of $H$ and $m$ is the
number of internal vertices of $H$; then, $m+2k_1-1\leq k$.
\end{thm}

Beineke, Oellermann, and Pippert \cite{BeinekeOellermannPippert} obtained the following algorithm for finding the Steiner $k$-median of a tree.

\noindent {\bf Algorithm 5.2.}  Finding the Steiner $k$-median of a tree $T$ of
order $n\geq k\geq 2$.

$(1)$ Construct a digraph $D$ having $T$ as underlying graph by
replacing each edge $e=uv$ by

\begin{itemize}
\item[] $(a)$ The arc $(u,v)$ if $T_v$ has at least $k$ vertices and the
number of vertices in $T_v$ exceeds the number of
vertices in $T_u$.

\item[] $(b)$ The arc $(v,u)$ if $T_u$ has at least $k$ vertices and the
number of vertices in $T_u$ exceeds the number of
vertices in $T_v$.

\item[] $(c)$ The symmetric pair of arcs $(u,v)$ and $(v,u)$ if $T_u$ and
$T_v$ have the same number of vertices or if $T_u$ and $T_v$
both have fewer than $k$ vertices.
\end{itemize}

$(2)$ If $D$ has symmetric pairs of arcs, then the subgraph $H$ of
$T$ induced by those edges corresponding to the symmetric pairs of arcs in $D$ are output since this is the Steiner $k$-median
of $T$; otherwise, let $H$ be the vertex with outdegree $0$.
Output $H$ since it is the Steiner $k$-median of $T$ and stop.

Figure 5.3 illustrates Algorithm 5.2 with $k=4$ and $5$. The
Steiner $k$-median is induced by the gray vertices.
\begin{figure}[!hbpt]
\begin{center}
\includegraphics[scale=0.9]{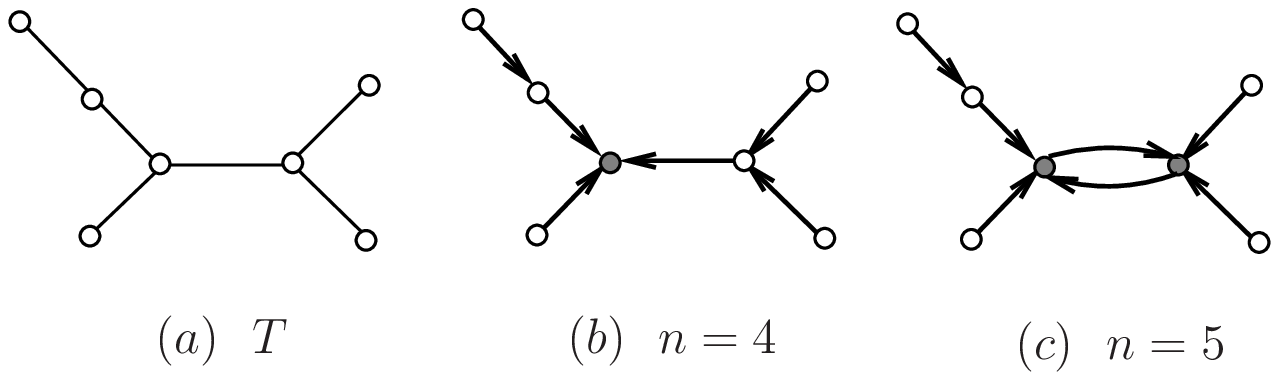}\\[0.5cm]
Figure 5.3: Graphs for Algorithm 5.2.
\end{center}\label{fig7}
\end{figure}
\begin{thm}{\upshape \cite{BeinekeOellermannPippert}}\label{th5-5}
If $T$ is a tree of order $n\geq k+1\geq 3$, then $M_k(T)\subseteq M_{k+1}(T)$.
\end{thm}

Yeh, Chiang, and Peng \cite{YehChiangPeng} showed that the containment relation between the Steiner $k$-median and Steiner $(k+1)$-median is also true for block graphs.
\begin{thm}{\upshape \cite{YehChiangPeng}}\label{th5-6}
If $G$ is a block graph with $|V(G)|\geq k+1\geq 3$, then
$$
M_{k}(T)\subseteq M_{k+1}(T).
$$
\end{thm}

It would be interesting to ask if Theorem \ref{th5-6} holds for distance-hereditary graphs, particularly since the distance-hereditary graphs are also Steiner distance hereditary. Yeh, Chiang, and Peng \cite{YehChiangPeng} showed that for each $k\geq 2$ there is a
distance-hereditary graph $G$ such that $M_k(G)\nsubseteq M_{k+1}(G)$.
\begin{figure}[!hbpt]
\begin{center}
\includegraphics[scale=0.8]{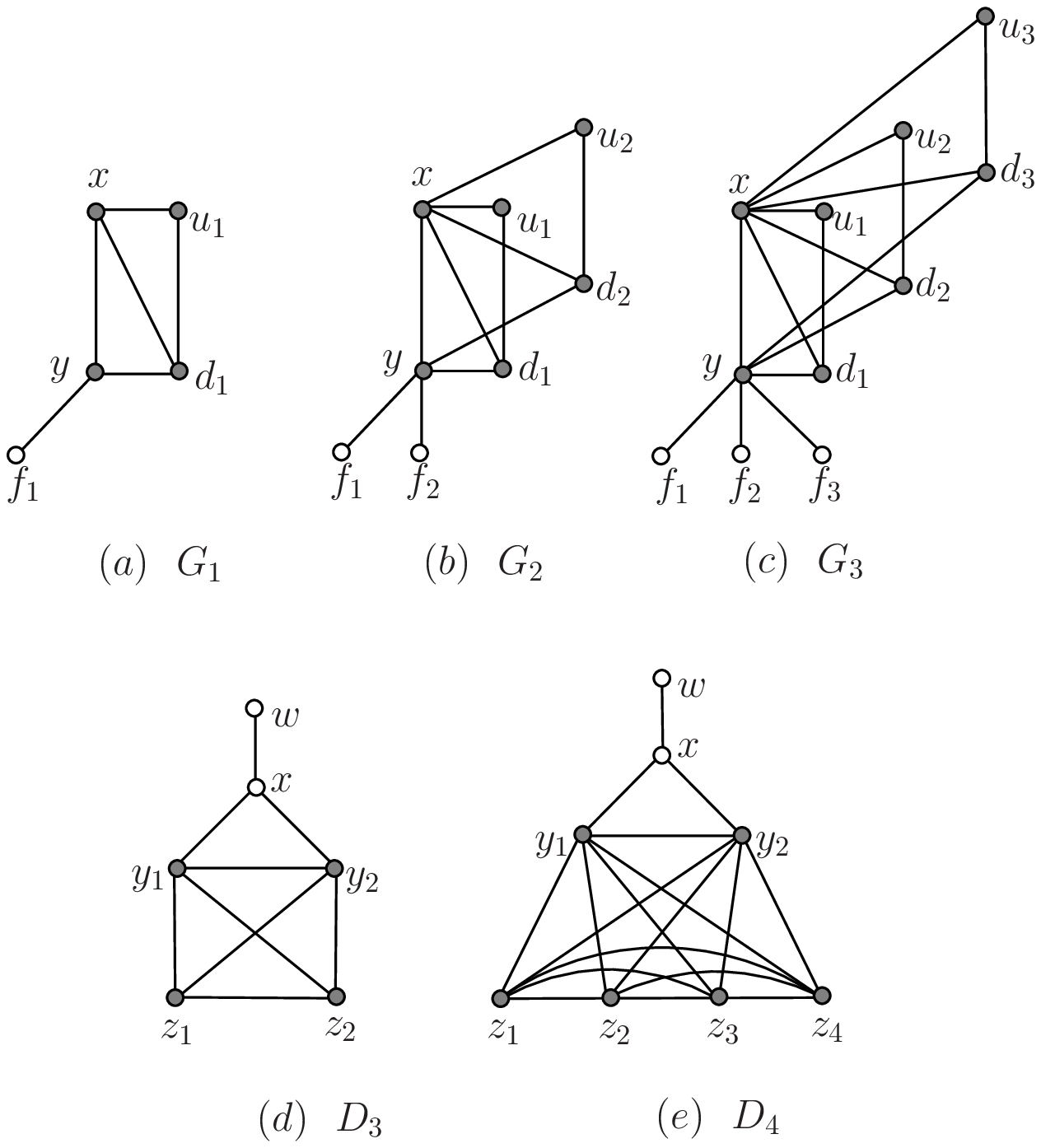}\\[0.5cm]
Figure 5.4: The graphs $G_i \ (1\leq i\leq 3)$ and $D_i \ (i=3,4)$.
\end{center}\label{fig7}
\end{figure}

Let $G_r$ be a graph with $V(G_r)=\bigcup_{i=1}^r\{x,y,u_i,d_i,f_i \}$
and
$$
E(G_r)=\bigcup_{i=1}^r\{xy,xu_i,xd_i,yd_i,yf_i,u_id_i\}.
$$
For $k\geq 3$, let $D_r$ be a graph with $V(D_r)=\{w,x,y_1,y_2\}\cup \{z_1,z_2,\ldots,z_{2r-4}\}$ and $E(D_r)=\{wx, xy_1, xy_2, y_1y_2\}\cup \{y_iz_j\,|\,1\leq i\leq 2~{\rm and}~1\leq j\leq 2r-4\}\cup \{z_iz_j\,|\,1\leq i < j\leq 2r-4\}$.
It is easy to see that $G_1$ and $D_r$, $r\geq 3$, are distance-hereditary graphs \cite{Howorka}. As examples, graphs $G_i \ (1\leq i\leq 3)$ and $D_i \ (i=3,4)$ are depicted in Figure 5.4.

\begin{thm}{\upshape \cite{BeinekeOellermannPippert}}\label{th5-5}
If $T$ is a tree of order $n\geq k+1\geq 3$, then $M_k(T)\subseteq M_{k+1}(T)$.
\end{thm}

\begin{pro}{\upshape \cite{YehChiangPeng}}\label{pro5-3}
$(1)$ For each positive integer $r\geq 2$,$M_2(G_r)\nsubseteq M_3(G_r)$.

$(2)$ For each distance-hereditary graph
$D_r$, $r\geq 3$, $M_r(D_r)\nsubseteq M_{r+1}(D_r)$.
\end{pro}

Yeh, Chiang, and Peng \cite{YehChiangPeng} derived the following result, and presented
a linear time algorithm for finding the Steiner $k$-median of a block graph.
\begin{thm}{\upshape \cite{YehChiangPeng}}\label{th5-7}
Let $G$ be a block graph with $|V (G)|\geq k\geq 2$ such that $G$ is not a complete graph. If vertex $x$ is not a cut-vertex of $G$, then either $x$ is not a vertex in the Steiner $k$-median of $G$, or $|V(G)|=k$.
\end{thm}

In the following algorithm we use the following notation: $\overrightarrow{G_c}$ has two kinds of arcs $u\leftarrow v$ and $u\leftrightarrow v$. For an arc
of the form $u\leftarrow v$ (resp. $u\rightarrow v$), we say that it is an \emph{out-edge} of vertex $v$ (resp. $u$). For the purpose of this algorithm
we assume that an arc of the type $u\leftrightarrow v$ is not an out-edge of $u$ or $v$.

\noindent {\bf Algorithm 5.3.}  Finding the Steiner median $(G,k)$.

{\bf Input:} A block graph $G$ with $|V(G)|\geq k\geq 2$.

{\bf Output:} The Steiner $k$-median $M_k(G)$ of $G$.

{\bf Begin}
\begin{itemize}
\item[] 1. if $n=|V(G)|$ or $G$ is a complete graph,
then $M_k(G)\leftarrow G$ and STOP.

\item[] 2. Let $G_c$ be the subgraph of $G$ induced by the cut-vertices of $G$.

\item[] 3. Construct a graph $\overrightarrow{G_c}$ with two kinds of arcs from $G_c$ as follows:

\item[] 4. for each edge $uv$ in $G_c$,

\ \ if $|V_{uv}|>|V_{vu}|$ and $|V_{uv}|+|V_{uv}^{=}|\geq k$,

\ \ \ \ then replace the edge $uv$ by an arc of the form $u\leftarrow v$,

\ \ \ \ \ else if $|V_{uv}|<|V_{vu}|$ and $|V_{vu}|+|V_{uv}^{=}|\geq k$,

\ \ \ \ \ \ then replace the edge $uv$ by an arc of the form $u\rightarrow v$,

\ \ \ \ \ \ else replace the edge $uv$ by an arc of the form $u\leftrightarrow v$.

\item[] 5. Let $M_k(G)$ be those vertices in $\overrightarrow{G_c}$ with no out-edges.
\end{itemize}
{\bf End.}

\begin{thm}{\upshape \cite{YehChiangPeng}}\label{th5-7}
Let $G=(V,E)$ be a block graph with $|V|\geq k\geq 2$. Algorithm 5.3 correctly
finds the Steiner $k$-median of $G$ in time $O(|V|+|E|)$.
\end{thm}

A block of a block graph $G$ that contains exactly one cut-vertex is called an \emph{end-block} of $G$.

Yeh, Chiang, and Peng \cite{YehChiangPeng} showed that one can easily find the Steiner $k$-distance of all vertices in a block graph in polynomial time.
\begin{thm}{\upshape \cite{YehChiangPeng}}\label{th5-7}
Suppose $x$ is a vertex in an end-block of a block graph $G$ with $|V(G)|\geq k\geq 2$, and $x$ is not a cut-vertex
of $G$. Let $y$ be a vertex of $G$ adjacent to $x$. Let $d_k'(y)$ (resp. $d_k(y)$) denote the Steiner $k$-distance of $y$ in $G-x$ (resp. $G$). Then
$$
d_k(y)=d_k'(y)+d_{k-1}'(y)+{|V(G)|-2 \choose k-2}.
$$
\end{thm}

Beineke, Oellermann, and Pippert \cite{BeinekeOellermannPippert} turned
their attention to the minimum value of the Steiner $k$-distances of the vertices of a graph.    In particular, for a connected graph $G$ of order $n\geq k$, let the Steiner $k$-median    value $m_k(G)$ be defined as $\min\{d_k(v)\,|\,v\in V(G)\}$. Sharp
bounds on the Steiner $k$-median value of trees of order $n\geq 2k-1$ are also established in \cite{BeinekeOellermannPippert}.
\begin{thm}{\upshape \cite{BeinekeOellermannPippert}}\label{th5-10}
If $T_n$ is a tree of order $n\geq 2k-1$, then
\begin{eqnarray*}
(k-1)\binom{n-1}{k-1}&\leq &m_k(T_n)\\[0.1cm]
&\leq&\left\{
\begin{array}{ll}
(n-1)\binom{n-1}{k-1}-2\sum_{j=1}^{(n-1)/2}\binom{n-j-1}{k-1},&\mbox{{\rm if}~$n$~{\rm is~odd};}\\[0.2cm]
(n-1)\binom{n-1}{k-1}-2\sum_{j=1}^{(n-2)/2}\binom{n-j-1}{k-1}-\binom{\frac{1}{2}n-1}{k-1},&\mbox{{\rm if}~$n$~{\rm is~even}.}
\end{array}\right.
\end{eqnarray*}
Furthermore, these bounds are sharp.
\end{thm}

\subsection{From Steiner centers to Steiner medians}

In the preceding two subsections, the focus has been on finding
the Steiner $k$-centers and Steiner $k$-medians of trees. In \cite{Oellermann}, it was shown
that, except for trees of small order, these two types of
``centers'' can be arbitrarily far apart.

Let $T$ be a tree of order $n \ (2\leq k\leq n)$. Suppose that $T$ has
at most $2k-2$ vertices. If $T$ had at most $k-1$ end-vertices, then by Algorithm 5.1, $C_k(T)=T$.  So, in this
case, $M_k(T)\subset C_k(T)$. Suppose that $T$ has at least $k$
end-vertices. Then, by Algorithm 5.1, $C_k(T)$ is obtained by
deleting the end-vertices of $T$. By Algorithm 5.2, $M_k(T)$
does not contain any end-vertex of $T$. Hence, once again,
$M_k(T)\subset C_k(T)$.

Oellermann \cite{Oellermann} turned her attention to trees having at least $2k-1$
vertices. For a given graph $G$ and subgraphs $F$ and $H$ of $G$,
the \emph{distance $d_G(F,H)$ between $F$ and $H$} is defined as $\min
\{d_G(u, v)\,|\,u\in V(F)~{\rm and}~v\in V(H)\}$.

Hendry \cite{Hendry} showed
that if $F$ and $H$ are any two graphs then there exists a
connected graph $G$ such that $C(G)\cong F$ and $M(G)\cong H$.
The graph constructed by Hendry \cite{Hendry} had the property that
$d_G(F,H)=1$. Holbert \cite{Holbert} showed that the distance between $C(G)$ and $M(G)$ can  be arbitrarily large. It was
shown in \cite{Oellermann} that the distance between the Steiner $k$-center and
Steiner $k$-median of a tree of sufficiently large order can be arbitrarily large. Moreover, it is shown that the structure of the
Steiner $k$-center and Steiner $k$-median can be prescribed
provided that they satisfy the conditions of Theorem \ref{th5-2}
and $(1)$ of Theorem \ref{th5-4}.
\begin{thm}{\upshape \cite{Oellermann}}\label{th5-11}
Let $T_1$ be any tree with at most $k-1$
end-vertices and let $T_2$ be isomorphic to $K_1$ or $K_2$. Let $d\geq 1$
be an integer. Then, there exists a tree $T$ with $C_k(T)=T_1$,
$M_k(T)=T_2$ and $d(T_1,T_2)=d$.
\end{thm}

Since the Steiner $k$-center and Steiner $k$-median of a
graph are both measures of centrality, it seems reasonable to
ask whether there are measures of centrality that allow each
vertex on a shortest path between the Steiner $k$-center and
Steiner $k$-median to belong to the ``center'' with respect to at
least one of these measures.

For $1\leq k\leq |V(G)|$, and $u\in V(G)$, Slater \cite{Slater} defined
$$
r_k(u)=\max\left\{\sum_{s\in S}
d(u,s)\,|\,S\subseteq V(G), \ |S|=k\right\}.
$$

The \emph{$t$-centrum} of $G$, denoted by $C(G;t)$, is defined to
be the subset of vertices $u$ in $G$ for which $r_t(u)$ is a
minimum. Thus, the vertices of $C(G;1)$ induce the center
and the vertices of $C(G;|V(G)|)$ induce the median of $G$.

Slater \cite{Slater} showed that for a tree $T$, if $u$ belongs to the
center and $v$ to the median of $T$, then the subgraph induced
by the vertices in $\bigcup_{t=1}^{|V(T)|}C(T;t)$ is a subtree of $T$ containing
the $u$-$v$ path.

Oellermann \cite{Oellermann2} turned her attention to general $k$.
Let $G$ be a connected graph of order $n \ (2\leq k\leq n)$ and let $\mathscr{P}_k$ be the
collection of all $k$-element subsets of $V(G)$. For $1\leq t\leq \binom{n-1}{k-1}$ and $u\in V(G)$, let the Steiner $(k,t)$-eccentricity of
$u$ be defined by
$$
e_G^{(k)}(u;t)=\max\left\{\sum_{s\in S}d(S)\,|\,S\subseteq \mathscr{P}_n, \ |S|=k~{\rm and}~u\in S~{\rm for~all}~s\in S\right\}.
$$

The Steiner $(k,t)$-center $C_k(G;t)$ is the subgraph
induced by the vertices of $G$ with minimum Steiner $(k,t)$-eccentricity. Thus, if $t=1$, then $C_k(G;t)=C_k(G)$, and if $t=\binom{n-1}{k-1}$, then $C_k(G;t)=M_k(G)$. The following result was established in \cite{Oellermann}:
\begin{thm}{\upshape \cite{Oellermann}}\label{th5-12}
If $T$ is a tree of order $n\geq k$ and $w$ is a vertex
on a shortest path from $C_k(T)$ to $M_k(T)$, then $w$ is a vertex of
$C_k(T;t)$ for some $t$, $1\leq t\leq \binom{n-1}{k-1}$.
\end{thm}

Whether or not Theorem \ref{th5-11} can be extended to general
graphs or even to other classes of graphs for which the
Steiner problem can be solved efficiently is still an open
problem.

\section{Steiner Distance Hereditary Graphs}

To be able to state the characterizations of
distance hereditary graphs given by Howorka \cite{Howorka}, we need the following terminology: An {\it induced path\/} of $G$
is a path that is an induced subgraph of $G$. Let $u,v\in V(G)$. Then, a {\it $u$-$v$ geodesic\/} is a shortest $u$-$v$ path.
Let $C$ be a cycle of $G$. A path $P$ is an {\it essential part\/} of $C$ if $P$ is a subgraph of $C$ and $\frac{1}{2}e(C)<e(P)<e(C)$. An edge of $G$ that joins two vertices of $C$ that are not adjacent in $C$ is called
a {\it diagonal\/} of $C$. We say that two diagonals $e_1,e_2$ of $C$ are {\it skew diagonals\/} if $C+e_1+e_2$ is homeomorphic with $K_4$.
\begin{thm}{\upshape \cite{Howorka}}\label{th6-1}
The following are equivalent:
\begin{itemize}
\item[] $(1)$ $G$ is distance hereditary;

\item[] $(2)$ Every induced path is a geodesic;

\item[] $(3)$ No essential part of a cycle is induced;

\item[] $(4)$ Each cycle of length at least $5$ has at least two diagonals and each $5$-cycle has a pair of skew diagonals;

\item[] $(5)$ Each cycle of $G$ of length at least $5$ has a pair of skew diagonals.
\end{itemize}
\end{thm}

In \cite{DayOellermannSwart}, it was pointed out that, in general, it appears to be a
difficult problem to determine the Steiner distance of a set of vertices in a graph. However, it was shown in \cite{DayOellermannSwart} that if $G$ is $k$-Steiner distance hereditary, then the Steiner distance of any set of $k$ vertices of $G$ can be determined efficiently. Furthermore, the following result is established in \cite{DayOellermannSwart}:
\begin{thm}{\upshape \cite{DayOellermannSwart}}\label{th6-2}
If $G$ is $2$-Steiner distance hereditary, then $G$ is $k$-Steiner distance hereditary for all $k\geq 3$.
\end{thm}

A vertex $v$ in a
graph $G$ is a \emph{true (false) twin} of a vertex $v'$ if $v$ and $v'$ have the same closed (respectively, open) neighborhood in $G$. A
\emph{twin} of a vertex $v$ is a vertex that is either a true or false twin of $v$.

Bandelt and Mulder \cite{BandeltMulder}, Hammer and Maffray \cite{HammerMaffray} derived the following result for distance hereditary graphs.
\begin{thm}{\upshape \cite{BandeltMulder, HammerMaffray}}\label{th6-3}
A graph $G$ of order $n$ is distance hereditary if and
only if there is a sequence of subgraphs $G_1,G_2,\ldots,G_{n-1}$ such that $G_1\cong K_2$ and for $2\leq i\leq n-1$, $G_i$ is obtained
by adding a new vertex $v$ as pendant, or twin of some vertex $v'$ of $G_{i-1}$.
\end{thm}

\subsection{Steiner geodetic sets of distance-hereditary graphs}

In general there is no relation between $g(G)$ and $sg(G)$. To see that $sg(G)$ can be much larger than $g(G)$ observe
that $g(K_{m,n})=4$ and $sg(K_{m,n})=m$ for $5\leq m\leq n$. The {\it convex hull\/} of a set $S$ of vertices (that is not necessarily convex) is the
smallest convex set of vertices that contains $S$. A vertex $v$ of a convex set $S$ is an extreme point of $S$ if $S\setminus \{v\}$ is still convex.

Oellermanna and Puertas \cite{OellermannPuertas} showed that if $G$ is a distance-hereditary graph, then the geodetic number never exceeds the Steiner geodetic
number.
\begin{thm}{\upshape \cite{OellermannPuertas}}\label{th6-4}
If $G$ is a distance-hereditary graph, then $g(G)\leq sg(G)$.
\end{thm}

A vertex is a \emph{contour vertex} if its eccentricity is at least as large
as the eccentricity of each of its neighbors. The collection of all contour vertices of $G$ is called its \emph{contour} and is
denoted by $Ct(G)$.

The next result is used to find minimum Steiner geodetic sets in distance-hereditary graphs.
\begin{lem}{\upshape \cite{OellermannPuertas}}\label{lem6-1}
Let $G$ be a distance-hereditary graph and $T$ a Steiner tree for $Ct(G)$. Then every vertex of $G$ not in $T$ is adjacent to some vertex of $T$.
\end{lem}

In \cite{OellermannPuertas}, Oellermanna and Puertas stated a result that describes a procedure for finding minimal Steiner geodetic sets in all distance-hereditary
graphs.
\begin{thm}{\upshape \cite{OellermannPuertas}}\label{th6-5}
Let $G$ be a distance-hereditary graph and suppose that $\overline{I}=V(G)\setminus I(Ct(G))$. Let $S=Ct(G)\cup \overline{I}$. Then
$S$ is a Steiner geodetic set for $G$ that is minimal in the sense that for all $v\in S$, $v\notin I(S-v)$.
\end{thm}

Except for some special cases, the Steiner geodetic set described in Theorem \ref{th6-5}, is a minimum Steiner geodetic set
as we now see.
\begin{thm}{\upshape \cite{OellermannPuertas}}\label{th6-6}
Suppose $G$ is a distance-hereditary graph and $S$ a minimum Steiner geodetic set for $G$. If $diam(G)\geq 3$, then $S$ contains all contour vertices of $G$.
\end{thm}

Let $G$ be a distance-hereditary graph. If $G$ has diameter at least $3$ it follows, from Theorem \ref{th6-6}, that a minimum
Steiner geodetic set $S$ contains $Ct(G)$. Let $\overline{I}=V(G)\setminus I(Ct(G))$. If $\overline{I}\nsubseteq S$, let $u$ be a vertex in $\overline{I}$ that does not belong to $S$.
Then there is a Steiner tree $T$ for $S$ that contains $u$. Since $G$ is Steiner distance hereditary, the induced subgraph $<V(T)>$
must contain a Steiner tree $H$ for $Ct(G)$ that does not contain $u$. Let $W=S\cap \overline{I}$. Then, by Lemma \ref{lem6-1}, every vertex of
$W$ is adjacent with some vertex in $H$. So the tree obtained from $H$ by adding an edge from every vertex in $W$ to some
vertex in $H$ is a tree that contains $S$ but not $u$ and thus has size less than the size of $T$. This is not possible. So every
vertex in $Ct(G)\cup \overline{I}$ belongs to every minimum Steiner geodetic set. By Theorem \ref{th6-5}, $I(Ct(G)\cup \overline{I})=V(G)$. Thus,
$Ct(G)\cup \overline{I}$ is a minimum Steiner geodetic set. If $G$ has diameter $1$, then $\overline{I}=\emptyset$ and a Steiner geodetic set for $G$ must
contain all (contour) vertices. So in this case $Ct(G)\cup \overline{I}$ is a minimum Steiner geodetic set. Similarly if $diam(G)=2$
and $rad(G)=1$, $Ct(G)\cup \overline{I}$ is a minimum Steiner geodetic set.
The next theorem summarizes these results.
\begin{thm}{\upshape \cite{OellermannPuertas}}\label{th6-7}
If $G$ is a distance-hereditary graph with

$(i)$ $diam(G)=1$; or

$(ii)$ $diam(G)=2$ and $rad(G)=1$; or

$(iii)$ $diam(G)\geq 3$;

then $sg(G)=|Ct(G)|+|V(G)\setminus I(Ct(G))|$ and $Ct(G)\cup (V(G)\setminus I(Ct(G)))$ is the unique minimum Steiner geodetic set.
\end{thm}

\subsection{$3$-Steiner distance
hereditary graphs}

Several new characterizations of $2$-Steiner distance hereditary graphs which lead to efficient algorithms that test
whether a graph is $2$-Steiner distance hereditary have been established; see \cite{BandeltMulder, DAtriMoscarini, HammerMaffray}.

A structural characterization of $3$-$SDH$ graphs is given in \cite{DayOellermannSwart2}. Suppose $C: v_1, v_2, \ldots,$ $v_{\ell}, v_1$ is a cycle in a graph $G$.
An edge of $G$ that joins two vertices of $C$ that are not adjacent on $C$ is called a \emph{diagonal} or a \emph{chord} of $C$. Two chords
$e_1$ and $e_2$ of $C$ are skew or crossing, if $C+e_1+e_2$ is homeomorphic to $K_4$.
\begin{thm}{\upshape \cite{DayOellermannSwart2}}\label{th6-8}
A graph $G$ is $3$-Steiner distance hereditary if and only if it is $2$-Steiner distance hereditary or if the following conditions hold:
\begin{itemize}
\item[] $(1)$ Every cycle $C$ : $v_1,v_2,\ldots, v_{\ell},v_1$ of length $\ell\geq 6$

\item[] \ \ \ \ \ $(a)$ has at least two skew diagonals, or, if $\ell=6$, then
$v_1v_3v_5v_1$ or $v_2v_4v_6v_2$ is a cycle in $<V(C)>$;

\item[] \ \ \ \ \ $(b)$ has no two adjacent vertices neither of which is on a diagonal of $C$.

\item[] $(2)$ $G$ does not contain an induced subgraph isomorphic to any of the following
(any subset of dotted edges may be included in the graph).
\end{itemize}
\end{thm}

A \emph{hole} is an induced
cycle of length at least $5$; a \emph{house} is a $5$-cycle with exactly one chord and a \emph{domino} is a $6$-cycle with exactly one chord
that joins two vertices distance $3$ apart on the cycle. The $HHD$-free graphs are characterized as those graphs for which
every cycle of length at least $5$ contains at least two chords; see \cite{Brandstadt}.

\begin{thm}{\upshape \cite{DayOellermannSwart2}}\label{th6-9}
$(1)$ If $G$ is a $3$-$SDH$ graph, then the contour of $G$ is a geodetic set.

$(2)$ If $G$ is a $HHD$-free graph, then the contour of $G$ is a geodetic set.
\end{thm}

We show here that the Steiner geodetic number is an upper bound for the geodetic number for $3$-$SDH$ graphs.
\begin{thm}{\upshape \cite{DayOellermannSwart2}}\label{th6-10}
If $G$ is a $3$-$SDH$ graph and $S\subseteq V(G)$, then $I(S)\subseteq I[S]$.
\end{thm}

\begin{cor}{\upshape \cite{DayOellermannSwart2}}\label{cor6-1}
If $G$ is a $3$-$SDH$ graph, then $g(G)\leq sg(G)$.
\end{cor}

We say that a $6$-cycle $v_1,v_2,v_3,v_4,v_5,v_6,v_1$
has a triangle of chords if $v_1,v_3,v_5, v_1$ or $v_2, v_4,v_6,v_2$ is a cycle in $G$.

Oellermann and Spinrad \cite{OellermannSpinrad} obtained the following polynomial
algorithm to test whether a graph is $3$-$SDH$.

\noindent {\bf Algorithm 6.1.} {\upshape \cite{OellermannSpinrad}} To test whether a given graph $G$ is $3$-$SDH$
\begin{itemize}
\item[] $(1)$ If $G$ is distance hereditary, then output ``$G$ is
$3$-$SDH$'' and halt;

\item[] $(2)$ If $G$ has a bad edge, then output ``$G$ is not $3$-$SDH$'' and halt;

\item[] $(3)$ If $G$ has a $7$-cycle without crossing chords or if $G$ has a $6$-cycle without  crossing chords and without a triangle of chords, then output ``$G$
is not $3$-$SDH$'' and halt;

\item[] $(4)$ If $G$ has any of the forbidden subgraphs of
Fig. 1 as an induced subgraph, then output
``$G$ is not $3$-$SDH$'' else output ``$G$ is $3$-$SDH$''.
\end{itemize}

Goddard \cite{Goddard} showed that if a graph is $k$-$SDH$, then it is $t$-$SDH$ for all $t\geq k$. The following algorithm, see \cite{DayOellermannSwart}, finds
the Steiner distance of a set $S$ of vertices in a $k$-$SDH$ graph. We use this algorithm to find the Steiner interval for $S$. We
say that a set $S$ of vertices of a graph $G$ is separated in an induced subgraph $H$ of $G$ that contains $S$ if the vertices of $S$
do not belong to the same component of $H$.

\noindent {\bf Algorithm 6.2.} {\upshape \cite{DayOellermannSwart}}
Algorithm to find the Steiner distance of a set $S$, of at least three vertices, in a $3$-$SDH$ graph $G$.
\begin{itemize}
\item label $V(G)\backslash S$ in arbitrary order $v_1,v_2,\ldots,v_m$

\item $G_1=G$

\item for $i=1$ to $m$

if $S$ is separated in $G_i-v_i$

then $G_{i+1}\longleftarrow G_i$

else $G_{i+1}\longleftarrow G_i-v_i$

\item $d_G(S)=|V(G_{m+1})|-1$
\end{itemize}

Eroh and Oellermann \cite{ErohOellermann} derived the following result.
\begin{thm}{\upshape \cite{ErohOellermann}}\label{th6-10}
If $G$ is a $3$-$SDH$ graph and if $v_m$, in the final step of the algorithm above, does not separate $S$ in $G_m$, then
$v_m\notin I(S)$.
\end{thm}

Thus, in order to find $I(S)$ in a $3$-$SDH$ graph, we may perform the given algorithm $m$ times, with each of the vertices
of $V(G)\setminus S$ in turn as the last vertex $v_m$ in the sequence of vertices input in the algorithm. If $v_m$ does not separate $S$ in
$G_m$, then it is not in $I(S)$; otherwise, it is.

\subsection{$k$-Steiner-distance-hereditary}

For $k\geq 2$ and $d\geq k-1$, Goddard \cite{Goddard} defined the property $S(k,d)$ as meaning
that for all sets $S$ of $k$ vertices with Steiner distance $d$, the distance of $S$ is
preserved in any connection for $S$. The property $S(k)$ as meaning
$S(k,d)$ for all $d$ is also defied in \cite{Goddard}.
Day, Oellermann, and Swart \cite{DayOellermannSwart} introduced the property and called such
graphs $k$-Steiner-distance-hereditary. Thus distance-hereditary graphs are
the ones obeying $S(2)$.

Day, Oellermann, and Swart \cite{DayOellermannSwart} conjectured
being $k$-Steiner-distance-hereditary implies being $(k+1)$-Steiner-distance-hereditary.
Goddard \cite{Goddard} showed that there is a partial converse:
\begin{thm}{\upshape \cite{Goddard}}\label{th6-12}
$(1)$ For all $k\geq 2$ it holds that $S(k,k)$ is equivalent to $S(k)$.

$(2)$ For all $k\geq 2$ it holds that $S(k)$ implies $S(k+1)$.

$(3)$ For all $k\geq 3$ it holds that $S(k)$ implies $S(k-1,d)$ for all $d\geq k$.
\end{thm}

\section{Steiner Intervals in Graphs}

The following result in \cite{KubickaKubickiOellermann} shows how the Steiner distance of an $k$-set can be found if its
$2$-intersection interval is nonempty.
\begin{thm}{\upshape \cite{KubickaKubickiOellermann}}\label{th7-1}
Let $S=\{u_1,u_2,\ldots,u_k\}$ be a set of $k\geq 2$ vertices of a graph $G$. If the
$2$-intersection interval of $S$ is nonempty and $x\in I_2(S)$, then $d(S)=\sum_{i=1}^rd(u_i,x)$.
\end{thm}

Since it can be determined in polynomial time  whether an $k$-set satisfies the
hypothesis of Theorem \ref{th7-1}, the Steiner distance for such an $k$-set can be found in
polynomial time. The result of Theorem \ref{th7-1} is best possible in the sense that we will  now discuss. Let $G$ be the graph shown in $(1)$ of Figure 7.1 and let $S=\{u_1,u_2,\ldots,u_k\}$. Then $x$ lies   on a shortest $u_i$-$u_j$ path for all $1\leq i<j\leq n$ except on a shortest $u_{k-1}$-$u_{k}$ path, but $x$ does not belong to a Steiner tree for $S$ and $d(S)\neq \sum_{i=1}^kd(u_i,x)$. In fact, $d(S)=2k-1$ but a tree with the least number of edges and containing $S$ and $x$ has $2k$ edges and
$\sum_{i=1}^kd(u_i,x)=2k$. This graph also serves to illustrate that there may be $k$-sets   for which the $2$-intersection intervals are empty.

Kubicka, Kubicki, and Oellermann \cite{KubickaKubickiOellermann} observed also that the $r$-intersection interval of $S$ for all $r$, $3\leq r\leq k$, contains
$z$ and is thus nonempty. Indeed, it is our belief that if the $r$-intersection interval
of some $k$-set $S$ is nonempty, then the $\ell$-intersection interval of $S$ is nonempty      for all $\ell>r$ and if $x\in I_{r}(S)$, then $x\in I_{\ell}(S)$. To prove this statement,      it would suffice to show that if $S$ is an $k$-set and $x\in I_{k-1}(S)$, then there exists   a Steiner tree for $S$ that contains $x$.
\begin{figure}[!hbpt]
\begin{center}
\includegraphics[scale=0.7]{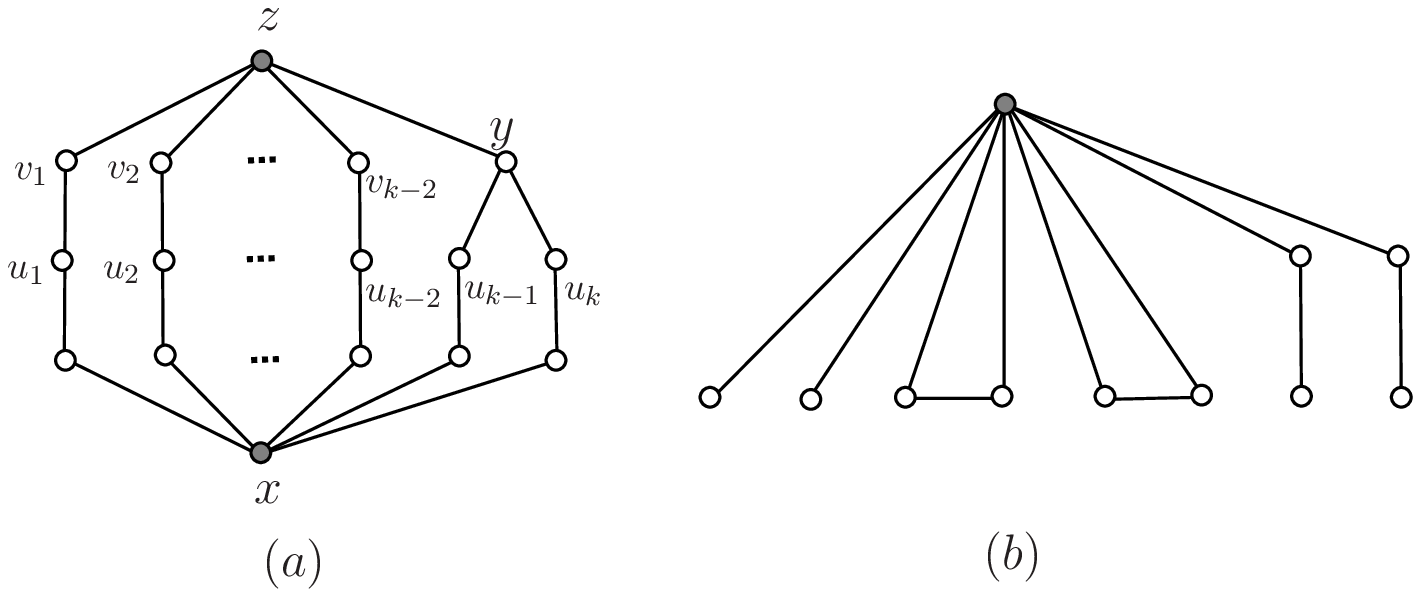}\\[0.5cm]
Figure 7.1: Graphs for Theorems \ref{th7-1} and \ref{th7-5}.
\end{center}\label{fig7}
\end{figure}

\subsection{Results for Steiner intervals}

Graphs for which the $2$-intersection interval of every $3$-set is nonempty have been
studied; see \cite{Mulder}. The only graphs which have the property that the $2$-intersection
interval of every $k$-set is nonempty, where $k>3$, are the stars. To see this, note that if
$G$ is a connected graph but not a star, then there exist two independent edges in $G$. If
we place the four end vertices of these two edges in an $k$-set, then the $2$-intersection
interval of this $k$-set is empty.

Kubicka, Kubicki, and Oellermann \cite{KubickaKubickiOellermann} investigated graphs with the property that all
$n$-sets have nonempty $3$-intersection intervals for some fixed $n\geq 4$.
\begin{thm}{\upshape \cite{KubickaKubickiOellermann}}\label{th7-2}
A graph $G$ has the property that the $3$-intersection interval of every $4$-set is
nonempty if and only if $G$ has no cycles of length other than $3$ or $5$.
\end{thm}

In order to characterize graphs with the property that the $3$-intersection intervals of
every $5$-set is nonempty, they \cite{KubickaKubickiOellermann} established the following result for trees which is of interest
in its own right.
\begin{thm}{\upshape \cite{KubickaKubickiOellermann}}\label{th7-3}
For a tree of order at least $2r-1$ the $r$-intersection interval of every
$(2r-1)$-set consists of exactly one vertex.
\end{thm}

\begin{thm}{\upshape \cite{KubickaKubickiOellermann}}\label{th7-4}
A graph $G$ has the property that the $3$-intersection interval of every $5$-set is
nonempty if and only if every block of $G$ is isomorphic to $K_2$ or $K_3$ and those blocks
isomorphic to $K_3$ are end-blocks.
\end{thm}

\begin{thm}{\upshape \cite{KubickaKubickiOellermann}}\label{th7-5}
For $k=6$, a graph $G$ has the property that the $3$-intersection intercal of
every $k$-set is nonempty if and only if its order is at least $k$ and, for some vertex $v$ every
component of $G-v$ is $K_1$ or $K_2$ (see $(2)$ of Figure 7.1).
\end{thm}

Kubicka, Kubicki, and Oellermann \cite{KubickaKubickiOellermann} characterized those graphs $G$ for which $I_r(S)$ is nonempty for every
$k$-set $S$, provided $k$ is sufficiently large in comparison to $r$. In order to present our
characterization we need the following characterization of $2$-connected graphs that
appears in \cite{DankelmannOellermannSwart} and is a stronger characterization of $2$-connected graphs given independently
by Gy\"{o}ri \cite{Gyori} and Lov\'{a}sz \cite{Lovasz}.
\begin{thm}{\upshape \cite{KubickaKubickiOellermann}}\label{th7-6}
A graph $G$ of order $n$ is $2$-connected if and only if for every two distinct
vertices $a,b$ in $V(G)$ there exists an ordering of the vertices of $G$, $a=x_1,x_2,\ldots,
x_n=b$ such that for each $\ell$ with $1\leq \ell\leq n$ the subgraphs induced by $\{x_1, x_2,\ldots, x_{\ell}\}$
and $\{x_{\ell}, x_{\ell+1},\ldots, x_{n}\}$ are connected.
\end{thm}

\begin{thm}{\upshape \cite{KubickaKubickiOellermann}}\label{th7-7}
Let $G$ be a graph of order $n\geq k$ and suppose $k\geq 2r$. Then $I_r(S)$ is nonempty
for every $k$-set $S$ of vertices of $G$ if and only if $G$ has a cut vertex $v$ such that every
component of $G-v$ has at most $r-1$ vertices.
\end{thm}

\subsection{Finding Steiner intervals in distance-hereditary graphs}

Let $G$ be a distance-hereditary graph and $S\subseteq V(G)$ with $|S|\geq 2$. Oellermann and Puertas \cite{OellermannPuertas} developed an algorithm for finding the Steiner
interval $I(S)$ of $S$ in $G$. Let $I_2(S)=\bigcup_{a,b\in S}
I[a,b]$. If $G$ is a graph and $\mathscr{S}\subseteq V(G)$, then $<\mathscr{S}>$ denotes the subgraph
induced by $\mathscr{S}$.
\begin{pro}{\upshape \cite{OellermannPuertas}}\label{pro7-1}
If $G$ is a distance-hereditary graph and $S$ a set of vertices of $G$, then $I(S)\subseteq I_2(S)$.
\end{pro}

\noindent{\bf Algorithm 7.1} {\upshape \cite{OellermannPuertas}} (Finding $I(S)$ for a set $S\subseteq V(G)$ where $G$ is a distance-hereditary graph and $|S|\geq 2$).
\begin{itemize}
\item[] $(1)$ Let $H=<I_2(S)>$. Since $H$ is connected and induced, $H$ is distance hereditary.

\item[] $(2)$ Mark every vertex of $H$ that belong to $S$ with $T$ and all other vertices with $F$. (We associate with each vertex $v$ of $H$
a set $s(v)$ of vertices of the original graph $H$. These sets will be used to construct the Steiner interval for $S$.) For
each $v\in V(H)$ initialize $s(v)\longleftarrow \{v\}$ and let $I(S)\longleftarrow \emptyset$.

\item[] $(3)$ If $|S|=2$, go to Step $7$.

\item[] $(4)$ While $|S|\geq 3$ and $H$ contains a pendant $u$, proceed as follows: if $u$ is marked $T$, let $I(S)\longleftarrow I(S)\cup s(u)$, mark the
neighbor $u'$ of $u$ with $T$, let $H\longleftarrow H-u$ and $S\longleftarrow (S-u)\cup \{u'\}$; otherwise, if $u$ is marked $F$ let $H\longleftarrow H-u$.
Now if $|S|\geq 3$ (but $H$ contains no pendants), go to Step $5$; otherwise, go to Step $7$ (since $|S|=2$).

\item[] $(5)$ While $|S|\geq 3$ and $H$ has a pair $u,v$ of twins which have both been marked $T$ (or $F$, respectively) delete the vertex with
the larger label (say $v$) from $H$ and $S$ and modify $s(u)$, i.e., $H\longleftarrow H-v$, $S \longleftarrow S-\{v\}$ and let $s(u)\longleftarrow s(u)\cup s(v)$.
Now if $|S|\geq 3$ (but no such twins remain), proceed to Step $6$; otherwise, go to Step $7$ (since $|S|=2$).

\item[] $(6)$ While $|S|\geq 3$ and $H$ contains a pair of twins, one marked $T$ (say $u$) and the other $F$ (say $v$), then delete the vertex
marked $F$, i.e., define $H\longleftarrow H-v$ and return to Step $4$.

\item[] $(7)$ If $|S|=2$, say $S=\{x,y\}$, let $I(S)\longleftarrow I(S)\cup s(u)_{u\in I_H[x,y]}$. Output $I(S)$ and stop.
\end{itemize}

\begin{thm}{\upshape \cite{OellermannPuertas}}\label{th7-8}
If $G$ is a distance-hereditary graph and $S\subseteq V(G)$, $|S|\geq 2$, then the set $I(S)$ output by the algorithm is
the Steiner interval for $S$.
\end{thm}

\begin{rem}{\upshape \cite{OellermannPuertas}}\label{rem7-1}
Step $1$ of the Algorithm is not essential but is useful and allows $I(S)$ to be found more efficiently if $I(S)$
represents only a small proportion of the vertices of the graph.
\end{rem}

\section{Steiner Distance Stable Graphs}

In this section, we summarize the known results
on Steiner distance stable graphs and independent Steiner distance stable graphs.

\subsection{Steiner distance stable graphs}

The next result, due to Goddard, Oellermann, and Swart \cite{GoddardOellermannSwart}, shows that if distances of $(s,m)$-sets in a connected graph are
preserved after the deletion of certain numbers of vertices and edges, then so are
distances preserved for $(s,d)$-sets where $d>m$.
\begin{thm}{\upshape \cite{GoddardOellermannSwart}}\label{th8-1}
If a connected graph $G$ is $k$-vertex $\ell$-edge $(s,m)$-Steiner distance stable, then
it is $k$-vertex $\ell$-edge $(s,m+1)$-Steiner distance stable.
\end{thm}

\begin{cor}{\upshape \cite{GoddardOellermannSwart}}\label{cor8-1}
If a connected graph is $k$-vertex $\ell$-edge $(s,m)$-Steiner distance stable, then it
is $k$-vertex $\ell$-edge $(s,n)$-Steiner distance stable for all $n\geq m$.
\end{cor}

The next theorem implies another result of this type.
\begin{thm}{\upshape \cite{GoddardOellermannSwart}}\label{th8-2}
If a connected graph $G$ is $k$-vertex $\ell$-edge $(s,m)$-Steiner distance stable,
$m\geq s\geq 3$, then $G$ is $k$-vertex $\ell$-edge $(s-1,m)$-Steiner distance stable.
\end{thm}

\begin{cor}{\upshape \cite{GoddardOellermannSwart}}\label{cor8-2}
If a connected graph is $k$-vertex $\ell$-edge $(s,m)$-Steiner distance stable,
$m\geq s\geq 3$, then it is $k$-vertex $\ell$-edge $(s',m)$-Steiner distance stable for all $s' \ (2\leq s'\leq s)$.
\end{cor}

In Theorem \ref{th8-1}, one can see that the condition that a connected graph is $k$-vertex $\ell$-edge $(s,m)$-Steiner distance stable is sufficient for the graph to be $k$-vertex $\ell$-edge $(s,m+1)$-Steiner distance stable. The next result shows that this condition in not necessary.
\begin{thm}{\upshape \cite{GoddardOellermannSwart}}\label{th8-3}
For any integers $k$, $s$ and $m$ such that $s\geq 2$, $2s-2\geq m\geq s$ and $k\geq 1$, there
exists a graph $G$ which is $k$-vertex $0$-edge $(s,m+1)$-Steiner distance stable, but not
$k$-vertex $0$-edge $(s,m)$-Steiner distance stable.
\end{thm}

If we let $m=s-1$ in the construction of the proof of Theorem \ref{th8-3}, we obtain a graph
that is $k$-vertex $0$-edge $(s,s)$-Steiner distance stable and not $k$-vertex $0$-edge
$(s-1,s-1)$-Steiner distance stable.

They showed that the converse of Theorem \ref{th8-2} does not hold.
\begin{thm}{\upshape \cite{GoddardOellermannSwart}}\label{th8-4}
For $s\geq 3$ there is a graph which is $1$-vertex $0$-edge $(s-1,s-1)$-Steiner
distance stable but not $1$-vertex $0$-edge $(s,s)$-Steiner distance stable.
\end{thm}

Since the graph $G$ of the proof of Theorem \ref{th8-4} is $1$-vertex $0$-edge $(s-1, s-1)$-Steiner
distance stable, it follows, by Theorem \ref{th8-1}, that $G$ is $1$-vertex $0$-edge $(s-1,s)$-Steiner
distance stable. Since $G$ is not $1$-vertex $0$-edge $(s,s)$-Steiner distance stable, it follows
that the converse of Theorem \ref{th8-2} does not hold in general.

Recall it was shown in \cite{EntringerJacksonSlater}, for a positive integer $k$, that a graph is $(k,0,\{2\})$-stable if
and only if it is $(0,k,\{2\})$-stable. So a graph is $k$-vertex $0$-edge $(2,2)$-Steiner distance
stable if and only if it is $0$-vertex $k$-edge $(2,2)$-Steiner distance stable. The next result
shows that the necessity of this condition has an extension to $(3,3)$-sets.
\begin{thm}{\upshape \cite{GoddardOellermannSwart}}\label{th8-5}
For a positive integer $k$, a graph $G$ is $k$-vertex $0$-edge $(3,3)$-Steiner distance
stable if it is $0$-vertex $k$-edge $(3,3)$-Steiner distance stable.
\end{thm}

The converse of Theorem \ref{th8-5} does not hold; see \cite{GoddardOellermannSwart}. Let $H_1\cong K_s-uv \ (s\geq 3)$ for some pair $u,v$
of vertices and $H_2\cong K_2$ where $V(H_2)=\{x,y\}$. Let $G$ be obtained from $H_1\cup H_2$ by
adding the edges $ux$ and $vy$. Then $G$ is $1$-vertex $0$-edge $(s,s)$-Steiner distance stable, but
G is not $0$-vertex $1$-edge $(s,s)$-Steiner distance stable. To see this let $z\in V(H_1)-\{u,v\}$.
Then $d_G(\{x,y,z\})= 3$ but $d_{G-xy}(\{x,y,z\})=4$.

Theorem \ref{th8-5} cannot be extended to $(s,s)$-sets for $s\geq 4$; see \cite{GoddardOellermannSwart}. To see this, let
$G\cong (K_2\cup K_{s-2}+K_1$, i.e., $G$ is obtained by joining a new vertex to every vertex in
disjoint copies of $K_2$ and $K_{s-2}$. Then $G$ is $0$-vertex $1$-edge $(s,s)$-Steiner distance stable
but $G$ is not $1$-vertex $0$-edge $(s,s)$-Steiner distance stable.

\subsection{Independent Steiner distance stable graphs}

Goddard, Oellermann, and Swart \cite{GoddardOellermannSwart} also focused their attention on independent sets of vertices of a graph. Their
first result shows that in a certain sense the problem of finding Steiner trees for sets of independent vertices is equivalent to the problem of finding the Steiner trees of sets of
vertices that are not necessarily independent.

Let $\prod_1$ be the problem of finding
a Steiner tree for a nonempty set of vertices of a connected graph and $\prod_2$ the problem
of finding a Steiner tree for a nonempty independent set of vertices of a connected
graph. Let $G$ be a connected graph and $S$ a nonempty set of vertices of $G$.
Suppose $G_1,G_2,\ldots, G_n$ are the
components of $<S>_G$. Let $R(G; S)$ be the graph with vertex set
$(V(G)-S)\cup \{v_1,v_2,\ldots,v_n\}$ (where $v_i$ corresponds to $G_i$, $1\leq i\leq n$) and edge set
$\{uv\,|\,uv\in E(G-S)\}\cup \{uv_i\,|\,u\in V(G-S)$ and $u$ is adjacent in $G$ to some vertex of $G_i\}$.
Thus $R(G;S)$ is the contraction of $G$ that results from the partition
$(\bigcup_{i=1}^nV(G_i))\cup \{u\,|\,u\in V(G-S)\}$.
\begin{thm}{\upshape \cite{GoddardOellermannSwart2}}\label{th8-6}
There is an (ejficient) algorithm that solves $\prod_1$, is and only if there is an
(efficient) algorithm that solves $\prod_2$.
\end{thm}

Theorem \ref{th8-6} also follows directly from the nearest vertex reduction test described by
Beasley \cite{Beasley}.

The concepts presented in Subsection 9.2 and Theorem \ref{th8-6} suggest the next topic. If $G$ is
a connected graph and $S$ an independent set of $s$ vertices of $G$ such that $d_G(S)=m$, then $S$ is called an \emph{$I(s,m)$-set}. A connected graph is defined to be \emph{$k$-vertex $\ell$-edge $I(s,m)$-Steiner distance stable} if, for every $I(s,m)$-set $S$ and every set $A$ of at most $k$ vertices of
$G-S$ and at most $\ell$ edges of $G$, $d_{G-A}(S)=m$.

The following result in \cite{GoddardOellermannSwart} establishes an analogue of
Theorem \ref{th8-1} with respect to $I(3,m)$-sets.
\begin{thm}{\upshape \cite{GoddardOellermannSwart}}\label{th8-7}
If $G$ is a $k$-vertex $\ell$-edge $I(3,m)$-Steiner distance stable graph $m\geq 4$, then
$G$ is a $k$-vertex $\ell$-edge $I(3,m+1)$-Steiner distance stable graph.
\end{thm}

It remains an open problem to determine if a $k$-vertex $\ell$-edge $I(s, m)$-Steiner    distance stable graph $m\geq 4$, is a $k$-vertex $\ell$-edge $I(s,m+1)$-Steiner distance      stable graph, where $s>4$.

\section{Extremal problems on Steiner diameter}

What is the minimal size of a graph of order $n$ and diameter $d$?
What is the maximal size of a graph of order $n$ and diameter $d$?
It is not surprising that these questions can be answered without
the slightest effort (see \cite{Bollobas}) just as the similar
questions concerning the connectivity or the chromatic number of a
graph. The class of maximal graphs of order $n$ and diameter $d$ is
easy to describe and reduce every question concerning maximal graphs
to a not necessarily easy question about binomial coefficient, as in
\cite{HomenkoOstroverhii, HomenkoStrok, Ore, Watkins}. Therefore, the authors studied the minimal
size of a graph of order $n$ and under various additional
conditions.

Erd\"{o}s and R\'{e}nyi \cite{ErdosRenyi} introduced the following problem.
Let $d,\ell$ and $n$ be natural numbers, $d<n$ and $\ell<n$. Denote
by $\mathscr{H}(n,\ell,d)$ the set of all graphs of order $n$ with
maximum degree $\ell$ and diameter at most $d$. Put
$$
e(n,\ell,d)=\min\{e(G):G\in \mathscr{H}(n,\ell,d)\}.
$$
If $\mathscr{H}(n,\ell,d)$ is empty, then, following the usual
convention, we shall write $e(n,\ell,d)=\infty$. For more details on
this problem, we refer to \cite{Bollobas, Bollobas2, ErdosRenyi, ErdosRenyiSos}.

Mao \cite{Mao3} considered the generalization of the above problem. Let
$d,\ell$ and $n$ be natural numbers, $d<n$ and $\ell<n$. Denote by
$\mathscr{H}_k(n,\ell,d)$ the set of all graphs of order $n$ with
maximum degree $\ell$ and $sdiam_k(G)\leq d$. Put
$$
e_k(n,\ell,d)=\min\{e(G):G\in \mathscr{H}_k(n,\ell,d)\}.
$$
If $\mathscr{H}_k(n,\ell,d)$ is empty, then, following the usual
convention, we shall write $e_k(n,\ell,d)=\infty$. From Theorem
\ref{th3-5}, we have $k-1\leq d\leq n-1$.

The following results can be easily proved.
\begin{pro}{\upshape \cite{Mao3, MaoWang}}\label{pro9-1}
$(1)$ For $2\leq \ell \leq n-1$ and $3\leq k \leq n$,
$e_k(n,\ell,n-1)=n-1$.

$(2)$ For three integers $n,d,\ell$ with $2\leq d\leq n-2$ and
$n-d+2\leq \ell \leq n-2$,
$e_3(n,\ell,d)=n-1$.
\end{pro}

\subsection{Results for small $k$}

If $sdiam_3(G)=2$, then
$n-2\leq \delta(G)\leq n-1$,
and hence $n-2\leq \Delta(G)\leq n-1$. So one can assume that $n-2\leq
\ell \leq n-1$ for $d=2$.
\begin{thm}{\upshape \cite{Mao3}}\label{th9-1}
$(1)$ For $\ell=n-1$, $e_3(n,\ell,2)={n\choose 2}-\frac{n-1}{2}$ for
$n$ odd; $e_3(n,\ell,2)={n\choose 2}-\frac{n-2}{2}$ for $n$ even.

$(2)$ For $\ell=n-2$, $e_3(n,\ell,2)={n\choose 2}-\frac{n}{2}$ for
$n$ even; $e_3(n,\ell,2)=\infty$ for $n$ odd.
\end{thm}

In \cite{Mao3}, Mao got the following results for $d=3$.
\begin{thm}{\upshape \cite{Mao3}}\label{th9-2}
$(1)$ For $\ell=n-1$, $e_3(n,n-1,3)=n-1$;

$(2)$ For $\ell=n-2$, $e_3(n,n-2,3)=2n-5$;

$(3)$ For $\ell=n-3$, $e_3(n,n-3,3)=2n-5$;

$(4)$ For $\ell=2$, $e_3(n,2,3)=3$ for $n=4$; $e_3(n,2,3)=5$ for
$n=5$; $e_3(n,2,3)=\infty$ for $n\geq 6$.

$(5)$ For $\frac{n}{2}\leq \ell\leq n-4$, $n\leq e_3(n,\ell,3)\leq
\ell(n-\ell)$.
\end{thm}

For $d=n-2,n-3,n-4$, Mao obtained the following results.
\begin{thm}{\upshape \cite{Mao3}}\label{th9-3}
$(1)$ For $n\geq 4$, $e_3(n,2,n-2)=n$.

$(2)$ For $n\geq 4$,
$$
e_3(n,3,n-2)=\left\{
\begin{array}{ll}
n+1&\mbox{{\rm if}~$n=4$,}\\
n&\mbox{{\rm if}~$n=5$,}\\
n-1&\mbox{{\rm if}~$n\geq 6$.}
\end{array}
\right.
$$

$(3)$ For $n\geq 5$ and $4\leq \ell\leq n-1$, $e_3(n,\ell,n-2)=n-1$.
\end{thm}

\begin{thm}{\upshape \cite{Mao3}}\label{th9-4}
$(1)$ For $n\geq 5$,
$$
e_3(n,2,n-3)=\left\{
\begin{array}{ll}
\infty&\mbox{{\rm if}~$n=5,6$,}\\
n&\mbox{{\rm if}~$n\geq 7$.}
\end{array}
\right.
$$

$(2)$ For $n\geq 5$,
$$
e_3(n,3,n-3)=\left\{
\begin{array}{ll}
\infty,&\mbox{{\rm if}~$n=5$,}\\
n+1&\mbox{{\rm if}~$n=6$,}\\
n&\mbox{{\rm if}~$n=7$,}\\
n-1&\mbox{{\rm if}~$n\geq 8$}.
\end{array}
\right.
$$

$(3)$ For $n\geq 5$,
$$
e_3(n,4,n-3)=\left\{
\begin{array}{ll}
{n \choose 2}-2&\mbox{{\rm if}~$n=5$,}\\
n+1&\mbox{{\rm if}~$n=6$,}\\
n-1&\mbox{{\rm if}~$n\geq 7$}.
\end{array}
\right.
$$

$(4)$ For $n\geq 6$ and $5\leq \ell\leq n-1$, $e_3(n,\ell,n-3)=n-1$.
\end{thm}

\begin{thm}{\upshape \cite{Mao3}}\label{th9-5}
$(1)$ For $n\geq 5$,
$$
e_3(n,2,n-4)=\left\{
\begin{array}{ll}
\infty&\mbox{{\rm if}~$5\leq
n\leq 9$,}\\
n&\mbox{{\rm if}~$n\geq 10$.}
\end{array}
\right.
$$

$(2)$ For $n\geq 6$,
$$
e_3(n,3,n-4)=\left\{
\begin{array}{ll}
\infty&\mbox{{\rm if}~$n=6$,}\\
n+3&\mbox{{\rm if}~$n=7$,}\\
n+2&\mbox{{\rm if}~$n=8$,}\\
n+1&\mbox{{\rm if}~$n=9$,}\\
n-1&\mbox{{\rm if}~$n\geq 10$}.
\end{array}
\right.
$$

$(3)$ For $n\geq 6$,
$$
e_3(n,4,n-4)=\left\{
\begin{array}{ll}
2n&\mbox{{\rm if}~$n=6$,}\\
n+2&\mbox{{\rm if}~$n=7$,}\\
n-1&\mbox{{\rm if}~$n\geq 8$}.
\end{array}
\right.
$$

$(4)$ For $n\geq 6$,
$$
e_3(n,5,n-4)=\left\{
\begin{array}{ll}
2n+1&\mbox{{\rm if}~$n=6$,}\\
n+2&\mbox{{\rm if}~$n=7$,}\\
n-1&\mbox{{\rm if}~$n\geq 8$}.
\end{array}
\right.
$$

$(5)$ For $n\geq 7$ and $6\leq \ell\leq n-1$, $e_3(n,\ell,n-4)=n-1$.
\end{thm}

Mao \cite{Mao3} also constructed a graph and gave an upper bound of $e_3(n,\ell,d)$
for general $\ell$ and $d$.
\begin{pro}\label{pro9-2}
For $4\leq d\leq n-1$ and $2\leq \ell\leq n-1$,
$$
e_3(n,\ell,d)\leq \frac{(n-d+1)(n-d+2)}{2}+d-3.
$$
\end{pro}

\subsection{Results for large $k$}

The following result is immediate.
\begin{pro}{\upshape \cite{MaoWang}}\label{pro9-3}
For $2\leq \ell \leq n$, $e_n(n,\ell,n-1)=n-1$.
\end{pro}

For $k=n-1$ and $k=n-2$, Mao and Wang \cite{MaoWang} derived the
following results.
\begin{pro}{\upshape \cite{MaoWang}}\label{pro9-4}
$(1)$ For $2\leq \ell \leq n-1$, $e_{n-1}(n,\ell,n-1)=n-1$.

$(2)$ For $2\leq \ell \leq n-1$, $e_{n-1}(n,\ell,n-2)=n+\ell-2$.
\end{pro}

\begin{pro}{\upshape \cite{MaoWang}}\label{pro9-5}
For $2\leq \ell \leq n-1$ and $n\geq 5$,
$$
e_{n-2}(n,\ell,n-2)=\left\{
\begin{array}{ll}
n,&\mbox{{\rm if}~$2\leq \ell \leq n-2$;}\\
n-1,&\mbox{{\rm if}~$\ell=n-1$.}
\end{array}
\right.
$$
\end{pro}

Let $P_j^i$ be a path of order $j$, where $1\leq i\leq r+2$. We call
the graph $K_1\vee (K_1\cup P_j^i)$ as a {\it
$(u_i,v_i,P_j^i)$-Fan}; see Figure 9.1 $(a)$. For $1\leq i\leq r$, we
choose $(u_i,v_i,P_2^i)$-Fan, and choose
$(u_{r+1},v_{r+1},P_{\ell-1}^{r+1})$-Fan and
$(u_{r+2},v_{r+2},P_{s}^{r+2})$-Fan. Let $H_n$ be a graph obtained
from the above $(r+2)$ Fans by adding the edges in
\begin{eqnarray*}
&&\{w_1^{r+1}w_{s}^{r+2},w_1^{r+2}w_{1}^{1}\}\cup \{w_2^{i}w_{1}^{i+1}\,|\,1\leq i\leq r-1\}\cup \{w_2^{r}w_{\ell-1}^{r+1}\}\\[0.1cm]
&&\cup \{v_{i}v_{i+1}\,|\,1\leq i\leq r+1\}\cup \{v_{r+2}v_{1}\};
\end{eqnarray*}
see Figure 9.1 $(b)$, where $4r+\ell+s+3=n$, $2\leq s\leq 5$ and $1\leq r\leq \frac{n-\ell-4}{4}$.
\begin{figure}[!hbpt]
\begin{center}
\includegraphics[scale=0.5]{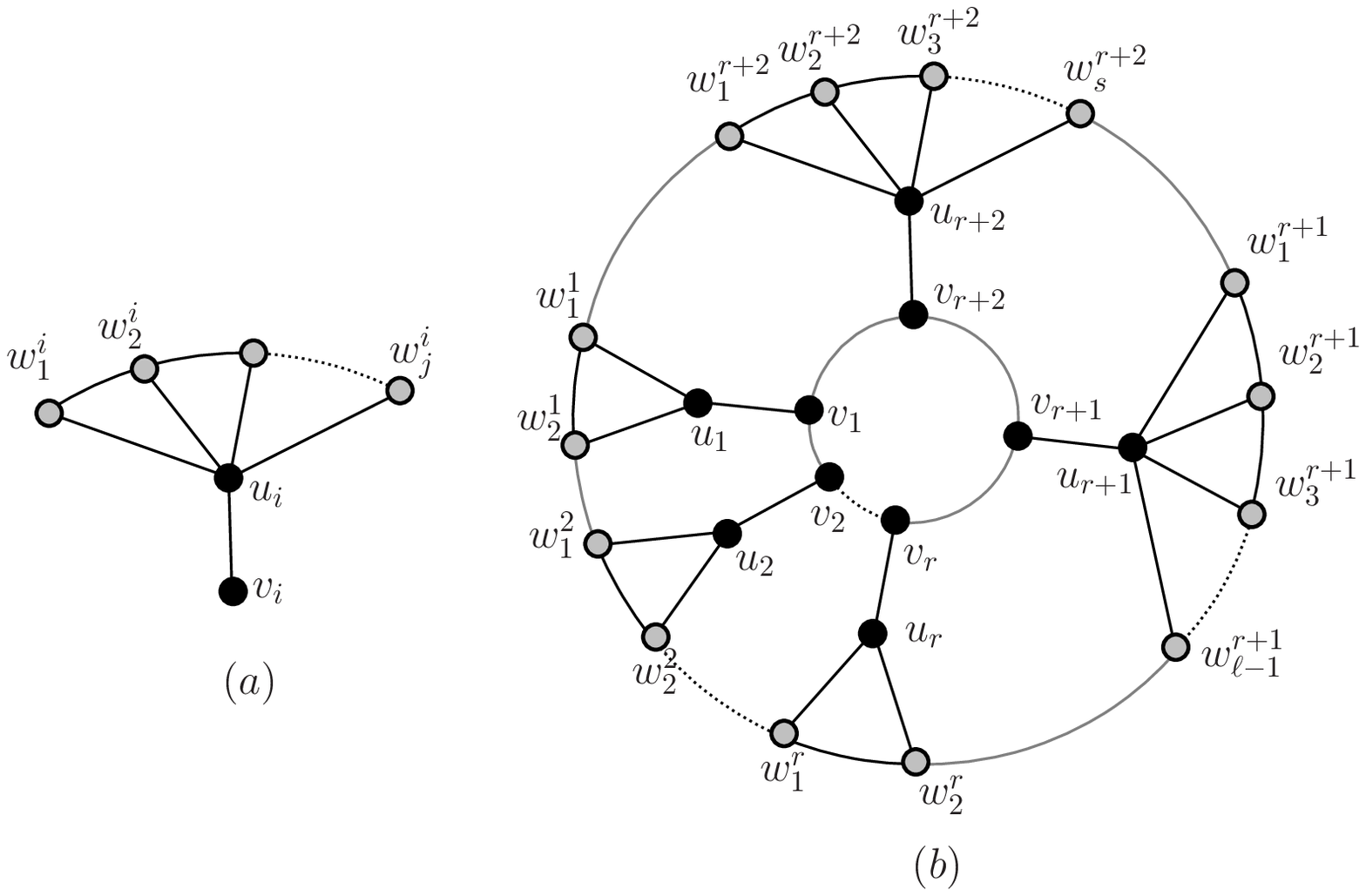}\\
Figure 9.1: Graphs for $(3)$
of Theorem \ref{th9-6}.
\end{center}\label{fig7}
\end{figure}

By the above graph class, Mao and Wang \cite{MaoWang} derived the
result in $(3)$ of Theorem \ref{th9-6} for $6\leq \ell \leq n-9$.
\begin{thm}{\upshape \cite{MaoWang}}\label{th9-6}
$(1)$ For $2\leq \ell \leq n-1$,
$$
e_{n-2}(n,\ell,n-1)=n-1.
$$

$(2)$ For $2\leq \ell \leq n-1$ and $n\geq 5$,
$$
e_{n-2}(n,\ell,n-2)=\left\{
\begin{array}{ll}
n,&\mbox{{\rm if}~$2\leq \ell \leq n-2$}\\
n-1,&\mbox{{\rm if}~$\ell=n-1$.}
\end{array}
\right.
$$

$(3)$ For $n-8\leq \ell \leq n-1$, $e_{n-2}(n,n-1-i,n-3)=2n-2$ for
$n\geq 5+i$ and $i=0,1$; $e_{n-2}(n,n-3-i,n-3)=2n-3$ for $n\geq
7+2i$ and $i=0,1$; $e_{n-2}(n,n-5-i,n-3)=2n-4$ for $n\geq 11+2i$ and
$i=0,1$; $e_{n-2}(n,n-7-i,n-3)=2n-5$ for $n\geq 15+2i$ and $i=0,1$.
For $6\leq \ell \leq n-9$,
$$
\frac{1}{2}(3n+\ell-3)\leq e_{n-2}(n,\ell,n-3)\leq
\frac{1}{2}(3n+\ell+s-5),
$$
where $2\leq s\leq 5$. Furthermore, if $s=2$, then
$e_{n-2}(n,\ell,n-3)=\frac{1}{2}(3n+\ell-3)$.
\end{thm}

Let $A_{32}$ be a minimally $4$-connected graph shown in Figure 9.2
$(a)$ (see \cite{Bollobas}, Page 18). We now give a graph $H_n$ of
order $n \ (n\geq 96)$ such that $\Delta(H_n)=\ell$ and
$sdiam_{n-3}(H_n)=n-4$ constructed by the following steps.
\begin{itemize}
\item[] \textbf{Step 1}: For each $i \ (1\leq i\leq x)$, we let $A_{32}^i$ be the copy of
$A_{32}$, where $n=32x+y$, $x=\lfloor n/32\rfloor$, and $0\leq y\leq
31$. Let $V(A_{32}^i)=\{u^i_{j}\,|\,1\leq j\leq 12\}\cup
\{v^i_{j}\,|\,1\leq j\leq 20\}$ such that $d_G(u^i_{j})=5$ for
$1\leq j\leq 12$, and $d_G(v^i_{j})=4$ for $1\leq j\leq 20$; see
Figure 9.2 $(a)$. Let $B_{32x}$ be a graph obtained from $A_{32}^i \
(1\leq i\leq x)$ by adding the edges in
$\{v_5^{i}v_1^{i+1}\,|\,1\leq i\leq x-1\}\cup
\{v_6^{i}v_2^{i+1}\,|\,1\leq i\leq x-1\}\cup
\{v_7^{i}v_3^{i+1}\,|\,1\leq i\leq x-1\}\cup
\{v_8^{i}v_4^{i+1}\,|\,1\leq i\leq x-1\}$; see Figure 9.2 $(b)$.
\begin{figure}[!hbpt]
\begin{center}
\includegraphics[scale=0.7]{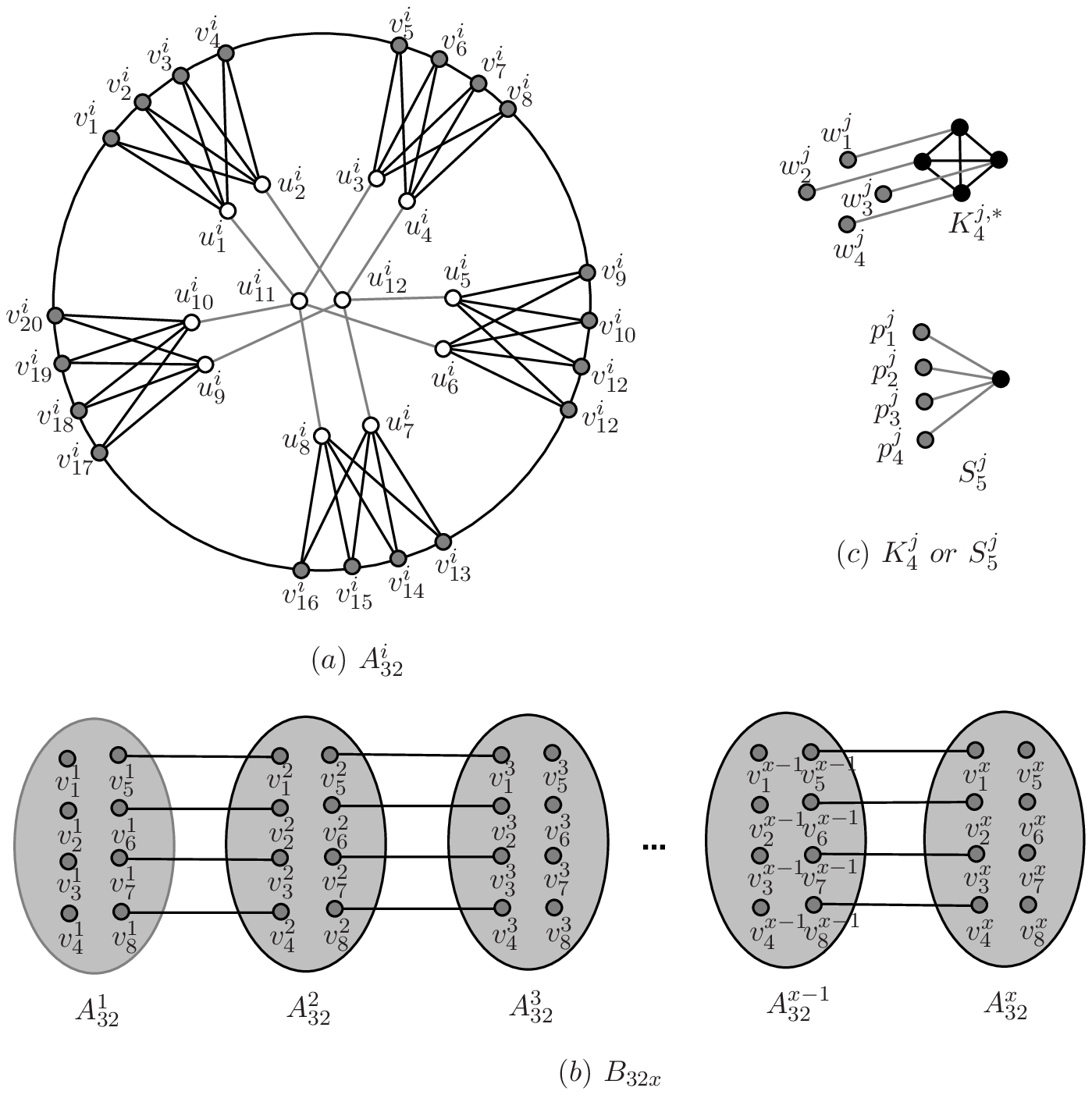}\\[0.5cm]
Figure 9.2: Graphs for $(3)$
of Theorem \ref{th9-7}.
\end{center}\label{fig7}
\end{figure}

\item[] \textbf{Step 2}: Let $y=4z+a$, where $z=\lfloor
y/4\rfloor$, $0\leq a\leq 3$. For each $j \ (1\leq j\leq z)$, we let
$K_4^{j}$ be the complete graph of order $4$. Furthermore, let
$K_4^{j,*}$ be the graph obtained from $K_4^j$ by adding four
pendant vertices $w_1^j,w_2^j,w_3^j,w_4^j$ with four pendant edges
such that another end vertex of each pendant edge is attached on
only one vertex in $K_4^j$; see Figure 9.2 $(c)$. For each $j \ (1\leq
j\leq a)$, we let $S_5^{j}$ be the star of order $5$ with its leaves
$p_1^j,p_2^j,p_3^j,p_4^j$. Since $n\geq 96$, it follows that
$A_{32}^1,A_{32}^2,A_{32}^3$ all exist. Set $S_1=\{v_j^{1}\,|\,1\leq
j\leq 4\}\cup \{v_j^{1}\,|\,9\leq j\leq 20\}\subseteq V(A_{32}^1)$,
and $S_2=\{v_j^{2}\,|\,9\leq j\leq 20\}\subseteq V(A_{32}^2)$, and
$S_3=\{v_j^{3}\,|\,9\leq j\leq 20\}\subseteq V(A_{32}^3)$. Then
$|S_1\cup S_2\cup S_3|=40$. If $n\equiv 0~(mod~32)$, then
$D_n=B_{32x}$. If $n\neq 0~(mod~32)$ and $n-32x\equiv 0~(mod~4)$,
then $D_n$ is a graph obtained from $B_{32x}$ and
$K_4^{1,*},K_4^{2,*},\ldots,K_4^{z,*}$ by identifying each vertex in
$S'=\{w_j^{i}\,|\,1\leq i\leq 4, \ 1\leq j\leq z\}$ and only one
vertex in $S_1\cup S_2\cup S_3$. Since $|S'|=4z<40=|S_1\cup S_2\cup
S_3|$, for any vertex in $S'$, we can find a vertex in $S_1\cup
S_2\cup S_3$ and then identify the two vertices. If $n\neq
0~(mod~32)$ and $n-32x\neq 0~(mod~4)$, then $D_n$ is a graph
obtained from $B_{32x}$, $K_4^{1,*},K_4^{2,*},\ldots,K_4^{z,*}$ and
$S_5^{1},S_5^{2},\ldots,S_5^{a}$ by identifying each vertex in
$S'=\{w_j^{i}\,|\,1\leq i\leq 4, \ 1\leq j\leq z\}\cup
\{p_j^{i}\,|\,1\leq i\leq 4, \ 1\leq j\leq a\}$ and only one vertex
in $S_1\cup S_2\cup S_3$. Since $|S'|=4z+4a\leq 28+12=40=|S_1\cup
S_2\cup S_3|$, for any vertex in $S'$, we can find a vertex in
$S_1\cup S_2\cup S_3$ and then identify the two vertices.

\item[] \textbf{Step 3}: Let $H_n$ be the graph $D_n$ by adding
$\ell-5$ edges between $u_{12}^{1}$ and $V(G)-u_{12}^{1}$.
\end{itemize}

By the above graph class, they derived the result in $(3)$
of Theorem \ref{th9-7}.
\begin{thm}{\upshape \cite{MaoWang}}\label{th9-7}
$(1)$ For $2\leq \ell \leq n-1$, $e_{n-1}(n,\ell,n-1)=n-1$.

$(2)$ For $2\leq \ell \leq n-1$ and $n\geq 4$,
$$
e_{n-3}(n,\ell,n-2)=\left\{
\begin{array}{ll}
n,&\mbox{{\rm if}~$2\leq \ell \leq \lfloor\frac{n}{2}\rfloor-1$};\\
&\mbox{{\rm ~~~or}~$\ell=\lfloor\frac{n}{2}\rfloor$~and~n~is~odd};\\
n-1,&\mbox{{\rm if}~$\lfloor\frac{n}{2}\rfloor+1\leq \ell \leq n-1$};\\
&\mbox{{\rm ~~~or}~$\ell=\lfloor\frac{n}{2}\rfloor$~and~n~is~even}.\\
\end{array}
\right.
$$

$(3)$ For $2\leq \ell \leq n-1$ and $n\geq 96$,
$$
2n-2-\lceil \ell/2\rceil\leq e_{n-3}(n,\ell,n-4)\leq 74\left\lfloor
\frac{n}{32}\right\rfloor+2i+\ell-9,
$$
where $n\equiv~i~(mod~32)$, $1\leq i\leq 31$.

$(4)$ If $2\leq \ell \leq n-1$, then $e_{n-3}(n,\ell,n-3)\geq
\max\{n-1+\lceil \ell/2\rceil,\frac{3n-\ell-5}{2}\}$. If $\lceil
n/2\rceil+1\leq \ell\leq n-1$, then $e_{n-3}(n,\ell,n-3)\leq
2n-\ell+1$. If $5\leq \ell\leq \lceil n/2\rceil$, then
$$
e_{n-3}(n,\ell,n-3)\leq (2\ell+3)\left\lfloor
\frac{n}{\ell+1}\right\rfloor+\ell+\left\{
\begin{array}{ll}
-8,&\mbox{{\rm if}~$n\equiv~0~(mod~\ell+1)$};\\
-5,&\mbox{{\rm if}~$n\equiv~1~(mod~\ell+1)$};\\
-2,&\mbox{{\rm if}~$n\equiv~2~(mod~\ell+1)$};\\
&\mbox{{\rm ~~~or}~$n\equiv~3~(mod~\ell+1)$};\\
2i-7,&\mbox{{\rm if}~$n\equiv~i~(mod~\ell+1)$},
\end{array}
\right.
$$
where $n=(\ell+1)x+i$ and $0\leq i\leq \ell$.
\end{thm}

Mao and Wang \cite{MaoWang} constructed a graph and gave an upper bound of $e_k(n,\ell,d)$
for general $k,\ell,d$.
\begin{thm}{\upshape \cite{MaoWang}}\label{th9-8}
Let $k,\ell,d$ be three integers with $2\leq k\leq n$, $2\leq
\ell\leq n-1$, and $k-1\leq d\leq n-1$.

$(1)$ If $d=k-1$, $\lceil \frac{n+1}{2}\rceil\leq k\leq n$, and
$\max\{n-k+1,\lceil\frac{n}{2}\rceil\}<\ell\leq n-1$, then
$$
\left\lceil \frac{\ell+(n-1)(n-k+1)}{2}\right\rceil \leq
 e_k(n,\ell,d)\leq \frac{(n-1)^2}{4}+\ell.
$$

$(2)$ If $2\leq k\leq d$, $k\leq d\leq n-1$, and $2+\lceil
\frac{n-d+k-3}{d-k+1}\rceil\leq \ell\leq n-1$, then
$$
e_k(n,\ell,d)=n-1.
$$
\end{thm}


\begin{thebibliography}{1}


\bibitem{AkiyamaHarary}
J. Akiyama and F. Harary, A graph and its complement with specified
properties, \emph{Internat. J. Math. \& Math. Sci.} 2(2) (1979),
223--228.

\bibitem{Ali}
P. Ali, The Steiner diameter of a graph with prescribed girth,
\emph{Discrete Math.} 313(12) (2013), 1322--1326.

\bibitem{AliBoalsSherwani}
H.H. Ali, A. Boals, and N.S. Sherwani, \emph{Distance stable
graphs}, Paper read at 2nd Internat. Conf. in Graph Theory,
Combinatorics, Algorithms and Applications, San Francisco, 1989.

\bibitem{AliDankelmannMukwembi}
P. Ali, P. Dankelmann, and S. Mukwembi, Upper bounds on the Steiner
diameter of a graph, \emph{Discrete Appl. Math.} 160 (2012),
1845--1850.


\bibitem{AliMukwembiDankelmann}
P. Ali, S. Mukwembi, and P. Dankelmann, Steiner diameter of $3$, $4$
and $5$-connected maximal planar graphs, \emph{Discrete Appl. Math.}
179 (2014), 222--228.




\bibitem{AnandChangatKlavzarPeterin}
B.S. Anand, M. Changat, S. Klav\v{z}ar, I. Peterin, Convex sets in
lexicographic product of graphs, \emph{Graphs Combin.} 28 (2012),
77--84.

\bibitem{Anthonisse}
J. M. Anthonisse, \emph{The rush in a directed graph}, Technical Report BN
9/71, Stiching Math. Centrum, Amsterdam, October 1971.

\bibitem{AouchicheHansen}
M. Aouchiche and P. Hansen, A survey of Nordhaus-Gaddum type
relations, \emph{Discrete Appl. Math.} 161 (2013), 466--546.

\bibitem{Atici}
M. Atici, Computational complexity of geodetic set, \emph{Internat.
J. Comput. Math.} 79 (2002), 587--591.

\bibitem{Avann}
S.P. Avann, Metric ternary distributive semi-lattices, \emph{Proc. Amer. Math. Sot.}
12 (l961), 407--414.

\bibitem{BandeltMulder}
H. J. Bandelt and H. M. Mulder, Distance-hereditary graphs, \emph{J.
Combin. Theory B} 41 (1986), 183--208.

\bibitem{Brandstadt}
A. Brandst\"{a}dt, V.B. Le, and J.P. Spinrad, \emph{Graph Classes: A Survey}, SIAM Monograph on Discrete Mathematics and Applications, Philidelphia,
1999.

\bibitem{BaoIgarashiOhring}
F. Bao, Y. Igarashi, and S.R. \"{O}hring, Reliable broadcasting in
product networks, \emph{Discrete Applied Math.} 83 (1998), 3--20.


\bibitem{Beasley}
J.E. Beasley, An SST-based algorithm for the Steiner problem in
graphs, \emph{Networks} 19 (1989) 1--16.

\bibitem{BeinekeOellermannPippert}
L.W. Beineke, O.R. Oellermann, R.E. Pippert, On the Steiner median
of a tree, \emph{Discrete Appl. Math.} 68 (1996), 249--258.

\bibitem{Bloom}
G.S. Bloom, A characterization of graphs of diameter two,
\emph{Amer. Math. Monthly} 95(1) (1988), 37--38.

\bibitem{Bollobas}
B. Bollob\'{a}s, Extremal Graph Theory,
Acdemic press, 1978.

\bibitem{Bollobas2}
B. Bollob\'{a}s, \emph{Graphs with a given diameter and maximal
valency and with a minimal number of edges}, in:``Combinatorial
Mathematics and its Applications'' (Welsh, D.J.A., ed.), Academic
Press, London and New York, 1971), 25--37.

\bibitem{Bondy}
J.A. Bondy and U.S.R. Murty,
{\it Graph Theory}, GTM 244, Springer, 2008.

\bibitem{BuckleyHarary}
F. Buckley and F. Harary, \emph{Distance in Graphs},
Addision-Wesley, Redwood City, CA (1990).

\bibitem{BuckleyMillerSlater}
F. Buckley, Z. Miller, and P.J. Slater, On graphs containing a given
graph as center, \emph{J. Graph Theory} 5 (1981), 427--434.


\bibitem{CaccettaSmyth}
L. Caccetta, W.F. Smyth, Properties of edge-maximal $K$-edge-connected $D$-critical graphs, \emph{J. Combin. Math. Comput. Comput.} 2 (1987), 111--131.

\bibitem{CaccettaSmyth2}
L. Caccetta, W.F. Smyth, $K$-edge-connected $D$-critical graphs of minimum order, \emph{Congr. Numer.} 58 (1987), 225--232.

\bibitem{CaccettaSmyth3}
L. Caccetta, W.F. Smyth, Restriction of vertices for maximum edge count in $K$-edge-connected $D$-critical graphs, \emph{Ars Combin} 26B (1988), 115--132.

\bibitem{CaceresMarquezPuertas}
J. Caceres, A. Marquez, and M.L. Puertas, Steiner distance and
convexity in graphs, \emph{European J. Combin.} 29 (2008), 726--736.

\bibitem{Calkin}
N.J. Calkin and H.S. Wilf, The number of independent sets in a grid
graph, \emph{SIAM J. Discrete Math.} 11(1) (1998), 54--60.

\bibitem{ChartrandLesniak}
G. Chartrand and L. Lesniak, \emph{Graphs \& Digraphs}, 2nd Edition,
Monterey, CA, 1986.


\bibitem{ChartrandOellermannTianZou}
G. Chartrand, O.R. Oellermann, S. Tian, and H.B. Zou, Steiner
distance in graphs, \emph{\'{C}asopis pro p\v{e}stov\'{a}n\'{i}
matematiky} 114 (1989), 399--410.

\bibitem{ChartrandOkamotoZhang}
G. Chartrand, F. Okamoto, and P. Zhang, Rainbow trees in graphs and
generalized connectivity, \emph{Networks} 55 (2010), 360--367.

\bibitem{ChartrandZhang}
G. Chartrand, P. Zhang, The Steiner number of a graph, \emph{Discrete Math.} 242 (2002), 41--54.

\bibitem{Cho}
J. Cho, \emph{A min-cost flow based min-cost rectilinear Steiner distance-preserving tree construction},
Proceedings of the 1997 international symposium on Physical design 82--87.


\bibitem{Chung}
F.R.K. Chung, \emph{Diameter of graphs: Old problems and new
results}, 18th Southeastern Conf. on Combinatorics, Graph Theory and
Computing (1987).


\bibitem{DankelmannEntringer}
P. Dankelmann and R. Entringer, Average distance, minimum degree,
and spanning trees, \emph{J. Graph Theory} 33 (2000), 1--13.


\bibitem{DankelmannOellermannSwart}
P. Dankelmann, O.R. Oellermann, and H.C. Swart, The average Steiner
distance of a graph, \emph{J. Graph Theory} 22(1) (1996), 15--22.

\bibitem{DankelmannSwartOellermann}
P. Dankelmann, H.C. Swart, and O.R. Oellermann, On the average
Steiner distance of graphs with prescribed properties,
\emph{Discrete Appl. Math.} 79 (1997), 91--103.



\bibitem{DankelmannSwartOellermann2}
P. Dankelmann, H. Swart, and O.R. Oellermann, (Kalamazoo, MI, 1996)
\emph{Bounds on the Steiner diameter of a graph}, Combinatorics,
Graph Theory, and Algorithms, Vol. I, II, New Issues Press,
Kalamazoo, MI, 1999, 269--279.


\bibitem{DAtriMoscarini}
A. D'Atri and M. Moscarini, Distance-hereditary graphs, Steiner
trees, and connected domination, \emph{SIAM J. Comput.} 17(3),
521--538.



\bibitem{DayOellermannSwart}
D.P. Day, O.R. Oellermann, and H.C. Swart, Steiner
distance-hereditary graphs, \emph{SIAM Discrete Math.} 7(3) (1994),
437--442.

\bibitem{DayOellermannSwart2}
D.P. Day, O.R. Oellermann, and H.C. Swart, A characterization of
$3$-Steiner distance hereditary graphs, \emph{Networks} 30(4)
(1997), 243--253.

\bibitem{DoyleGraver}
J.K. Doyle and J.E. Graver, Mean distance in a graph, \emph{Discrete
Math.} 17 (1977), 147--154.

\bibitem{DobryninKochetova} A. Dobrynin, A. Kochetova, Degree distance of a graph:
        A degree analogue of the Wiener index,
        {\it J. Chem. Inf. Comput. Sci.\/} {\bf 34} (1994) 1082--1086.

\bibitem{DuLyuuHsu}
D.Z. Du, Y.D. Lyuu, and D.F. Hsu, Line digraph iteration and
connectivity analysis of de Bruijn and Kautz graphs, \emph{IEEE
Trans. Comput.} 42 (1993), 612--616.


\bibitem{EntringerJacksonSlater}
R.C. Entringer, D.E. Jackson, and P.J. Slater, Geodetic connectivity
of graphs, \emph{IEEE Trans. Circuits Sys.} 24 (1977), 460--463.

\bibitem{EntringerJacksonSnyder}
R.C. Entringer, D.E. Jackson, and P.A. Snyder, Mean distance in a
graph, \emph{Czech. Math. J.} 26 (1976), 283--296.



\bibitem{ErdosPachPollackTuza}
P. Erd\"{o}s, J. Pach, R. Pollack, Z. Tuza, Radius, diameter, and
minimum degree, \emph{J. Combin. Theory B} 47(1989), 73--79.

\bibitem{ErdosRenyi}
P. Erd\"{o}s and A. R\'{e}nyi, On a problem in the theory of graphs
(in Hungarian), \emph{Publ. Math. Inst. Hungar. Acad. Sci.} 7(1962).

\bibitem{ErdosRenyiSos}
P. Erd\"{o}s, A. R\'{e}nyi, and V.T. S\'{o}s, On a problem of graph
theory, \emph{Studia Sci. Math. Hungar.} 1 (1966), 215--235.

\bibitem{ErohOellermann}
L. Eroh, O.R. Oellermann, Geodetic and Steiner geodetic sets in
$3$-Steiner distance hereditary graphs, \emph{Discrete Math.} 308
(2008), 4212--4220.

\bibitem{FavaronKouiderMaheo}
O. Favaron, M. Kouider, M. Mah\'{e}o, Edge-vulnerability and mean distance,
\emph{Networks} 19 (1989), 493--504.

\bibitem{FarberJamison}
M. Farber, R.E. Jamison, Convexity in graphs and hypergraphs,
\emph{SIAM J. Alg. Disc. Math.} 7(3) (1986), 433--444.

\bibitem{FinkLuzarSkrekovski} J. Fink, B. Lu\v{z}ar, R. \v{S}krekovski,
        \emph{Some remarks on inverse Wiener index problem},
        Discrete Appl. Math. 160 (2012) 1851--1858.


\bibitem{Freeman}
L.C. Freeman, A set of measures of centrality based upon betweenness,
\emph{Sociometry} 40 (1977), 35--41.

\bibitem{FulekMoricPritchard}
R. Fulek, F. Mori\'{c}, D. Pritchard,
\emph{Diameter bounds for planar graphs},
Discrete Math. 311(2011), 327--335.

\bibitem{FurtulaGutmanKatanic} B. Furtula, I. Gutman, V. Katani\'{c},
        Three-center Harary index and its applications,
        \emph{Iranian J. Math. Chem.} 7(1) (2016), 61--68.


\bibitem{GareyJohnson}
M.R. Garey and D.S. Johnson, \emph{Computers and Intractibility: A Guide to
the Theory of $NP$-Completeness}, Freeman \& Company, New York, 1979.

\bibitem{GirvanNewman}
M. Girvan, M.E.J. Newman, Community structure in social and biological networks, \emph{Proc. Natl. Acad. Sci. USA} 99 (2002), 7821--7826.

\bibitem{Goddard}
W. Goddard, A note on Steiner-distance-hereditary graphs, \emph{J.
Combin. Math. Combin. Comput.} 40 (2002), 167--170.



\bibitem{GoddardOellermann}
W. Goddard, O. R. Oellermann, Distance in graphs, in:
        M. Dehmer (Ed.), {\it Structural Analysis of Complex Networks\/},
        Birkh\"{a}user, Dordrecht, 2011, pp. 49--72.

\bibitem{GoddardOellermannSwart}
W. Goddard, O.R. Oellermann, H.C. Swart, Steiner distance stable
graphs, \emph{Discrete Math.} 132(1994), 65--73.

\bibitem{GoddardOellermannSwart2}
W. Goddard, O.R. Oellermann and H.C. Swart, A new approach to
distance stable graphs, \emph{J. Combin. Math. Combin. Comput.} 8
(1990), 209--220.


\bibitem{Gologranc}
T. Gologranc, Steiner convex sets and Cartesian product, \emph{Bull.
Malays. Math. Sci. Soc.}, DOI 10.1007/s40840-016-0312-8.

\bibitem{GutmanSDD} I. Gutman,
On Steiner degree distance of trees, \emph{Appl. Math. Comput.} 283
(2016), 163--167.

\bibitem{GutmanFurtulaLi} I. Gutman, B. Furtula, X. Li, Multicenter Wiener indices
        and their applications, J. Serb. Chem. Soc. 80 (2015) 1009--1017.

\bibitem{GutmanYeh} I. Gutman, Y.N. Yeh,
        \emph{The sum of all distances in bipartite graphs},
        Math. Slovaca 45 (1995) 327--334.

\bibitem{GutmanYehChen} I. Gutman, Y.N. Yeh, J.C. Chen,
        \emph{On the sum of all distances in graphs},
        Tamkang J. Math. 25 (1986) 83--86.

\bibitem{Gyori}
E. Gy\"{o}ri, \emph{On division of graphs to connected graphs},
Combinatorics, North-Holland, New York, 1978, 485--494.

\bibitem{Hakimi}
S.L. Hakimi, Steiner's problem in graph and its implications,
\emph{Networks} 1 (1971), 113--133.


\bibitem{Hammack}
R. Hammack, W. Imrich, and S. Klav\u{z}ar,
\emph{Handbook of product graphs}, Secend edition, CRC Press, 2011.

\bibitem{HammerMaffray}
P.L. Hammer and F. Maffray, Completely separable graphs,
\emph{Discrete Appl. Math.} 27 (1990), 85--100.

\bibitem{BandeltMulder}
H.-J. Bandelt, H.M. Mulder, Distance-hereditary graphs, \emph{J. Combin. Theory Ser. B} 41 (1986), 182--208.

\bibitem{HararyRobinson}
F. Harary, R.W. Robinson, The diameter of a graph and its complement, \emph{Amer. Math. Monthly} 92 (1985), 211--212.

\bibitem{HellSeyffarth}
P. Hell, K. Seyffarth, Largest planar graphs of diameter two and
fixed maximum degree, \emph{Discrete Math.} 111 (1993), 313--322.

\bibitem{HenningOellermannSwart}
M.A. Henning, O.R. Oellermann, H.C. Swart, On the Steiner radius and
Steiner diameter of a graph, \emph{Ars Combin.} 29 (1990), 13--19.

\bibitem{Hendry}
G.R.T. Hendry, On graphs with prescribed median $I$, J.
Graph Theory 9 (1985), 477--481.

\bibitem{HenningOellemannSwart2}
M.A. Henning, O.R. Oellemann, and H.C. Swart, \emph{On vertices with
maximum Steiner eccentricity in graphs}, Graph Theory,
Combinatorics, Algorithms and Applications (eds. Y. Alavi, F. R. K.
Chung, R. L. Graham, and D. F. Hsu.), SIAM Publications,
Philadelphia (1991), 393--403.

\bibitem{Holbert}
K.S. Holbert, A note on graphs with distant center and median,
Recent studies in graph theory, V.R. Kulli (Editor),
Vishwa, 1989.

\bibitem{HomenkoOstroverhii}
N.P. Homenko, N.A. Ostroverhii, Diameter-critical graphs (in
Russian), \emph{Ukrainian Math. J.} 22(1970), 637--646.


\bibitem{Howorka}
E. Howorka, A characterization of distance hereditary graphs,
\emph{Q. J. Math. Oxford} 28 (1977), 417--420.




\bibitem{Hsu}
D.F. Hsu, \emph{On container width and length in graphs, groups, and
networks}, IEICE Transaction on Fundamentals of Electronics,
Communications and Computer Science, E77-A (1994), 668--680.

\bibitem{Hsu2}
D.F. Hsu and T. {\L}uczak, Note on the $k$-diameter of $k$-regular
$k$-connected graphs, \emph{Discrete Math.} 133(1994), 291--296.

\bibitem{HwangRW}
F.K. Hwang, D.S. Richards, and P. Winter, \emph{The Steiner
Tree Problem}, North-Holland, Amsterdam, 1992.


\bibitem{Itai}
A. Itai and M. Rodeh, The multi-tree approach to reliability in
distributed networks, \emph{Inform. Comput.} 79 (1988), 43--59.

\bibitem{Johnsson}
S.L. Johnsson and C.T. Ho, Optimum broadcasting and personaized
communication in hypercubes, \emph{IEEE Trans. Comput.} 38(9)(1989),
1249--1268.

\bibitem{Jordan}
C. Jordan, Sur les assemblages de lignes, \emph{J. Reine. Angew.
Math.} 70 (1869), 185--190.

\bibitem{Kovse}
M. Kov\v{s}e, Vertex decomposition of Steiner Wiener index and Steiner betweenness
centrality, arXiv:1605.00260 [math.CO] 2016.

\bibitem{Ku}
S. Ku, B. Wang, and T. Hung, Constructing edge-disjoint spanning
trees in product networks, \emph{IEEE Trans. Parall. Distr. Sys.}
14(3) (2003), 213--221.

\bibitem{KubickaKubickiOellermann}
E. Kubicka, G. Kubicki, O.R. Oellermann, Steiner intervals in graphs,
\emph{Discrete Appl. Math.} 81 (1998), 181--190.

\bibitem{LepovicGutman} M. Lepovi\'{c}, I. Gutman,
      A collective property of trees and chemical trees,
      J. Chem. Inf. Comput. Sci. 38(1998) 823--826.

\bibitem{Levi}
A.Y. Levi, Algorithm for shortest connection of a group of graph
vertices, \emph{Sov. Math. Dokl.} 12(1971), 1477--1481.

\bibitem{LiMaoGutman}
X. Li, Y. Mao, I. Gutman, The Steiner Wiener index of a graph,
\emph{Discuss. Math. Graph Theory} 36(2)(2016)  455--465.

\bibitem{LiMaoGutman2}
X. Li, Y. Mao, I. Gutman, Inverse problem on the Steiner Wiener
index, \emph{Discuss. Math. Graph Theory}, in press.

\bibitem{Lovasz}
L. Lov\'{a}sz, A homology theory for spanning trees of a graph,
\emph{Acta Math. Acad. Sci. Hungar.} 30 (1977) 241--251.

\bibitem{Lovasz2}
L. Lov\'{a}sz, \emph{Combinatorial Problems and Exercises},
Akademiai Kiado, Budapest, 1979.

\bibitem{Mao}
Y. Mao, Path-connectivity of lexicographical product  graphs,
\emph{Int. J. Comput. Math.} 93(1)(2016), 27--39.

\bibitem{Mao2}
Y. Mao, The Steiner diameter of a graph, \emph{Bull. Iran. Math.
Soc.} 43(2)(2017), 439--454.


\bibitem{Mao3}
Y. Mao, Steiner $3$-diameter, maximum degree and size of a graph,
arXiv:1703.04974 [math.CO] 2017.

\bibitem{Mao4}
Y. Mao, Steiner Harary index, \emph{Kragujevac J. Math.}
42(1)(2018), 29--39.

\bibitem{MaoDas}
Y. Mao, K.C. Das, Steiner Gutman index, {\it MATCH Commun. Math.
Comput. Chem.}, in press.



\bibitem{MaoMelekianCheng}
Y. Mao, C. Melekian, E. Cheng, \emph{The Steiner $(n-3)$-diameter of
a graph}, arXiv:1703.03984 [math.CO] 2017.

\bibitem{MaoChengWang} Y. Mao, E. Cheng, Z. Wang, \emph{Steiner distance in product networks},
        arXiv:1703.01410 [math.CO] 2017.

\bibitem{MaoWang}
Y. Mao, Z. Wang, Steiner diameter, maximum degree and size of a graph,
 arXiv:1704.04695 [math.CO] 2017.

\bibitem{MaoWangWangWang} Y. Mao, Z. Wang, X. Wang, C. Wang,
Nordhaus-Guddum type results for the Steiner Harary of graphs,
\emph{Iranian J. Math. Chem.} 8(2) (2017), 181--198.

\bibitem{MaoWangGutman}
Y. Mao, Z. Wang, and I. Gutman, Steiner Wiener index of graph
products, \emph{Trans. Combin.} 5(3) (2016), 39--50.

\bibitem{MaoWangGutmanKlobucar}
Y. Mao, Z. Wang, I. Gutman, and A. Klobu\v{c}ar, Steiner degree
distance, \emph{MATCH Commun. Math. Comput. Chem.} 78(1) (2017),
221--230.

\bibitem{MaoWangGutmanLi}
Y. Mao, Z. Wang, I. Gutman, H. Li, Nordhaus-Gaddum-type results for
the Steiner Wiener index of graphs, \emph{Discrete Appl. Math.}
219(2017), 167--175.

\bibitem{MaoWangXiaoYe}
Y. Mao, Z. Wang, Y. Xiao, C. Ye, Steiner Wiener index and
connectivity of graphs, \emph{Utilitas Math.} 102 (2017), 51--57.

\bibitem{MeyerPradhan}
F.J. Meyer, and D.K. Pradhan, Flip  trees, \emph{IEEE Trans.
Comput.} 37(3) (1987), 472--478.

\bibitem{Miller}
Z. Miller, Medians and distance sequences of graphs, Ars
Comb. 15 (1983), 177--179.

\bibitem{Mulder}
H.M. Mulder, \emph{The interval function of a graph},
Mathematical Centre Tracts 132, Amsterdam, 1980.

\bibitem{Nebesky}
L. Neb\'{e}sky, Median graphs, \emph{Comment. Math. Univ. Carolinae} 12 (1971), 317--325.

\bibitem{NewmanGirvan}
M.E.J. Newman, M. Girvan, Finding and evaluating community structure in networks, \emph{Phys. Rev. E} 69 (2004), 026113.

\bibitem{NielsenOellermann}
M.H. Nielsen, O.R. Oellermann, Steiner trees and convex geometries,
\emph{SIAM J. Discrete Math.} 23(2)(2009), 680--693.

\bibitem{Oellermann}
O.R. Oellermann, From steiner centers to steiner medians, \emph{J.
Graph Theory} 20(2) (1995), 113--122.

\bibitem{Oellermann2}
O.R. Oellermann, \emph{On Steiner centers and Steiner
medians of graphs}, Networks 34(4)(1999), 258--263.

\bibitem{OellermannPuertas}
O.R. Oellermann, M.L. Puertas, Steiner intervals and Steiner
geodetic numbers in distance-hereditary graphs, \emph{Discrete
Math.} 307 (2007), 88--96.


\bibitem{OellermannSpinrad}
O. Oellermann, and J.P. Spinrad, A polynomial algorithm for testing
whether a graph is $3$-Steiner distance hereditary, \emph{Inform.
Proc. Lett.} 55(3) (1995), 149--154.

\bibitem{OellermannTian}
O.R. Oellermann and S. Tian, Steiner centers in graphs, \emph{J.
Graph Theory} 14(5)(1990), 585--597.

\bibitem{Ore}
O. Ore, Diameter in graphs, \emph{J. Combin. Theory} 5 (1968),
75--81.


\bibitem{Pelayo}
I. Pelayo, Comment on ``The Steiner number of a graph'' by G.
Chartrand and P. Zhang [Discrete Mathematics 242 (2002) 41-54],
\emph{Discrete Math.} 280 (2004), 259--263.

\bibitem{Plesnik}
J. Plesn\'{\i}k, On the sum of all distances in a graph or digraph,
\emph{J. Graph Theory} 8 (1984), 1--24.


\bibitem{Rouvray1} D.H. Rouvray, Harry in the limelight: The life and
        times of Harry Wiener, in: D. H. Rouvray, R. B. King (Eds.),
        {\it Topology in Chemistry -- Discrete Mathematics of
        Molecules\/}, Horwood, Chichester, 2002, pp. 1--15.

\bibitem{Rouvray2} D.H. Rouvray, The rich legacy of half century of
        the Wiener index, in: D. H. Rouvray, R. B. King (Eds.),
        {\it Topology in Chemistry -- Discrete Mathematics of
        Molecules\/}, Horwood, Chichester, 2002, pp. 16--37.


\bibitem{Seyffarth}
K. Seyffarth, Maximal planar graphs of diameter two, \emph{J. Graph
Theory} 13 (1989), 619--648.

\bibitem{Straffin}
P.D. Straffin, Letter to editor, \emph{Amer. Math. Monthly} 93 (1986), 76.

\bibitem{Slater}
P.J. Slater, Medians of arbitrary graphs, J. Graph Theory 4 (1980), 389--392.

\bibitem{Soltes}
L. Soltes, Transmission in graphs: a bound on vertex removing,
\emph{Math Slov.} 41 (1991), 11--16.



\bibitem{TomescuMelter}
I. Tomescu, R.A. Melter, On distances in chromatic graphs, Quart.
\emph{J. Math. Oxford} 40(2)(1989), 475--480.

\bibitem{Van de Vel}
M.J.L. Van de Vel, \emph{Theory of Convex Structures}, North-Holland, Amsterdam, 1993.

\bibitem{Wagner} S. Wagner,
        A note on the inverse problem for the Wiener index,
        \emph{MATCH Commun. Math. Comput. Chem.} 64 (2010), 639--646.

\bibitem{Wagner2} S. Wagner,
        \emph{A class of trees and its Wiener index},
        Acta Appl. Math. 91, 119--132.

\bibitem{WagnerWangYu} S.G. Wagner, H. Wang, G. Yu,
        \emph{Molecular graphs and the inverse Wiener index problem},
        Discrete Appl. Math. 157 (2009) 1544--1554.


\bibitem{WangMaoDas}
Z. Wang, Y. Mao, K.C. Das, Nordhaus-Guddum type results for the Steiner degree distance of
graphs, submitted.


\bibitem{WangMaoDas2}
Z. Wang, Y. Mao, K.C. Das, Nordhaus-Guddum type results for the
Steiner Gutman index of graphs, submitted.

\bibitem{WangMaoLiYe}
Z. Wang, Y. Mao, H. Li, C. Ye, The Steiner $4$-diameter of a graph,
arXiv:1702.05681 [math.CO] 2017.

\bibitem{WangMaoChengMelekian}
Z. Wang, Y. Mao, E. Cheng, C. Melekian, Steiner distance in corona and cluster graphs,
arXiv:1702.05681 [math.CO] 2017.

\bibitem{Watkins}
M.E. Watkins,
\emph{A lower bound for the number of vertices of a graph},
Amer. Math. Monthly 74 (1976), 297.

\bibitem{Wiener} H. Wiener,
        \emph{Structural determination of paraffin boiling points},
        J. Am. Chem. Soc. 69 (1947), 17--20.

\bibitem{Winter}
P. Winter, \emph{Steiner problems in networks: a survey},
Networks 17 (1987), 129--167.

\bibitem{Xu}
S.J. Xu,
\emph{Some parameters of graph and its complement},
Networks 65(2)(1987), 197--207.



\bibitem{YehChiangPeng}
H.-G. Yeh, C.-Y. Chiang, S.-H. Peng, Steiner centers and Steiner medians of graphs,
\emph{Discrete Mathematics} 308 (2008), 5298--5307.

\bibitem{YehGutman}
Y. Yeh and I. Gutman, On the sum of all distances in composite
graphs, {\em Discrete Math.} {\bf 135} (1994), 359--365.

\bibitem{ZhangWu} L. Zhang, B. Wu, The Nordhaus--Gaddum-type inequalities
        for some chemical indices, \emph{MATCH Commun. Math. Comput. Chem.}
        54 (2005), 189--194.
\end{thebibliography}
\end{document}